\documentclass[10pt,  aps, pra, notitlepage, superscriptaddress, longbibliography, showpacs,nofootinbib, floatfix]{revtex4-1}

\usepackage{amsmath}
\usepackage{graphicx}
\usepackage{amsfonts}
\usepackage{amsbsy}
\usepackage{bbm}
\usepackage{dsfont}
\usepackage{color}
\usepackage{mathrsfs}
\usepackage{soul}
\usepackage{enumitem}

\newcommand{\argmax}{\mathrm{argmax}}
\newcommand{\argmin}{\mathrm{argmin}}

\newcommand{\Var}{\mathrm{Var}}
\newcommand{\Covar}{\mathrm{Cov}}
\newcommand{\Tran}{\mathrm{T}}

\newcommand{\x}{\mathbf{x}}

\newcommand{\z}{\mathbf{z}}

\newcommand{\cmatrix}{  \mathbf{C}  }
\newcommand{\J}{  \mathbf{J}  }
\newcommand{\C}{  \mathbf{C}  }
\newcommand{\Cov}{   \ensuremath{\boldsymbol{\Sigma}}  }
\newcommand{\A}{  \ensuremath{\mathbf{A}}   }

\newcommand{\I}{\mathbf{I}} 
\newcommand{\1}{\mathds{1}}
\newcommand{\thetav}{ \boldsymbol{\theta}  }
\newcommand{\thetaest}{  \ensuremath{\boldsymbol{\hat{\theta}}}   }
\newcommand{\Thetav}{  \ensuremath{\boldsymbol{\Theta}}   }

\newcommand{\vv}{\ensuremath{\mathbf{v}}}
\newcommand{\nullv}{\ensuremath{\mathbf{0}}}

\newcommand{\muv}{  \ensuremath{\boldsymbol{\mu}}   }

\newcommand{\tv}{\mathbf{t}}
\newcommand{\Zm}{ \ensuremath{ \mathbf{Z}}}
\newcommand{\epsilonv}{  \ensuremath{\boldsymbol{\epsilon}}   }
\newcommand{\rmd}{\mathrm{d}}
\newcommand{\rme}{\mathrm{e}}
\newcommand{\rmi}{\mathrm{i}}
\newcommand{\Tr}{\mathrm{Tr}}

\makeindex

\begin{document}

\title{Exact results on high-dimensional linear regression via statistical physics} 

%

\author{Alexander Mozeika}
\affiliation{London Institute for Mathematical Sciences, Royal Institution, 21 Albemarle Street, London W1S 4BS, UK. }

\author{Mansoor 	Sheikh}
\affiliation{Saddle Point Science Ltd., 10 Lincoln Street, York YO26 4YR, UK. }

\author{Fabian Aguirre-Lopez}
\affiliation{Universit\'{e} Paris-Saclay, CNRS, LPTMS, 91405, Orsay, France. }

\author{Fabrizio Antenucci}
\affiliation{Saddle Point Science Ltd., 10 Lincoln Street, York YO26 4YR, UK. }

\author{Anthony CC Coolen}
\affiliation{London Institute for Mathematical Sciences, Royal Institution, 21 Albemarle Street, London W1S 4BS, UK. }
\affiliation{Dept of Biophysics, Radboud University, 6525AJ Nijmegen, The Netherlands. }
\affiliation{Saddle Point Science Ltd., 10 Lincoln Street, York YO26 4YR, UK. }

\date{\today}

\begin{abstract}
It is clear that conventional statistical inference protocols need to be revised to deal correctly with the high-dimensional data that are now common.
Most recent studies aimed at achieving this revision rely  on powerful approximation techniques,  that call for rigorous results against which they can be tested. 
In this context,  the simplest case of high-dimensional linear regression has acquired significant new relevance and attention. 
In this paper we use the statistical physics perspective on inference to derive a number of new exact results for  linear regression in the high-dimensional regime.
\end{abstract}

\maketitle

\section{Introduction\label{section:intro}}
The advent of modern  high-dimensional data poses a significant challenge to statistical inference.
The latter is  understood well in the conventional regime of growing sample size with constant dimension.
For high-dimensional data, where the dimension  is of the same order as the sample size, the foundations of inference methods 
are still fragile, and even the simplest scenario of   
linear regression \cite{Rencher2008} has to be revised \cite{Huber2009}.   
The study of linear regression (LR)  in the high-dimensional  regime  has  recently attracted  significant attention 
in the mathematics \cite{Bayati2011,  Karoui2013, Karoui2016, Wainwright2019} and statistical 
physics communities \cite{Advani2016, Barbier2019, Coolen2020}.
The statistical physics framework is naturally suited to deal with high-dimensional data.

While the connection between statistical physics and information theory 
was established a while ago by Jaynes \cite{Jaynes1957},  
the approach has more recently been extended also to information processing \cite{Nishimori2001}  
and machine learning \cite{Bishop2006}.     
In the statistical physics framework, the free energy encodes  statistical properties of inference, akin to cumulant generating functions in  statistics, 
but its direct  computation via high-dimensional integrals is often difficult. This led to the development of  
 several non-rigorous methods, such as the mean-field approximation,   
the replica trick and the cavity method \cite{Mezard1987}. 
Message passing in particular, which can be seen as algorithmic implementation of the latter~\cite{Mezard2009}, 
has emerged as an efficient 
analysis tool for statistical inference in high dimensions~\cite{Donoho2016, Gerbelot2020, Barbier2020}.

\begin{figure*}
 \setlength{\unitlength}{0.89mm}
 \begin{picture}(200,90)
 
\put(0,0){\includegraphics[height=90\unitlength]{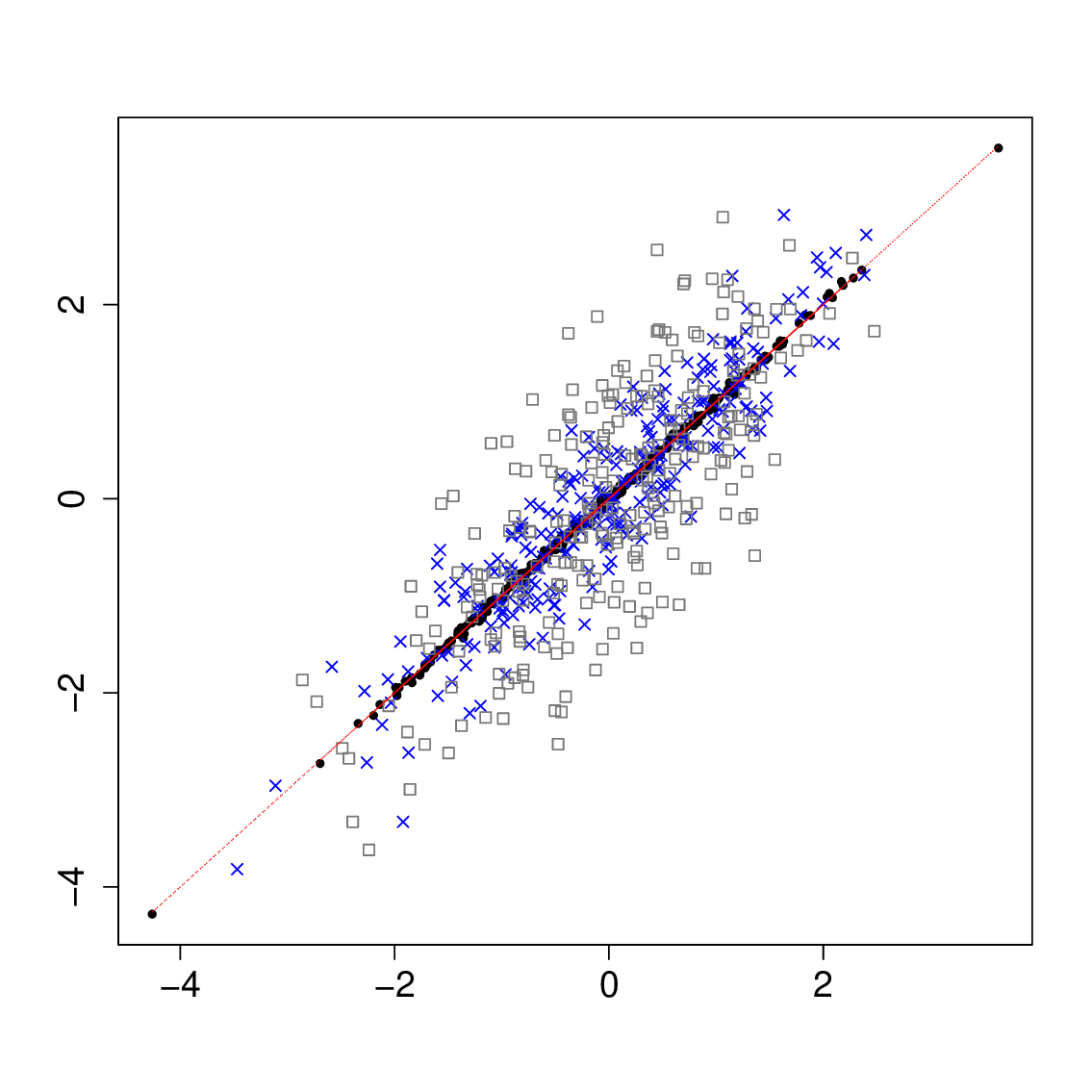}} 

\put(100,0){\includegraphics[height=90\unitlength]{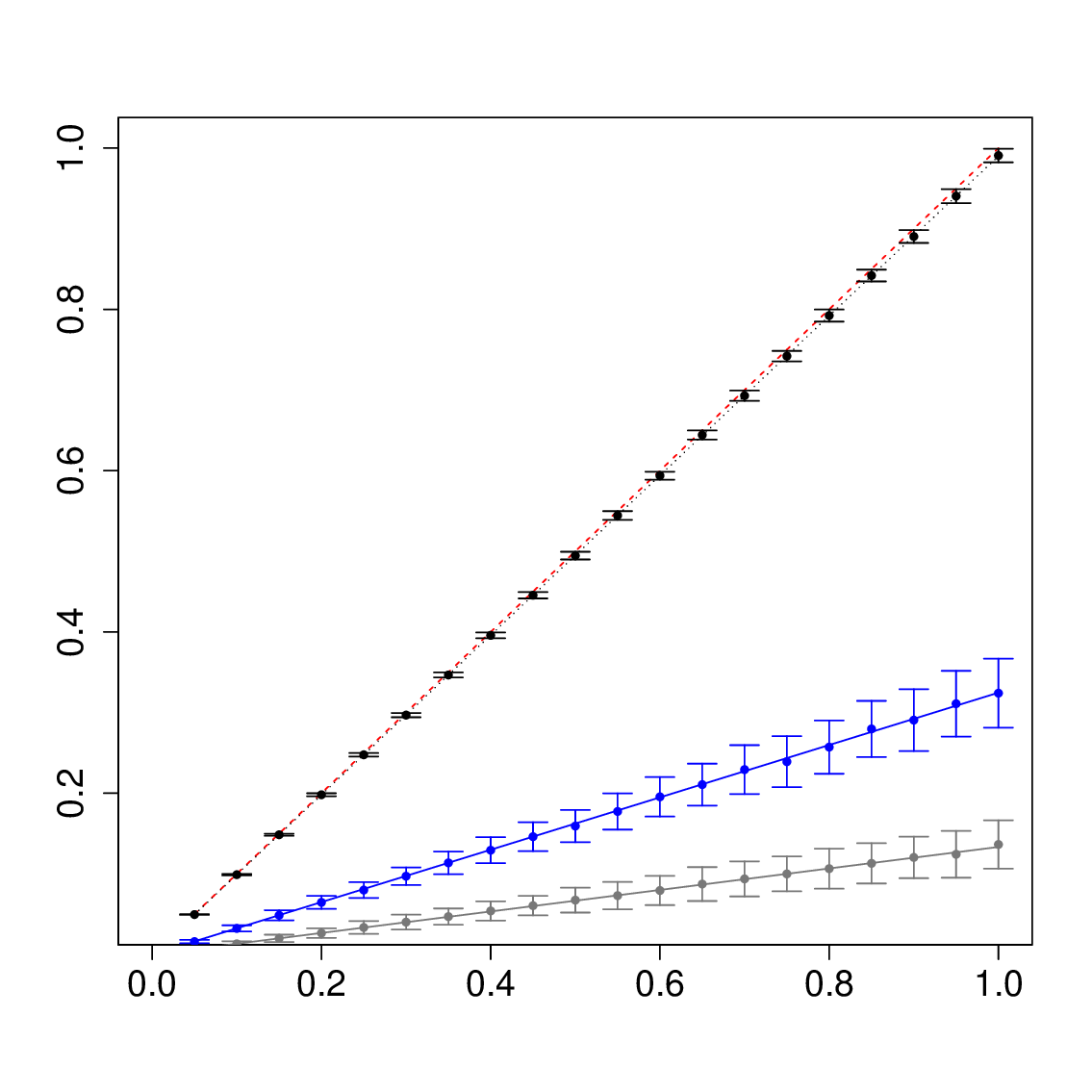}}

\put(-3,77){$(a)$} \put(97,77){$(b)$}

 
  \put(-1 ,44){ $\hat{\theta}$}  \put(43 ,2){$\theta_0$} 
    
  \put(99 ,44){$\hat{\sigma}^2$ }  \put(150 ,2){$\sigma_0^2$}       
\end{picture}
 \caption{High-dimensional phenomena in   inference  with  the linear regression (LR)  model in the teacher-student scenario (see text). 
Comparison between parameters inferred with maximum likelihood (ML) $(\hat{\theta},\hat{\sigma})$ and true values $(\theta_0,\sigma_0)$. 
 (a)  Plot of the ordered set $\hat{\theta}(\theta_0)=\{(\hat{\theta}(1),\theta_0(1)),\ldots,(\hat{\theta}(d),\theta_0(d))\}$ for $d/N \in \{0.01, 0.675, 0.867\}$, represented by  symbols $\{\bullet,\times, \square \}$,  with  $N \in \{26000, 385, 300\}$. For each value of $d/N$ the rows of $\Zm$ were sampled from $\mathcal{N}(\nullv, \I_d)$, $\epsilonv$  was sampled from $\mathcal{N}(\nullv, \sigma_0^2\I_N)$, with $\sigma_0^2=0.1$, and $\thetav_0$ was sampled from $\mathcal{N}(\nullv, \I_d)$.  (b)  Plot of $\hat{\sigma}^2$ versus $\sigma^2_0$,  represented by points connected by lines, for $d/N \in \{0.01, 0.675, 0.867\}$ (top to bottom). Each point, together with $\pm$ one standard-deviation error-bars, represents an average over $250$ samples.  Note that in both plots the diagonal line corresponds to perfect inference. }
 \label{figure:over-fitting-in-LR} 
 \end{figure*}

Most rigorous results on high-dimensional LR were obtained upon 
 assuming {uncorrelated} data \cite{Karoui2013, Donoho2016,  Barbier2019,  Barbier2020}, possibly with  
{sparsity} of parameters \cite{Bayati2011, Wainwright2019}. 
Recently,  correlations in {sampling} were analyzed in \cite{Gerbelot2020}
for rotationally invariant data matrices.
In all these studies, however, the parameters of the noise in  the data  were assumed \emph{known}, 
unlike the standard statistical setting where they are  
usually inferred \cite{Rencher2008}. 
The exact prescription of the noise strength is unwelcome, since it is artificially removing an important source of overfitting in 
realistic applications of regression. In high-dimensional LR, inference protocols can mistake noise for signal, reflected in increased under-estimation of the noise and 
over-estimation  of the magnitude of other model parameters (see  Figure \ref{figure:over-fitting-in-LR}).   

In this paper we derive new exact  
results for the high-dimensional regime of Bayesian LR which complement the aforementioned rigorous studies, using  
the  statistical physics  formulation of inference.

\textit{Statement of the problem and  preview of results}\textendash  
We consider Bayesian inference  of   the LR model,  $\tv=\Zm\thetav+\sigma\epsilonv$,
where $\tv\!\in\!{\rm I\!R}^N$ and $\Zm\!\in{\rm I\!R}^{N\times d}$ are observed 
and the parameters  $\thetav\!\in\!{\rm I\!R}^d$ and $\sigma\!\in\!{\rm I\!R}^+$ are to be inferred, 
with $\epsilonv$ denoting zero-average noise. 
We adopt a \textit{teacher-student} scenario \cite{engel2001statistical,ZdeborovaKrzakala2016review}:
the teacher  samples independently the rows of $\Zm$ from some probability distribution $P(\z)$
and then uses the  LR model to obtain $\tv$ with  the \emph{true} parameters $(\thetav_0, \sigma_0) $.
We assume that the student then applies the Bayes formula to try to  infer $(\thetav, \sigma)$
assuming a  Gaussian\footnote{The distribution  of   a  Gaussian (or Normal)  random variable  $\x\in\!{\rm I\!R}^d$, with mean $\muv\in\!{\rm I\!R}^d$ and covariance $\Cov\in{\rm I\!R}^{d\times d}$,   is given by the density  $\mathcal{N}(\x \vert \muv , \Cov)=\frac{\rme^{-\frac{1}{2} (\x-\muv)^\Tran\Cov^{-1} (\x-\muv) }}{\vert2\pi\Cov\vert^{d/2}}$.} prior $\mathcal{N}(\thetav\vert\nullv, \eta^{-1}\I_d)$ for $\thetav$, and a generic prior $P(\sigma^2)$ for 
the noise  parameter $\sigma^2$.
Specifically, we do not consider the case where the observations are coming from an unknown source 
and/or where one needs to do model selection.

We  map the LR inference problem
onto a Gibbs-Boltzmann  distribution  with inverse  `temperature' $\beta$. This allows us to investigate properties of different inference protocols. In particular, 
{maximum a posteriori} (MAP) inference
is obtained for $\beta\! \to\! \infty$ and $\eta\! >\! 0$, 
{maximum likelihood} (ML) inference
for $\beta\! \to\! \infty$ and $\eta\!=\!0$,
and marginalization inference  for $\beta\! =\! 1$.
We will refer to `ML (MAP) at finite temperature'
for the case  of $\eta\!=\!0$ ($\eta>0$) and $\beta$ finite.

The \emph{high-dimensional} regime is obtained for $(N,d)\!\rightarrow\!(\infty, \infty)$ with fixed ratio  $\zeta\!=\!d/N\!\in\!(0,\infty)$.
We will henceforth indicate this limit as $(N,d)\rightarrow \infty$, to simplify notation.   
Note that in order to keep $\tv$  finite in the $(N,d)\rightarrow \infty$ limit, 
the matrix $\Zm$ has to be  replaced with $\Zm/\sqrt{d}$ (unless of course we impose a suitable sparsity condition).

Within the above setting  we obtain the following results: (i) If $\sigma^2$ is known and the  distributions of $\Zm$ and $\epsilonv$ are  Gaussian, we compute the 
distribution  of the MAP and ML estimators of $\thetav$; (ii) The ML estimator $\hat{\sigma}^2_{\text{ML}}$ of the noise parameter $\sigma^2$ 
is \emph{self-averaging} as $(N,d)\!\rightarrow\!\infty$ (i.e. its variance  is vanishing\footnote{Here we adopt the  definition from statistical physics of disordered systems which states  that some density, such as  free energy, average energy, etc., is self-averaging if  its variance is vanishing in the thermodynamic limit~\cite{Mezard1987}. In statistics this phenomena is also referred to as having a `fully concentrated measure' \cite{vershynin2018high,Wainwright2019}. }  
 in this limit), 
for any  distributions of $\Zm$ and $\epsilonv$. We  bound the likelihood of deviations of  $\hat{\sigma}^2_{\text{ML}}$ 
from its mean for  Gaussian noise $\epsilonv$;  (iii) We compute the 
characteristic function of the \textit{mean square error}   $\frac{1}{d}\vert\vert\thetav_{\!0}-\thetaest_{\rm ML}[\mathscr{D}]\vert\vert^2$ for the ML estimator at  finite $(N,d)$, 
where $\thetav_0$ are the \emph{true}  parameters;  (iv)  We determine average and variance 
of  the  free energy, associated with  Gibbs-Boltzmann distribution of Bayesian LR,  of  ML inference for the finite     $\beta$ and  finite  $(N,d)$.  The ML free energy density is self-averaging as $(N,d)\!\rightarrow\!\infty$
if the eigenvalue \emph{spectrum} of the empirical covariance matrix $\Zm^\Tran\Zm/N$ is self-averaging. For Gaussian  $\epsilonv$ and   $\Zm$, we  recover 
the results obtained by the replica method in~\cite{Coolen2020};   (v) If the  true parameters $\thetav_0$ are 
independent \emph{random} variables, we derive average  and  variance  of the  free energy  of MAP inference  for finite  $\beta$ 
 and $(N,d)$.  The MAP free energy is shown to be self-averaging  if the spectrum of $\Zm^\Tran\Zm/N$ is self-averaging as $(N,d)\!\rightarrow\!\infty$. 

In the following subsections we describe how the above results were obtained, with full mathematical details relegated to  the Appendix. 

\section{Statistical Physics and Bayesian Inference}
 We assume that we  observe a data  sample   of  `input-output' pairs $\{(\z_1, t_1), \ldots, (\z_N,t_N)\}$, 
 where $(\z_i,t_i) \in{\rm I\!R}^{d+1}$, generated randomly and independently from 
\begin{eqnarray}
P(t, \z\vert\Thetav)= P(t\vert \z, \Thetav)  P(\z), \label{def:P(t, z)}
\end{eqnarray}
with parameters $\Thetav$ that are unknown to us.   If we assume a prior distribution $P(\Thetav)$,  
then the distribution of $\Thetav$, given the data, follows from the Bayes formula
\begin{eqnarray}
P(\Thetav\vert \mathscr{D} )&=&\frac{ P(\Thetav) \prod_{i=1}^NP(t_i\vert \z_i, \Thetav) }{ \int\!\rmd\tilde{\Thetav} ~P(\tilde{\thetav}) \left\{\prod_{i=1}^NP(t_i\vert \z_i, \tilde{\Thetav}) \right\} }\label{def:posterior}. 
\end{eqnarray}
Here $\mathscr{D} \!=\!\{\tv,\Zm\}$, with  $\tv\!=\!(t_1, \ldots, t_N)$, and $\Zm\!=\!(\z_1, \ldots,\z_N)$ is 
an $N\times d$ matrix.  In Bayesian language, expression (\ref{def:posterior}) is the {posterior} distribution of $\Thetav$, given the {prior} distribution   
$P(\Thetav)$ and the observed  data $\mathscr{D}$. 

The simplest way to use (\ref{def:posterior}) for inference  is  to compute  the \emph{maximum a posteriori} (MAP) estimator 
\begin{eqnarray}
\hat{\Thetav}_{\text{MAP}} \!\left[  \mathscr{D}  \right]
&=& \argmin_{\Thetav}\,E\!\left(\Thetav\vert\mathscr{D}\right) \label{def:MAP-estimator},
\end{eqnarray}
in which  the so-called Bayesian {likelihood}  function 
\begin{eqnarray}
E\left(\Thetav\vert \mathscr{D} \right)&=&-\sum_{i=1}^N\log P(t_i\vert \z_i, \Thetav)-\log P\!\left(\Thetav\right) \label{def:E}
\end{eqnarray}
consists of a first term, the log-likelihood used also in \emph{maximum likelihood} (ML) inference, 
and a second term, that acts as a {regularizer}.  
Bayesian  inference can thus be seen as a generalization of MAP inference, and MAP inference generalizes ML inference. 

The  \emph{square error} $\frac{1}{d}\vert\vert\Thetav_{\!0}-\hat{\Thetav}[\mathscr{D}]\vert\vert^2$, with the Euclidean norm $\vert\vert\cdots\vert\vert$ and 
the \emph{true} parameters $\Thetav_{0}$ underlying the data, is often used
to quantify the quality of inference in (\ref{def:MAP-estimator}).
Its first moment is the 
\emph{mean square error} (MSE)  $ \frac{1}{d}\langle\langle\vert\vert\Thetav_{\!0}-\hat{\Thetav}[\mathscr{D}]\vert\vert^2\rangle_{\mathscr{D}}\rangle_{\Thetav_0}$.  
Furthermore, the posterior mean 
\begin{eqnarray}
\hat{\Thetav}\!\left[  \mathscr{D}  \right]&=& \int\!\rmd\Thetav~ P(\Thetav\vert \mathscr{D}) \Thetav \label{eq:MMSE-estimator}
\end{eqnarray}
(the \textit{marginalization} estimator) is the \emph{minimum} MSE (MMSE) estimator in the
{Bayes optimal case}, i.e. when prior distribution and model likelihood are known \cite{ZdeborovaKrzakala2016review}.

%
The above approaches to Bayesian inference  can be unified conveniently in a \emph{\em single} statistical physics  (SP) formulation via the Gibbs-Boltzmann  distribution  
\begin{eqnarray}
P_{\beta}(\Thetav\vert\mathscr{D})&=&\frac{\rme^{-\beta E\left(\Thetav\vert\mathscr{D}\right)}}{Z_{\beta}[\mathscr{D}]}  \label{def:Gibbs},
\end{eqnarray}
with the normalization constant, or  `partition function' $Z_{\beta}[\mathscr{D}]=\int\!\rmd \Thetav~\rme^{-\beta E\left(\Thetav\vert\mathscr{D}\right)}$. For $\beta=1$ this is  the \emph{evidence} term  of Bayesian inference.  
In statistical physics language, (\ref{def:E})  plays the role of `energy'  in (\ref{def:Gibbs}) and $\beta$ is the  (fictional) inverse temperature.  
The temperature can be interpreted as a noise amplitude in stochastic gradient 
descent minimization of $E(\Thetav\vert\mathscr{D})$ \cite{Coolen2020}. 
Properties  of the  system   (\ref{def:Gibbs})  follow upon evaluating the  `free energy' 
\begin{eqnarray}
 F_{\beta}[\mathscr{D}] &=& -\frac{1}{\beta}\log Z_{\beta}[\mathscr{D}] \label{def:F}.
\end{eqnarray}

The estimators (\ref{def:MAP-estimator})   and  (\ref{eq:MMSE-estimator}) are recovered  from  the average $\int\!\rmd\Thetav~P_\beta(\Thetav\vert \mathscr{D}) \Thetav$ by taking the zero `temperature'  limit  $\beta\rightarrow\infty$, or by setting $\beta=1$,  respectively.  
This follows upon observing that for  $\beta=1$ the distribution (\ref{def:Gibbs}) and the posterior  (\ref{def:posterior}) 
are identical,  and that $\hat{\Thetav}_{\text{MAP}} \!\left[\mathscr{D}\right]=\lim_{\beta\rightarrow\infty}   \int\!\rmd\Thetav~P_{\beta}(\Thetav\vert\mathscr{D}) \Thetav$ by the Laplace argument\footnote{ In this work  we will mainly rely on the identities $\lim_{M\rightarrow\infty}-\frac{1}{M}\log\int \rmd \x\,\rme^{-M\phi(\x)}=\phi(\x_0) $,  where $\x_0=\argmin_{\x} \phi(\x)$,  and  $\lim_{M\rightarrow\infty}\int \rmd \x \frac{\rme^{-M\phi(\x)}}{\int \rmd \tilde{\x}\, \rme^{-M\phi(\tilde{\x})}} g(\x)=g(\x_0) $ for sufficiently  smooth and well behaved  functions $\phi, g$ of $\x\in\!{\rm I\!R}^p$ with $p= O(M^0)$~\cite{DeBruijn1981}.}. We note that the interpretation of the MAP estimator  (\ref{def:MAP-estimator}) 
in  the SP framework  (\ref{def:Gibbs}) is that $\hat{\Thetav}_{\text{MAP}} \!\left[\mathscr{D}\right]$ is  
the `{ground} state' of  the system.  Regarding the formulation (\ref{def:Gibbs})
we stress that, even though the (inverse) temperatures $\beta=1$ and $\beta\rightarrow\infty$ 
are the most common for inference, the benefits of the generic `thermal' noise  
are well known for optimisation problems in general~\cite{Kirkpatrick1983} and for Bayesian inference in particular~\cite{Iba1999}. 

The Kullback-Leibler  (KL) `distance'~\cite{Cover2012} between the distribution  $P(t, \z\vert\Thetav)$ and its empirical counterpart  $\hat{P}(t,\z\vert \mathscr{D})=N^{-1}\sum_i\delta(t\!-\!t_i)\delta(\z\!-\!\z_i)$, given by  %
\begin{eqnarray}
D(\hat{P}[\mathscr{D}]\vert\vert P_{\Thetav})&=&\!\int\! \rmd t\,  \rmd \z\,  \hat{P}(t,\z\vert \mathscr{D}) \log \Big( \frac{\hat{P}(t,\z\vert \mathscr{D})}{P(t, \z\vert\Thetav)}\Big)~~ \label{def:KL}
\end{eqnarray}
can also be used to obtain the ML  estimator, via  $\hat{\Thetav}_{\text{ML}} \!\left[  \mathscr{D}  \right]=\argmin_{\Thetav} D(\hat{P}\vert\vert P_{\Thetav})$. 
Furthermore, since  $N D(\hat{P}[\mathscr{D}]\vert\vert P_{\Thetav})  = E(\Thetav\vert \mathscr{D})\!+\!\log P(\Thetav)\!-\!N S(\hat{P}[\mathscr{D}])$ where the last term, minus the Shannon entropy of $\hat{P}[\mathscr{D}]$, is independent of $\Thetav$,  the MAP estimator can be obtained via $\hat{\Thetav}_{\text{MAP}}[  \mathscr{D}]=\argmin_{\Thetav}\big\{N D(\hat{P}\vert\vert P_{\Thetav})-\log P(\Thetav)\big\}$. 

Finally, the KL distance (\ref{def:KL}) can also be used to define the difference   $\Delta D(\Thetav,\Thetav_{\!0}\vert   \mathscr{D})=D(\hat{P}[\mathscr{D}]\vert\vert P_{\Thetav}) - D(\hat{P}[\mathscr{D}]\vert\vert P_{ \Thetav_{0}})$, where $\Thetav_{0}$ are the true parameters responsible for the data, which served as a useful measure  of over-fitting   
in ML inference~\cite{Coolen2017}, and was recently extended to MAP inference in generalized linear models~\cite{Coolen2020}.  
Both latter studies used the SP framework, equivalent to (\ref{def:F}), to study \emph{typical} 
(as opposed to \textit{worst-case}) properties of inference in the \emph{high-dimensional} regime 
via the average free energy  $\left\langle F_{\beta}[\mathscr{D}] \right\rangle_{\mathscr{D} }/N$
as computed by the replica method~\cite{Mezard1987}.

\section{Bayesian Linear Regression}
In  Bayesian linear regression (LR) with Gaussian priors, also called  \emph{ridge regression}, 
it is assumed that the data $(\z_i,t_i)$ are for all $i$ sampled independently from the distribution $\mathcal{N}(t\vert \thetav.\z, \sigma^2)P(\z)$, so  the energy (\ref{def:E}), with $\Thetav\equiv\{\thetav,\sigma^2\}$, is given by 
\begin{eqnarray}
E\left(\thetav,\!\sigma^2\vert \mathscr{D} \right)
&=&\frac{1}{2\sigma^2}\left\vert\left\vert \tv -\Zm\thetav\right\vert\right\vert^2  +   \frac{1}{2} \eta \left\vert\left\vert\thetav\right\vert\right\vert^2     +\frac{1}{2}N\log(2\pi\sigma^2)-\log P(\sigma^2)
 \label{eq:LR-E-1},
\end{eqnarray}
%
where $\eta\geq0$ is  the {hyper-parameter} for the Gaussian prior  $P(\thetav)$ and $P(\sigma^2)$ is a generic prior.
The  {true} parameters of $\mathscr{D}$ are written as $\thetav_{0}$ and $\sigma_0^2$, i.e. we assume that 
$\tv=\Zm\thetav_{\!0}+\epsilonv$ with the noise vector $\epsilonv$ being sampled from some  distribution, e.g. the 
multivariate Gaussian $\mathcal{N}(\nullv, \sigma_0^2\I_N)$,  with mean $\nullv$ and covariance $ \sigma_0^2\I_N$. 

The energy function can be expressed as 
%
\begin{eqnarray}
E\left(\thetav,\!\sigma^2\vert\mathscr{D}\right) &=&  \frac{\big(\thetav  - \J_{\sigma^2\eta}^{-1}  \Zm^\Tran\tv  \big)^\Tran  \J_{\sigma^2\eta}  \big(\thetav  - \J_{\sigma^2\eta}^{-1}  \Zm^\Tran\tv  \big)}{2\sigma^2}   + \frac{\tv^\Tran\left(\I_N-\Zm \J^{-1}_{\sigma^2\eta} \Zm^{T}\right)\tv}{2\sigma^2} \nonumber\\
&&~~~~~~~~~~~~~~~~~~~~~~~~~~~~~~~~~~~~~~ +\frac{N\log(2\pi\sigma^2) -2\log P(\sigma^2)}{2}\label{eq:LR-E-2} 
\end{eqnarray}
%
where we defined the $d\!\times\! d$  matrix $\J=\Zm^\mathrm{T}\Zm$,  with elements  $[\J]_{k\ell}=\sum_{i=1}^N z_i(k) z_i(\ell) $,  and 
 its `regularized' version $\J_{\sigma^2\eta}  =\J+\sigma^2 \,\eta\, \I_d$.  
The distribution (\ref{def:Gibbs})  is now
\begin{eqnarray}
P_{\beta}(\thetav,\sigma^2\vert \mathscr{D}) \label{eq:LR-Gibbs} &=&\frac{P_{\beta}\!\left(\thetav\vert\,\sigma^2,\! \mathscr{D}\right)\!
\rme^{-\beta\left[ F_{\beta,\sigma^2}\left[\mathscr{D}\right]  + \frac{1}{2}N\!\log(2\pi\sigma^2)    -\log\! P(\sigma^2) \right]}  }{\int_{0}^\infty \rmd\tilde{\sigma}^2\, \rme^{-\beta \left[F_{\beta,\tilde{\sigma}^2 }\left[\mathscr{D}\right]      +\frac{1}{2}N\log\tilde{\sigma}^2      -\log P(\tilde{\sigma}^2)\right]}   }   ,  \nonumber
\end{eqnarray}
where $P_{\beta}\left(\thetav\vert\sigma^2,\mathscr{D}\right)$ is  the Gaussian distribution 
\begin{eqnarray}
P_{\beta}(\thetav\vert\,\sigma^2,\mathscr{D})&=& \mathcal{N}\big(   \thetav\big\vert\,  \J_{\sigma^2\eta}^{-1}  \Zm^\Tran\tv ,     \,\sigma^2 \beta^{-1}  \J_{\sigma^2\eta}^{-1}   \big)   \label{eq:LR-Gibbs-cond.}
\end{eqnarray}
 with mean $ \J_{\sigma^2\eta}^{-1}  \Zm^\Tran\tv$ and covariance $\sigma^2 \beta^{-1}  \J_{\sigma^2\eta}^{-1} $. We have also defined the \emph{conditional} free energy
\begin{eqnarray}
 F_{\beta,\sigma^2}[\mathscr{D}]
&=&\frac{d}{2\beta} +\frac{1}{2\sigma^2}\tv^\Tran\big(\I_N-\Zm \J^{-1}_{\sigma^2\eta} \Zm^{T}\big)\tv 
-\frac{1}{2\beta} \log\vert 2\pi \rme \sigma^2 \beta^{-1}  \J_{\sigma^2\eta}^{-1} \vert  \label{eq:LR-F-cond.},
\end{eqnarray}
while the full free energy associated with (\ref{eq:LR-Gibbs})   is given by 
\begin{eqnarray}
F_{\beta}[\mathscr{D}]
&=&-\frac{1}{\beta}\log\int\!\rmd\thetav\,  \rmd\sigma^2\, \rme^{-\beta E\left(\thetav,\sigma^2\vert\mathscr{D}\right)  } \label{eq:LR-F}\\
 &=&-\frac{1}{\beta}\log\!\!\int_{0}^\infty\!\!\! \rmd\sigma^2 \rme^{-\beta \left[F_{\beta,\sigma^2}\left[\mathscr{D}\right]  
 +\frac{N}{2}\log(2\pi\sigma^2) - \log P(\sigma^2)    \right]   } \nonumber .
\end{eqnarray}
For $\beta\rightarrow\infty$ the free energy is via the Laplace argument given  by  
$F_{\infty}[\mathscr{D}] =\min_{\thetav,\sigma^2} E\left(\thetav,\sigma^2\vert\mathscr{D}\right)$.
 %
%
$F_{\infty}\left[\mathscr{D}\right]$ is the ground state energy of (\ref{eq:LR-Gibbs}). The 
 ground state 
$\big\{\thetaest[\mathscr{D}], \hat{\sigma}^2[\mathscr{D}]\big\}\!=\!\argmin_{\thetav,\sigma^2} E\left(\thetav,\sigma^2\vert\mathscr{D}\right)$ is given by 
\begin{eqnarray}
\thetaest\left[  \mathscr{D}\right]&=& \J_{\sigma^2\eta}^{-1}  \Zm^\Tran\tv \label{eq:LR-theta-MAP},
\end{eqnarray}
i.e. the mean of (\ref{eq:LR-Gibbs-cond.}), and the solution  of the equation  
\begin{eqnarray}
\sigma^2&=& \frac{1}{N}\big\vert\big\vert \tv -\Zm \thetaest\left[  \mathscr{D}\right]   \big\vert\big\vert^2 + \frac{2\sigma^4}{N}\frac{\partial}{\partial\sigma^2}\log P(\sigma^2) , \label{eq:LR-sigma-MAP}
\end{eqnarray}
corresponding to the MAP estimators of the parameters\footnote{If the inverse-$\chi^2$ distribution is used as a prior for $\sigma^2$, then the MAP estimator for the latter is given by $\sigma^2=\frac{1}{N+\nu+2}+\frac{1}{N+\nu+2}\vert\vert \tv -\Zm\,  \thetaest[ \mathscr{D}] \vert\vert^2$ which suggests that the hyper-parameter $\nu$ has to be \emph{extensive} in order to remain relevant for large $N$. }. From the second line in (\ref{eq:LR-F}) we infer 
 \begin{eqnarray}
F_{\infty}[\mathscr{D}] 
&=&\min_{\sigma^2}\Big[  F_{\infty,\sigma^2\!}\!\left[\mathscr{D}\right]\!+\!\frac{N\log(2\pi\sigma^2)}{2}\!-\! \log\! P(\sigma^2) \Big]\label{eq:LR-F-large-beta-2},~~~~ 
\end{eqnarray}
(again via the Laplace argument),  as well as  for $(N,d)\rightarrow\infty$ the  free energy density  $f_{\beta}\left[\mathscr{D}\right]=\frac{1}{N}F_{\beta}\left[\mathscr{D}\right]$ at any $\beta$:
\begin{eqnarray}
f_{\beta}\!\left[\mathscr{D}\right]
\! &=& \! \min_{\sigma^2}\Big[\!\frac{F_{\beta,\sigma^2}\left[\mathscr{D}\right]}{N}+\frac{\log(2\pi\sigma^2)}{2} -\frac{\log P(\sigma^2)}{N}  \Big].\label{eq:LR-F/N} 
 \end{eqnarray}

For $\beta=1$  the distribution (\ref{eq:LR-Gibbs})  can be used to compute the  MMSE estimators of $\thetav$ and $\sigma^2$, given by the averages 
\begin{eqnarray}
\int_{0}^\infty \!\! \rmd\thetav\,\rmd\sigma^2~P_{\beta}(\thetav,\sigma^2\vert \mathscr{D})\, \thetav &=&     \langle\J_{\sigma^2\eta}^{-1}  \Zm^\Tran\tv \rangle_{\sigma^2}          \label{eq:LR-MMSE-estimators}\\
\int_{0}^\infty \!\! \rmd\thetav\,\rmd\sigma^2~ P_{\beta}(\thetav,\sigma^2\vert \mathscr{D})\,\sigma^2&=&     \langle\sigma^2 \rangle_{\sigma^2},\nonumber  
\end{eqnarray}
where the short-hand  $\langle\cdots\rangle_{\sigma^2}$ refers to averaging over the following marginal
of the distribution (\ref{eq:LR-Gibbs}):
\begin{eqnarray}
P_{\beta}(\sigma^2\vert \mathscr{D}) \! &=& \! \frac{ 
\rme^{-\beta \left[ F_{\beta,\sigma^2}\left[\mathscr{D}\right]  +\frac{N}{2}\log(2\pi\sigma^2) -\log P(\sigma^2)    \right]  }  }{\int_{0}^\infty\! \rmd\tilde{\sigma}^2\, \rme^{-\beta \left[F_{\beta,\tilde{\sigma}^2}\left[\mathscr{D}\right]  +\frac{N}{2}\log(2\pi\tilde{\sigma}^2)  -\log   P(\tilde{\sigma}^2)   \right] } }   \label{eq:LR-sigma-distr}.  ~~~~~
\end{eqnarray}
If  the density $F_{\beta,\sigma^2}\left[\mathscr{D}\right]/N$ is self-averaging then  for $(N,d)\rightarrow\infty$ this marginal is dominated by the solution of (\ref{eq:LR-F/N}). The 
dominant value of $\thetav$ in  (\ref{eq:LR-MMSE-estimators}) is (\ref{eq:LR-theta-MAP}), 
but with $\sigma^2$ being the solution of the following equation, which for $\beta=1$ gives the MMSE estimators, and which  
recovers the MAP estimators  (\ref{eq:LR-theta-MAP}) and (\ref{eq:LR-sigma-MAP}) when $\beta\rightarrow\infty$: 
%
\begin{eqnarray}
\sigma^2&=& \frac{\beta}{(\beta\!-\!\zeta)}\frac{1}{N}\big\vert\big\vert \tv -\Zm \thetaest\left[  \mathscr{D}\right]   \big\vert\big\vert^2 -\frac{\sigma^4\eta}{(\beta\!-\!\zeta)}\frac{1}{N}\Tr\big[  \J^{-1}_{\sigma^2\eta}  \big]       + \frac{2\sigma^4\beta}{(\beta-\zeta)N}\frac{\partial}{\partial\sigma^2}\log P(\sigma^2). \label{eq:LR-sigma-beta}
\end{eqnarray}
%

The free energies (\ref{eq:LR-F-cond.}) and (\ref{eq:LR-F})  obey  the Helmholtz free energy relations. In particular, with $E\left(\thetav\vert\mathscr{D}\right)\!=\!E\left(\thetav,\!\sigma^2\vert\mathscr{D}\right)\!-\!\frac{1}{2}N\log(2\pi\sigma^2)\!+\!\log P(\sigma^2)$ we get
\begin{eqnarray}
F_{\beta,\sigma^2}\left[\mathscr{D}\right]&= &E_\beta[   \mathscr{D}]-T \,\mathrm{S}_\beta[\mathscr{D}]\label{eq:F}, 
\end{eqnarray}
where $T=1/\beta$,  with the average energy 
\begin{eqnarray}
E_\beta[\mathscr{D}]&=&\int\!\rmd\thetav~P_{\beta}(\thetav\vert\sigma^2,   \mathscr{D}) E(\thetav\vert   \mathscr{D}), \label{def:internal-E}
\end{eqnarray}
 and with the \emph{differential entropy} 
\begin{eqnarray}
 \mathrm{S}_\beta[\mathscr{D}]&=&-\int\!\rmd\thetav~P_\beta(\thetav\vert\sigma^2,   \mathscr{D}) \log  P_\beta(\thetav\vert\sigma^2,   \mathscr{D}) 
 \label{def:S}. 
\end{eqnarray}
 In the free energy (\ref{eq:LR-F-cond.}) we have
\begin{eqnarray}
 E_\beta[   \mathscr{D}]&=&\frac{d}{2\beta} +\frac{1}{2\sigma^2}\tv^\Tran\left(\I_N-\Zm \J^{-1}_{\sigma^2\eta} \Zm^{T}\right)\tv\label{eq:LR-int.-E}
\\
 \mathrm{S}_\beta[\mathscr{D}]&=&\frac{1}{2} \log\vert 2\pi\rme\sigma^2 \beta^{-1}  \J_{\sigma^2\eta}^{-1} \vert\nonumber.
\end{eqnarray}
 Furthermore, the average energy can be written as 
\begin{eqnarray}
 E_\beta[\mathscr{D}]=\frac{d}{2\beta} + \min_{\thetav} E(\thetav,\sigma^2\vert   \mathscr{D}).
\end{eqnarray}
We stress that the formulation of LR Bayesian inference via a Gibbs-Boltzmann distribution is not new, 
see, e.g., \cite{ZdeborovaKrzakala2016review,Coolen2017}. 
However, unlike most previous works, here we also consider the case
of unknown $\sigma$ and we keep the temperature generic for most of the analysis, instead of limiting ourselves to the 
familiar cases $T\in\{0,1\}$. 
In the case in which the noise parameter $\sigma^2$ is {known},  i.e. $P(\sigma^2)\!=\!\delta(\sigma^2\!-\!\sigma_0^2)$, 
our free energy expression reduces to $F_{\beta}[\mathscr{D}]\!=\!F_{\beta,\sigma^2_0}[\mathscr{D}]\!-\!\frac{N}{2\beta}\log(2\pi\sigma_0^2)$ 
and $P_{\beta}(\thetav,\sigma^2\vert \mathscr{D})=P_{\beta}(\thetav\vert\sigma_0^2, \mathscr{D})\delta(\sigma^2\!-\!\sigma_0^2)$.

\subsection{Distribution of estimators $\thetaest_{\text{\rm MAP}}$ and $\thetaest_{\text{\rm ML}}$   }
If the noise parameter $\sigma^2$ is independent of the realization of  the data $ \mathscr{D}$, e.g.  $\sigma^2$ is known or {self-averaging} as $(N,d)\!\rightarrow\!\infty$,  and the noise  $\epsilonv$ has Gaussian statistics $\mathcal{N}(\nullv, \sigma_0^2\I_N)$, 
the distribution of the MAP estimator  (\ref{eq:LR-theta-MAP}) is 
\begin{eqnarray}
P(\thetaest)
&=& \left\langle \mathcal{N}\Big(\thetaest\, \vert \,  \J^{-1}_{\sigma^2\eta}  \J\thetav_{\!0} ,\,  \sigma_0^2\J^{-2}_{\sigma^2\eta}   \J    \Big)\right\rangle_{\Zm} \label{eq:MAP-theta-distr-1}.
\end{eqnarray}
For $\eta=0$, i.e. ML inference, and  without averaging over $\Zm$, this recovers Theorem 7.6b in \cite{Rencher2008}.  
To probe the  $(N,d)\!\rightarrow\!\infty$ regime 
we rescale $z_i(\mu)\!\to\! z_i(\mu)/\sqrt{d}$ with now $z_i(\mu)\!=\!{\mathcal O}(1)$. 
This gives $\J\!=\! \C/\zeta$ and $\J_{\sigma^2\eta} = \C_{\zeta\sigma^2\eta}/\zeta$, 
with the sample covariance matrix  $\left [\C\right]_{\mu\nu}\!=\!N^{-1}\sum_{i=1}^N z_i(\mu) z_i(\nu)$ 
and $\C_{\zeta\sigma^2\eta}\!=\!\C\!+\!\zeta\sigma^2\eta \I$, 
so
\begin{eqnarray}
P(\thetaest)&=& \left\langle \mathcal{N}\big(\thetaest\,\vert \, \C^{-1}_{\zeta\sigma^2\eta}\C\thetav_{\!0},\, \zeta \sigma_0^2 \C^{-2}_{\zeta\sigma^2\eta}\C\big)\right\rangle_{\Zm} 
\label{eq:MAP-theta-distr-2}
\end{eqnarray}
Furthermore,  for a Gaussian sample  with {true} covariance matrix $\Cov$, i.e. if each  $\z_i$ in $\Zm$ is drawn independently from    $\mathcal{N}(\nullv,\Cov)$,  the distribution of $\thetaest$ for any {finite} $(N,d)$ is  the Gaussian mixture
\begin{eqnarray}
P(\thetaest)&=& \int\!\!\rmd\C~ \mathcal{W}\!\left(\C\vert\Cov/N,d,N\right) \mathcal{N}\Big(\thetaest\,\vert \, \C^{-1}_{\zeta\sigma^2\eta}[\C]\,\C\thetav_{\!0},\, \zeta\sigma_0^2  \C^{-2}_{\zeta\sigma^2\eta}[\C]\,\C  \Big).~~
\label{eq:MAP-theta-distr-z-Norm}
\end{eqnarray}
The integral is over all symmetric positive definite $d\times d$ matrices, and $\mathcal{W}(\C\vert\Cov/N,d,N)$ is the Wishart distribution, which is  non-singular when  $d\leq N$.  Note that  (\ref{eq:MAP-theta-distr-z-Norm})  also represents  the distribution of `ground states' of  (\ref{eq:LR-Gibbs-cond.}).

For $\eta=0$ the distribution (\ref{eq:MAP-theta-distr-z-Norm}) becomes the multivariate Student's $t$-distribution with $N\!+\!1\!-\!d$ degrees of freedom:
%
\begin{eqnarray}
%
   P(\thetaest)&=&    \frac{  \Gamma\left(\frac{N+1}{2}\right)             }{  \Gamma\left(\frac{N+1-d}{2}\right)} \left\vert  \frac{(1-\zeta+1/N)\Cov}{  \pi\left(  N+1-d  \right) \zeta \sigma_0^2    }  \right\vert^{\frac{1}{2}}  \left[ 1\!  + \! \left( \thetaest \! -\! \thetav_{\!0}   \right)^\Tran \frac{(1-\zeta+1/N)\Cov}{(N+1-d)\zeta\sigma_0^2}\!\left( \thetaest \! -\! \thetav_{\!0}   \right) \right]^{-\frac{N+1}{2}}   ~~~\label{eq:ML-theta-distr}
\end{eqnarray}
%
The vector of {true} parameters  $\thetav_{0}$ is the  mode  and $[\zeta\sigma_0^2/(1\!-\!\zeta\!-\!N^{-1})]\Cov^{-1}$ 
is the covariance matrix of  (\ref{eq:ML-theta-distr}).  In the  regime $(N,d)\rightarrow \infty$ one can 
recover from  (\ref{eq:ML-theta-distr}) the moments of the multivariate Gaussian suggested by the replica method~\cite{Coolen2020}.  
In this regime one indeed finds that any \emph{finite} subset of  components of $\thetaest_{\text{ML}}$   
is  described by a Gaussian distribution~\cite{Nadarajah2005,Kibria2006}.  

\subsection{Statistical properties of the estimator $\hat{\sigma}^2_{\text{\rm ML}}$   }
\begin{figure*}
 \setlength{\unitlength}{0.89mm}
 \begin{picture}(200,90)
 
 \put(0,0){\includegraphics[height=90\unitlength]{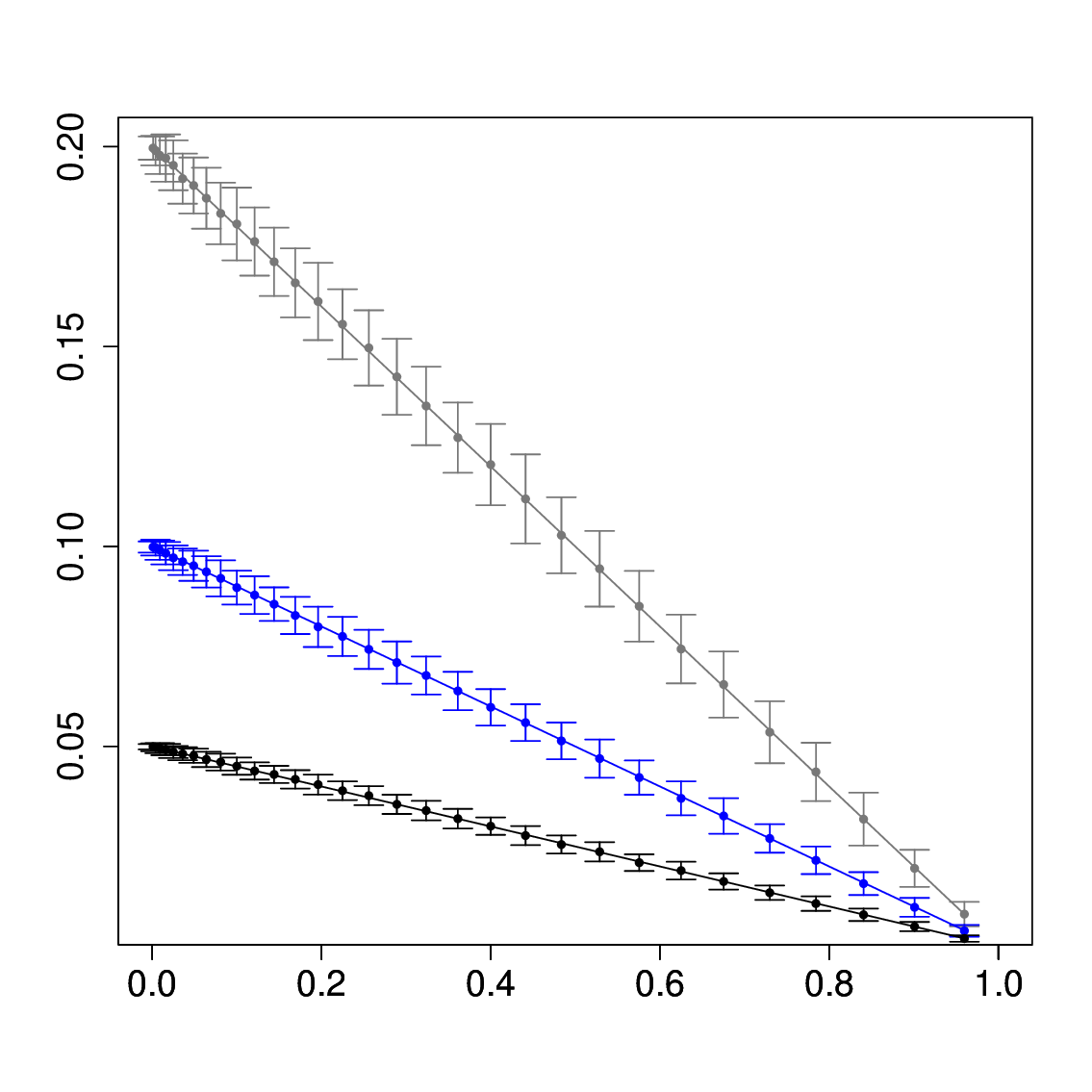}}

\put(100,0){\includegraphics[height=90\unitlength]{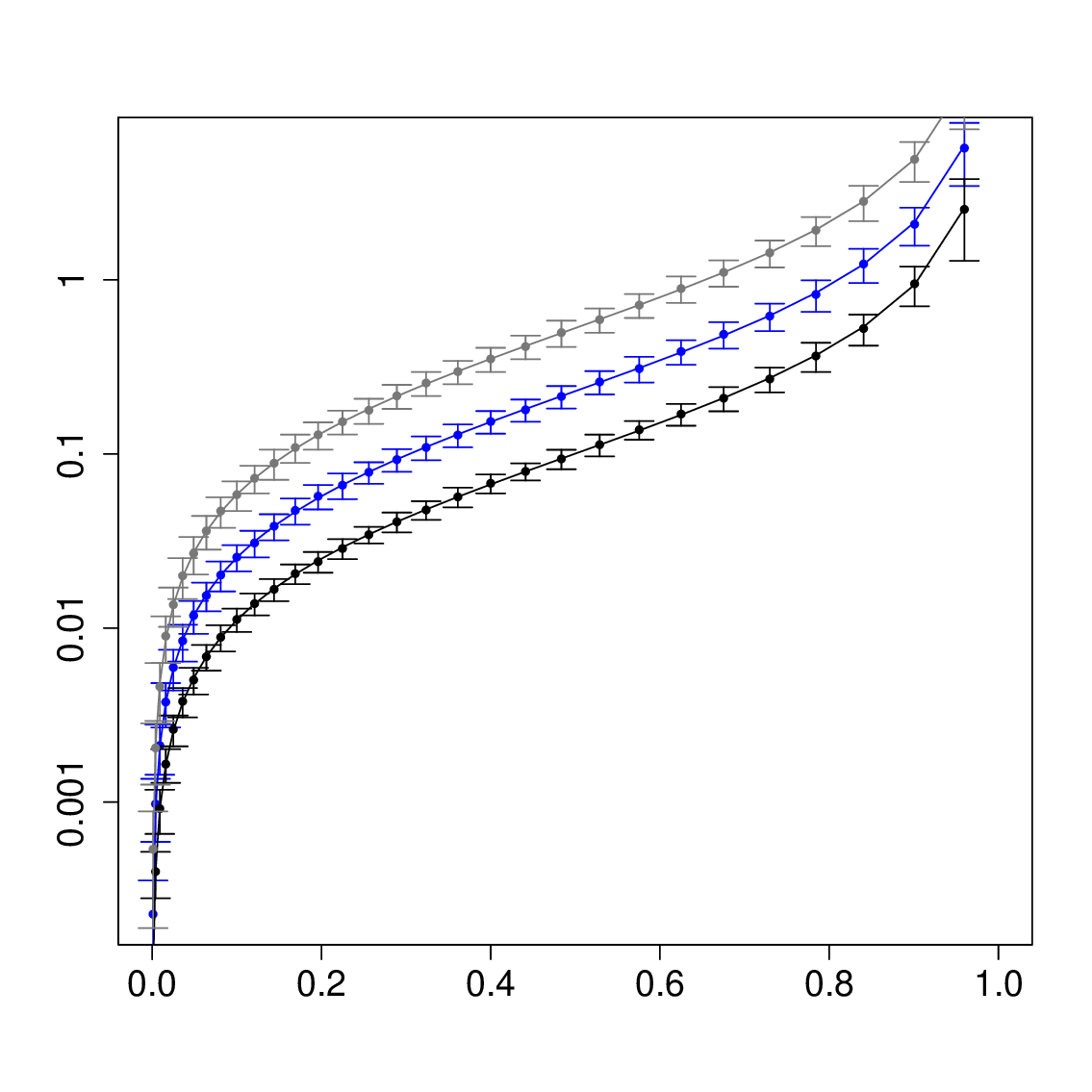}}
\put(-3,77){$(a)$} \put(97,77){$(b)$}
 
  \put(-2 ,50){ $\sigma^2$}  
    
  \put(98 ,40){\rotatebox{90}{
\small{$\frac{1}{d}\vert\vert\thetav_{\!0}\!-\!\thetaest\vert\vert^2$}
} }  
\put(43 ,3){$d/N$} 
          
\put(142 ,3){$d/N$}       
\end{picture}
 \caption{Theoretical predictions for ML inference of the LR model $\tv=\Zm\thetav+\sigma\epsilonv$,
where $\tv\!\in\!{\rm I\!R}^N$ and $\Zm\!\in{\rm I\!R}^{N\times d}$,  in the high-dimensional regime $0<d/N<1$. For each sample the rows of $\Zm$ were sampled from the Gaussian $\mathcal{N}\left(\nullv,\Cov\right)$. The covariance matrix  $\Cov$ is such that  $[\Cov]_{\nu,\nu}=1$,  $[\Cov]_{\nu, \nu+1}=[\Cov]_{\nu+1, \nu}=\epsilon$ with $ 0\leq \epsilon<1$ for $\nu$ odd, and  $[\Cov]_{\nu_1, \nu_2}=0$ for all other $\nu_1\neq \nu_2$.  The density of eigenvalues  of  $\Cov$ is given exactly by  $\rho(\lambda)=\frac{1}{2}\delta(\lambda\!-\!1\!-\!\epsilon)+\frac{1}{2}\delta(\lambda\!-\!1\!+\!\epsilon)$ for any  even  $d$. For each sample the noise vector $\epsilonv$ was sampled from $\mathcal{N}\left(\nullv,\sigma_0^2\I_N\right)$. For each value of $d/N$  the (true) parameter $\thetav_0$ was sampled  from $\mathcal{N}\left(\nullv,\I_d\right)$  only once. (a)  Inferred noise parameter $\sigma^2$ as a function of $d/N$, where $(d,N)$ pairs range from $(10,10^4)$ to $(310,323)$,  plotted for  true value of noise $\sigma_0^2\in\{0.05, 0.10, 0.20\}$ (top to bottom). Solid lines are  the  averages predicted by the theory in (\ref{eq:ML-noise-estimator-moments})  and symbols, with $\pm$ one standard-deviation error-bars, are empirical averages over the 250 samples of data $\{\tv,\Zm\}$.  Error-bars are consistent with the variance predicted by the theory in (\ref{eq:ML-noise-estimator-moments}). 
%
(b) MSE as a function of  $d/N$ plotted for $\epsilon\in\{ 0, 0.75, 0.9\}$ (bottom to top) and $\sigma_0^2=0.1$. Solid lines correspond  to  the theoretical prediction  $\frac{\zeta\sigma_0^2   }{1-\zeta-N^{-1}}    \frac{ 1}{1-\epsilon^2}$ for average MSE, computed via  (\ref{eq:ML-MSE-mean}), and symbols , with $\pm$ one standard-deviation error-bars, are empirical averages over the 250 samples. Note the logarithmic scale of the vertical axis.    Error-bars are consistent with the variance,  $\frac{ 2  \zeta^2\sigma_0^4   }{(1-\zeta)^2 } \frac{1+\epsilon^2}{d(1+\epsilon)^2(1-\epsilon)^2}$, as predicted by (\ref{eq:ML-MSE-mean}). 
}
 \label{figure:estimators-theory} 
 \end{figure*}

For $\eta=0$ the estimator (\ref{eq:LR-theta-MAP})  simplifies considerably   to    
\begin{eqnarray}
\thetaest_{\text{ML}}[\mathscr{D}] &=& \left(\Zm^\Tran\Zm\right)^{-1}  \Zm^\Tran\tv
\end{eqnarray}
 giving us, via  (\ref{eq:LR-sigma-MAP}), the ML noise estimator 
\begin{eqnarray}
\hat{\sigma}^2_{\text{ML}} &=&\frac{1}{N} \epsilonv^\Tran \big(\I_N-\Zm\,    \left(\Zm^\Tran\Zm\right)^{-1}   \Zm^\Tran  \big)\epsilonv \label{eq:ML-noise-estimator}. 
\end{eqnarray}
In particular, if the  noise  $\epsilonv$ originates from a distribution with  mean $\nullv$ and covariance $\sigma_0^2\,\I_N$, 
the mean and variance of  $\hat{\sigma}^2_{\text{ML}}$ are 
\begin{eqnarray}
&&\left\langle\hat{\sigma}^2_{\text{ML}}\right\rangle_{\epsilonv} =\sigma_0^2\left(1\!-\!\zeta\right),~~~~  \Var\left(\hat{\sigma}^2_{\text{ML}}\right) =\frac{2\sigma_0^4}{N}\left(1\!-\!\zeta\right).~~~   \label{eq:ML-noise-estimator-moments}
%
%
\end{eqnarray}
Hence for  $(N,d)\rightarrow\infty$  the noise estimator   (\ref{eq:ML-noise-estimator})  is independent of $\Zm$ and  self-averaging (see Figure  \ref{figure:estimators-theory}).  
Furthermore, for finite $(N,d)$ and $\delta\!>\!0$  the probability of finding an extreme value of $\hat{\sigma}^2_{\text{ML}}\notin \mathcal{I}_{\sigma_0,\delta}$, where  $\mathcal{I}_{\sigma_0,\delta}\equiv\big(\sigma_0^2(1\!-\!\zeta)\!-\!\delta, \sigma_0^2(1-\zeta)\!+\!\delta\big)$,  is 
%
\begin{eqnarray}
%
 %
 \mathrm{Prob} \left[ \hat{\sigma}^2_{\text{ML}}\notin \mathcal{I}_{\sigma_0,\delta}     \right]&=&\mathrm{Prob} \left[N\hat{\sigma}^2_{\text{ML}}\leq N\left(\sigma_0^2(1-\zeta)-\delta\right)\right] + \mathrm{Prob} \left[N\hat{\sigma}^2_{\text{ML}}\geq N\left(\sigma_0^2(1-\zeta)+\delta\right)\right] \nonumber\\
&&\leq   \big\langle\rme^{-\frac{1}{2}\alpha     \left\vert\left\vert\tv-\Zm\thetaest_{\text{ML}}[\mathscr{D}]\right\vert\right\vert^2      }\big\rangle_{\mathscr{D}}  ~\rme^{\frac{1}{2}\alpha N\left(\sigma_0^2(1-\zeta)-\delta\right)}\nonumber\\
&& ~~~~~~~~~~~+\big\langle\rme^{\frac{1}{2}\alpha   \left\vert\left\vert\tv-\Zm\thetaest_{\text{ML}}[\mathscr{D}]\right\vert\right\vert^2   }\big\rangle_{\mathscr{D}}  
 ~\rme^{-\frac{1}{2}\alpha N\left(\sigma_0^2(1-\zeta)+\delta\right)} \label{eq:ML-noise-estimator-tails}.
\end{eqnarray}
%
Assuming that the noise $\epsilonv$ is Gaussian, described by  $\mathcal{N}(\nullv, \sigma_0^2\I_N)$, the  moment-generating  function (MGF)
\begin{eqnarray}
\big\langle\rme^{\frac{1}{2}\alpha    \left\vert\left\vert\tv-\Zm\thetaest_{\text{ML}}[\mathscr{D}]\right\vert\right\vert^2     }\big\rangle_{\mathscr{D}} & =& \rme^{-\frac{N}{2}\left(1- \zeta\right)  \log  \left(1 - \alpha\sigma_0^2\right)} \label{eq:ML-noise-estimator-MGF}  
\end{eqnarray}
is  independent of $\Zm$, allowing us to estimate the fluctuations of $\hat{\sigma}^2_{\text{ML}}$ for $\delta\in\left(0,\sigma_0^2 (1-\zeta) \right)$ via the inequality 
\begin{eqnarray}
%
 %
 \mathrm{Prob} \left[  \hat{\sigma}^2_{\text{ML}}   \notin \mathcal{I}_{\sigma_0,\delta}     \right]
&\leq&
 \sum_{s=\pm1}  \rme^{ -\frac{1}{2}N   \big[    {  \left( 1-\zeta \right) \log 
 \left( {\frac { 1-\zeta}{
  1-\zeta +s\delta/\sigma_0^2 }} \right) +s\delta/\sigma_0^2  }   \big]     }\label{eq:ML-noise-estimator-tails-ub}
\end{eqnarray}

For $\alpha=2\rmi a$ with $a\in{\rm I\!R}$, the MGF  (\ref{eq:ML-noise-estimator-MGF})  becomes the  \textit{characteristic function} (CF)
\begin{eqnarray}
 \big\langle\rme^{\rmi a  \left\vert\left\vert\tv-\Zm\thetaest_{\text{ML}}[\mathscr{D}]\right\vert\right\vert^2  }\big\rangle_{\mathscr{D}} 
 &=&  \left(1 - \rmi a\,2\sigma_0^2\right)^{-\frac{1}{2}N\left(1- \zeta\right)}     \label{eq:ML-noise-estimator-CF} 
\end{eqnarray}
of the random variable $\big\vert\big\vert\tv\!-\!\Zm\thetaest_{\text{ML}}[\mathscr{D}]\big\vert\big\vert^2$.  
Note that (\ref{eq:ML-noise-estimator-CF}) is the CF of the {gamma} distribution  
(see  Theorem 7.6b in \cite{Rencher2008}), with mean $N\sigma_0^2(1-\zeta)$ and 
variance  $N2\sigma_0^4(1- \zeta) $. 
Mean and variance  of $\hat{\sigma}^2_{\text{ML}}$  are  $\sigma_0^2(1\!-\! \zeta)$ and  
$ 2\sigma_0^4\left(1\!-\! \zeta\right) /N$,  respectively.  
For $\sigma_0=1$ we obtain that $N \hat{\sigma}^2_{\text{ML}}$ is 
a chi-square distribution with $N(1-\zeta)$ degrees of freedom, as expected from Cochran's theorem \cite{cochran_1934}.

Finally, it follows from (\ref{eq:ML-noise-estimator-moments}) and (\ref{eq:LR-sigma-beta}) that  
the finite temperature   ML  noise estimator   in the high-dimensional regime is given by  
\begin{eqnarray}
\hat{\sigma}^2_{\text{ML}} =\frac{\beta}{\beta-\zeta}\sigma_0^2(1-\zeta) \label{eq:ML-noise-estimator-finite-T}.
\end{eqnarray}
We observe that for $\beta=1$ we obtain unbiased estimation of  $\sigma^2$.  
The latter suggests that `thermal' noise, controlled by $\beta$, 
is beneficial for  the ML  inference of $\sigma^2$. However,  that this is indeed the case, and that the value $\beta=1$ is `special', is not  \emph{a priori} obvious for this model. Our development confirms the result obtained in evaluating the average free energy with the replica method~\cite{Coolen2020}.

\subsection{ Statistical properties  of MSE  in ML inference  }
Using the distribution  (\ref{eq:ML-theta-distr}) and with the eigenvalues $ \lambda_1(\Cov)\!\leq\! \lambda_2 (\Cov)\!\leq\!\cdots\!\leq\! \lambda_d(\Cov)$ of the true (population) covariance matrix $\Cov$, the CF  of  the MSE, defined as $\frac{1}{d}\vert\vert\thetav_{\!0}-\thetaest_{\rm ML}[\mathscr{D}]\vert\vert^2$  for finite $(N,d)$,  can be written as
\begin{eqnarray}
\big\langle   \rme^{\rmi\alpha  \vert\vert\thetav_{\!0}-\thetaest_{\text{ML}} [\mathscr{D}]\vert\vert^2 }\big\rangle_{\mathscr{D}}
&=& \int_{0}^\infty\!\!\rmd\omega~\Gamma_{N+1-d}\,(\omega)
\prod_{\ell=1}^d \!\left(\! 1 \! -\!  \frac{\rmi \alpha 2    \zeta\sigma_0^2   }{\omega(1\!-\!\zeta\!+\!N^{-1}) \lambda_\ell (\Cov)} \!  \right)^{\!-\frac{1}{2}} \!\!\!\!\!\label{eq:ML-MSE-CF},~~~~~~
\end{eqnarray}
with the gamma  distribution $\Gamma_{\nu}\,(\omega)=\frac{\nu^{\nu/2}}{2^{\nu/2}\Gamma(\nu/2)}\omega^{\frac{\nu-2}{2}}\rme^{-\frac{1}{2}\nu\omega}$ for $\nu>0$.
%
%
The last term  in (\ref{eq:ML-MSE-CF}) is the product of  CFs of gamma  distributions with the same `shape' parameter $1/2$, but  different `scale'  parameters  $2\zeta\sigma_0^2 /\omega(1\!-\!\zeta\!+\!N^{-1}) \lambda_\ell (\Cov)$.
From (\ref{eq:ML-MSE-CF}) one obtains mean and variance of MSE:
\begin{eqnarray}
\mu(\Cov) =\frac{1}{d}\big\langle   \vert\vert\thetav_{\!0}\!-\!\thetaest_{\text{ML}}[\mathscr{D}]\vert\vert^2 \big\rangle_{\mathscr{D}}
&=&\frac{\zeta\sigma_0^2   }{1\!-\!\zeta\!-\!N^{-1}}    \frac{  \Tr[\Cov^{-1}] }{d}, \nonumber
\\
 \Var \Big(\frac{1}{d}\vert\vert\thetav_{\!0}\!-\!\thetaest_{\text{ML}}[\mathscr{D}]\vert\vert^2\Big)
&=& \frac{ 2  \zeta^2\sigma_0^4   }{(1\!-\!\zeta)^2 }  \frac{ \Tr[\Cov^{-2}]}{d^2}.   \label{eq:ML-MSE-mean}
\end{eqnarray}
The latter gives us the condition for self-averaging of the MSE, i.e.
$ \Var \big(\frac{1}{d}\vert\vert\thetav_{\!0}-\thetaest_{\text{ML}}[\mathscr{D}]\vert\vert^2\big)\rightarrow0$ as $(N,d)\rightarrow\infty$.  We note that  (\ref{eq:ML-MSE-mean}) suggests that MSE  is dominated by  the smallest eigenvalue of  the true covariance $\Cov$  and hence it can grow with an  increase in the covariate correlations  (see Figure  \ref{figure:estimators-theory}).

 We finally consider deviations of $\frac{1}{d} \vert\vert\thetav_{\!0}\!-\!\thetaest_{\text{ML}}[\mathscr{D}]\vert\vert^2  $ from 
 the mean $\mu(\Cov)$ given in (\ref{eq:ML-MSE-mean}).  
 The probability of observing the event event  
 $\frac{1}{d}\vert\vert\thetav_{\!0}\!-\!\thetaest_{\text{ML}}[\mathscr{D}]\vert\vert^2 \!\notin\! \mathcal{I}_{\mu,\delta} $, where $\mathcal{I}_{\mu,\delta}\equiv \big( \mu(\Cov) \!-\!\delta,\mu(\Cov) \!+\!\delta    \big)$ for $\delta>0$, 
  is bounded from above as follows 
  %
\begin{eqnarray}
%
%
\mathrm{Prob} \left[  \frac{1}{d}\vert\vert\thetav_{\!0}-\thetaest_{\text{ML}}[\mathscr{D}]\vert\vert^2  \notin \mathcal{I}_{\mu,\delta}     \right]&=&\mathrm{Prob} \left[  \vert\vert\thetav_{\!0}-\thetaest_{\text{ML}}[\mathscr{D}]\vert\vert^2         \leq d\left( \mu\left(\Cov\right)  -   \delta\right)\right] \nonumber\\
&& + \mathrm{Prob} \left[   \vert\vert\thetav_{\!0}-\thetaest_{\text{ML}}[\mathscr{D}]  \vert\vert^2 \geq d\left(  \mu\left(\Cov\right)  +\delta\right)       \right] \nonumber\\
&&~~~~~\leq\mathrm{C}_{-}\rme^{-N \Phi_{-}\left[\alpha, \mu(\lambda_d ), \delta\right]  }  + \mathrm{C}_{+}\rme^{-N\Phi_{+}\!\left[\alpha, \mu(\lambda_1 ), \delta\right]}    \label{eq:ML-MSE-dev-ub}.
\end{eqnarray}
%
with some small $\alpha>0$ and  positive  constants $\mathrm{C}_{\pm}$.  For the rate function $\Phi_{-}[\alpha, \mu(\lambda_d ), \delta]$ to be positive  for arbitrary small  $\delta$ it is sufficient that   $\mu(\lambda_d)\geq1$, where $\mu(\lambda)=\zeta\sigma_0^2 /(1\!-\!\zeta)\lambda$, while  for $\mu(\lambda_d)<1$  for this to happen the $\delta$ values must satisfy $\delta>1\!-\!\mu(\lambda_d)$.  The rate function $\Phi_{+}[\alpha, \mu(\lambda_1 ), \delta]$ is positive  for any  $\delta\in(0,\mu(\lambda_1))$. 

\subsection{Statistical properties of the free energy }
We consider the free energy (\ref{eq:LR-F-cond.}) for finite inverse temperature  $\beta$ and finite  $(N,d)$.   
Assuming that the  noise  $\epsilonv$ has  mean $\nullv$ and covariance $\sigma_0^2\,\I_N$, and that the parameter $\sigma^2$ 
is independent of $\mathscr{D}$,   the average free energy is 
%
\begin{eqnarray}
%
%
 \big\langle F_{\beta, \sigma^2}[\mathscr{D}] \big\rangle_{\mathscr{D}}
&=&\frac{d}{2\beta}+\frac{1}{2\sigma^2}      \thetav_{\!0}^\Tran     \big\langle\big(\J-\J \J^{-1}_{\sigma^2\eta} \J\big) \big\rangle_{\Zm} \thetav_{\!0}    + \frac{\sigma_0^2}{2\sigma^2} \big(N-\big\langle\Tr[\J \J^{-1}_{\sigma^2\eta}] \big\rangle_{\Zm} \big)\nonumber\\
&&   -\frac{1}{2\beta} \big\langle \log\vert 2\pi\rme\,\sigma^2 \beta^{-1}  \J_{\sigma^2\eta}^{-1} \vert \big\rangle_{\Zm} \label{eq:LR-cond-F-aver}
\end{eqnarray}
%
Under the same assumptions, the \emph{variance} of  $F_{\beta, \sigma^2}[\mathscr{D}]$ can be obtained by  exploiting  the   Helmholtz  free energy representation  $F_{\beta, \sigma^2}[\mathscr{D}] = E_\beta[\mathscr{D}] - T    \mathrm{S}_\beta[\mathscr{D}]$, giving us
\begin{eqnarray}
 \Var \left(F_{\beta, \sigma^2}\left[\mathscr{D}\right]\right)&=&\Var \left(E_{\beta}[\mathscr{D}]\right)+T^2 \Var\left(  \mathrm{S}_{\beta}[\mathscr{D}]  \right)
-2T\,\Covar\left( E_{\beta}[ \mathscr{D}], \mathrm{S}_{\beta}[\mathscr{D}]  \right)\label{eq:F-Var-identity}.
\end{eqnarray}
The full details on  each term in (\ref{eq:F-Var-identity}) are found in the Appendix.

\subsubsection{Free energy of ML inference }
For $\eta=0$ and after transforming $z_i(\mu) \to z_i(\mu)/\sqrt{d}$ for all $(i,\mu)$, with $z_i(\mu)={\mathcal O}(1)$, expression (\ref{eq:LR-cond-F-aver}) gives the average free energy density   
%
\begin{eqnarray}
%
 %
 \Big\langle \frac{F_{\beta, \sigma^2}\left[\mathscr{D}\right]}{N} \Big\rangle_{\mathscr{D}}
 &=&\frac{1}{2} \frac{\sigma_0^2}{\sigma^2} (1\!-\!\zeta)+\frac{\zeta}{2\beta}\log\Big(\frac{\beta}{2\pi\sigma^2\zeta}\Big) +\frac{\zeta}{2\beta}\int\!\rmd \lambda~ \rho_d (\lambda)\log (\lambda), \label{eq:LR-cond-F/N-aver-ML}
\end{eqnarray}
%
where we  defined the average eigenvalue density $\rho_d (\lambda)=\left\langle\rho_d (\lambda\vert\Zm )\right\rangle_{\Zm} $ of the empirical covariance matrix, with
\begin{eqnarray}
\rho_d (\lambda \vert\Zm)&=&\frac{1}{d}\sum_{\ell=1}^d\delta\big[\lambda-\lambda_\ell(\Zm^\Tran\Zm/N)\big].
\label{eq:eigenvalue-density}
\end{eqnarray}
The variance of free energy density is  
%
\begin{eqnarray}
%
 %
 \Var \Big(\frac{F_{\beta,\sigma^2}\left[\mathscr{D}\right]}{N}\Big)
 &=&\Var \Big(\frac{E[\mathscr{D}]}{N}\Big)+T^2 \Var\Big(\frac{\mathrm{S}(P[\mathscr{D}]) }{N}   \Big)\nonumber\\
 &=&   \frac{\zeta^2}{4\beta^2}\!\int\!\rmd \lambda\,\rmd \tilde{\lambda}~ C_d(\lambda, \tilde{\lambda})\log (\lambda) \log (\tilde{\lambda}) 
 +\frac{\sigma_0^4(1\!-\!\zeta)}{2\sigma^4N}  \label{eq:LR-cond-F/N-Var-ML}
\end{eqnarray}
%
where  we defined the correlation function  $C_d(\lambda, \tilde{\lambda})=\big\langle\rho_d (\lambda\vert\Zm )\rho_d (\tilde{\lambda}\vert\Zm )\!\big\rangle_{\Zm} -\big\langle\rho_d (\lambda\vert\Zm )\!\big\rangle_{\Zm}\! \big\langle \rho_d (\tilde{\lambda}\vert\Zm )\!\big\rangle_{\Zm}$. 

Clearly, if  $ \int\!\rmd \lambda\,\rmd \tilde{\lambda}~ C_d(\lambda, \tilde{\lambda}) f (\lambda,\tilde{\lambda})\!\rightarrow\! 0$ as $(N,d)\rightarrow\infty$,  for any smooth function $f(\lambda,\tilde{\lambda})$,  then the free energy density $f_{\beta}[\mathscr{D}] =F_{\beta}\left[\mathscr{D}\right]/N $ is  self-averaging.    

\begin{figure}[t]
%
 \setlength{\unitlength}{0.262mm}
\begin{picture}(350,280)
\put(0,0){\includegraphics[width=310\unitlength, height=280\unitlength]{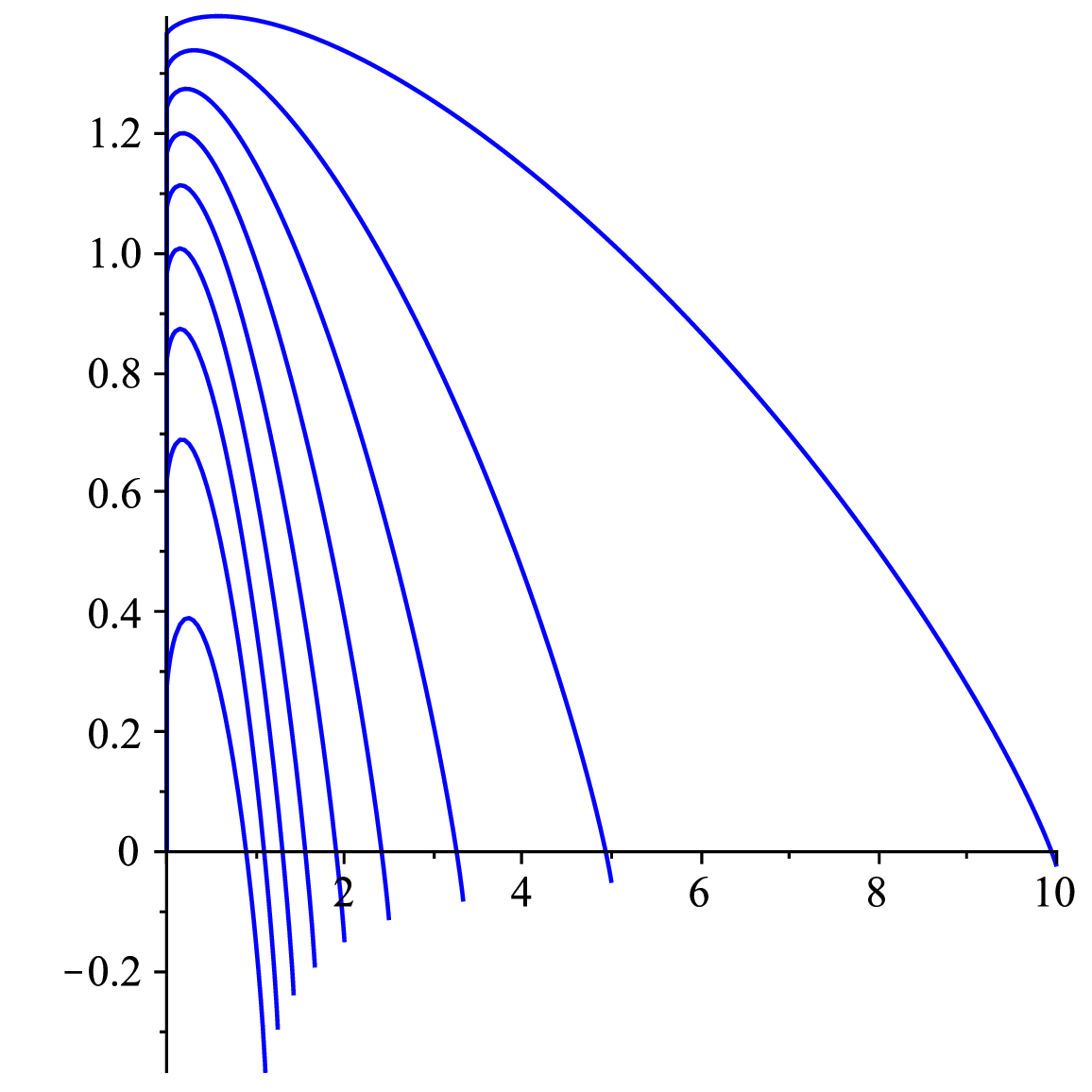}}
\put(-2,150){$f_\beta$} 
\put(180,8){$1/\beta$}
\end{picture}
%
\caption{Asymptotic free energy density $f_\beta=\lim_{N\to\infty}  f_{\beta}[\mathscr{D}]$ of finite temperature ML inference as a function of temperature $T=1/\beta$, plotted  for $\zeta\in\{1/10, 2/10,  \ldots, 9/10\}$ (from right to left)  
in the high-dimensional regime where $N,d\rightarrow\infty$ with fixed ratio $\zeta=d/N$. For $\beta\rightarrow\infty$ it approaches  the value $\frac{1}{2}\log[2\pi \rme\sigma_0^2(1-\zeta)]$.   For $\beta\rightarrow\zeta$ it approaches  $\frac{1}{2\zeta}[\zeta\log(1\!-\!\zeta)\!-\!\log(1\!-\!\zeta)\!-\!\zeta]$, and for $\beta\!\in\!(0,\zeta)$, i.e. in the high `temperature' region $T\in(1/\zeta, \infty)$, the free energy density is $-\infty$.  Here the  true noise  parameter is $\sigma_0^2=1$ and  the true data covariance  matrix  is $\I_d$.\hspace*{-2mm}    \label{figure:3d}
}
\end{figure}

Finally, if we use  (\ref{eq:LR-cond-F/N-aver-ML}) in the  free energy density (\ref{eq:LR-F/N})  for $\eta=0$, and assume Gaussian data with true population covariance matrix $\Cov=\I_d$, then  we find   
 %
%
 \begin{eqnarray}
&&\lim_{N\rightarrow\infty} f_{\beta}[\mathscr{D}]~=\nonumber\\
&&~~~\left\{\!\!\begin{array}{lr}
\frac{ \beta-\zeta}{2\beta}   \log  \Big( {\frac {2\pi \sigma_0^2\left( 1-\zeta \right) }{\beta-\zeta}} \Big)+          \frac{ \log \beta +1}{2}   -\frac {1}{2\beta}\Big(\zeta \log \zeta  \!+\!( 1 \!-\!\zeta)  \log ( 1\!-\!\zeta) \!+\!2\zeta\Big)~~{\rm if}~\beta\!>\!\zeta \\
-\infty \hspace*{110mm} {\rm if}~ \beta\!\in\!(0,\zeta)
 \end{array}\!\!\right.,
~~~~ \label{eq:LR-F/N-aver-ML-Cov-I}
\end{eqnarray}
%
%
with  the convention $0\log0\!=\!0$. 
%
%
Since for $\lambda\!\in\![a_{-}, a_{+}]$ and $0\!<\!\zeta\!<\!1$ the eigenvalue spectrum 
$\rho_d (\lambda\vert\Zm )$ converges to  $(2\pi\lambda\zeta)^{-1}\sqrt{(\lambda\!-\!a_{-})(a_{+}\!-\!\lambda)}$ in a distributional sense as $(N,d)\!\to\!\infty$ ~\cite{Gotze2004},  with $a_{\pm}\!=\!(1\!\pm\!\sqrt{\zeta})^2$, 
 the free energy density  is self-averaging. Its values are plotted versus the temperature in Figure \ref{figure:3d}.
Furthermore, the average free energy density   (\ref{eq:LR-F/N-aver-ML-Cov-I}) is identical to that of \cite{Coolen2020}. 
Since $\lim_{\beta\downarrow\zeta}\lim_{N\rightarrow\infty} f_{\beta}[\mathscr{D}]$ is finite,  the system exhibits a {zeroth}-order phase  
transition~\cite{Chavanis2002} at $T=1/\zeta$. 

\subsubsection{Free energy of MAP inference }
We next assume that the true parameters  	$\thetav_{0}$ are drawn at random,  with mean $\nullv$ and covariance matrix $S^2\I_d$.  As before we rescale $z_i(\mu)\to z_i(\mu)/\sqrt{d}$ where $z_i(\mu)={\mathcal O}(1)$, and define 
  $\J=\C/\zeta$ (so that $\C=\Zm^\Tran\Zm/N$)  and $\C_{\zeta\sigma^2\eta} = \zeta \J_{\sigma^2\eta} $.
Then the average of (\ref{eq:LR-cond-F-aver}) over $\thetav_{0}$ becomes 
\begin{widetext}
\begin{eqnarray}
\Big\langle \! \Big\langle      \frac{F_{\beta, \sigma^2}\left[\mathscr{D}\right]}{N}  
\Big\rangle_{\!\mathscr{D}}\Big\rangle_{\!\thetav_{0}}
 &=& 
\frac{\zeta}{2\beta}\!+\! \frac{1  }{2}\!\int\!\! \rmd \lambda\,\rho_d (\lambda)\! \left[     \frac{S^2   \zeta\eta \lambda}{\lambda\!+\!\zeta\sigma^2\eta}
 +  \frac{\sigma_0^2}{\sigma^2}\! \left[1\!-\!     \frac{ \zeta  \lambda}{\lambda\!+\!\zeta\sigma^2\eta}
  \right]  +\frac{\zeta}{\beta}\!    \log \left(\lambda\!+\! \zeta\sigma^2\eta   \right) \right]\nonumber\\
   &&-\frac{\zeta}{2\beta}\!   \log\!\left( 2\pi\rme\sigma^2\! \beta^{-1}\!  \zeta \right)
 \label{eq:LR-cond-F/N-aver-MAP} .
 \end{eqnarray}
 \end{widetext}
Furthermore,  using (\ref{eq:F-Var-identity}),  we obtain, under the same assumptions,   that $ \Var (   F_{\beta,\sigma^2}\left[\mathscr{D}\right]/N) $ is of the form (see Appendix):
\begin{eqnarray}
\Var \Big( \!  \frac{F_{\beta,\sigma^2}\!\!\left[\mathscr{D}\right]}   {N} \!    \Big)  \!  
=\!\!\!\int\!\!\rmd \lambda\,\rmd \tilde{\lambda}\,C_d(\lambda,\! \tilde{\lambda})\Phi(\lambda,\! \tilde{\lambda})\!+\! {\mathcal O}\!\left(\!\frac{1}{N}\!\right)\label{eq:LR-cond-F/N-Var-MAP}.
%
\end{eqnarray}
Hence for $\eta>0$ the conditional free energy is self-averaging with respect to the realization of the true parameter if the spectrum $\rho_d (\lambda\vert\Zm )$ is self-averaging (since then $C_d(\lambda,\tilde{\lambda})\to 0$ as $(N,d)\to\infty$). The latter, under the same assumptions, is the condition for  the MAP estimator (\ref{eq:LR-sigma-beta}) of the noise $\sigma^2$ to be self-averaging (see Appendix)  and hence the  free energy (\ref{eq:LR-F/N}) is self-averaging if $\rho_d (\lambda\vert\Zm )$ is self-averaging.   This is the case e.g. for Gaussian data with covariance matrix $\Cov=\I_d$.

 \section{Discussion and outlook}
 
In this paper we derived exact results for the Bayesian model (\ref{def:posterior}) of the linear regression $\tv=\Zm\thetav+\sigma\epsilonv$,
where $\tv\!\in\!{\rm I\!R}^N$ and $\Zm\!\in{\rm I\!R}^{N\times d}$.    
Mapping this to a Gibbs-Boltzmann  distribution   (\ref{eq:LR-Gibbs}), with finite (inverse) `temperature'  $\beta$, 
allowed us to investigate properties of several  inference protocols~\cite{Rencher2008} 
for finite $N$ (sample size),  $d$ (dimension)  and in the (high-dimensional)  limit $(N,d)\rightarrow\infty$. In particular we studied statistical properties of  free energy which is the main object of interest  in statistical physics approaches 
to inference (see~\cite{Coolen2020} and references therein). 

If the noise strength $\sigma^2$ is known and the  distributions of the data $\Zm$ and the noise $\epsilonv$ are Gaussian, 
the distribution  of the MAP estimator $\thetaest_{\text{MAP}}$ of $\thetav$ 
is the Gaussian mixture  (\ref{eq:MAP-theta-distr-z-Norm}), for any   
finite $(N,d)$.  
We used  (\ref{eq:MAP-theta-distr-z-Norm}) to show that the distribution of  ML estimator $\thetaest_{\text{ML}}$  is  
the Student's $t$-distribution (\ref{eq:ML-theta-distr}).  The consequence of this is that its marginal,  
the univariate Student distribution (which can be used in statistical hypothesis testing to calculate p-values), 
has `fat' tails for  finite  $(N,d)$. However, any marginal of (\ref{eq:ML-theta-distr}) that describes a finite number of components of $\thetaest_{\text{ML}}$ has  a  Gaussian form when $(N,d)\rightarrow\infty$. 

Also, for any choice for the distributions of $\Zm$ and $\epsilonv$, 
the ML estimator  $\hat{\sigma}^2_{\text{ML}}$ of the noise parameter $\sigma^2$ is \emph{self-averaging} , i.e. its variance is vanishing as $(N,d)\rightarrow\infty$. Furthermore, deviations of  $\hat{\sigma}^2_{\text{ML}}$ from its mean, estimated by the bound in (\ref{eq:ML-noise-estimator-tails-ub}), are exponentially suppressed  in  $(N,d)$ for Gaussian $\epsilonv$.  As a consequence the inference of $\hat{\sigma}^2_{\text{ML}}$ is almost deterministic even for 
moderate values of  $(N,d)$.   This result is independent of  $\Zm$.

We used the distribution of the ML estimator (\ref{eq:ML-theta-distr}) to derive  
the characteristic function of the MSE (\ref{eq:ML-MSE-CF}).   The latter was used to derive  
the  mean and variance (\ref{eq:ML-MSE-mean}), giving a  condition 
for the MSE to be  self-averaging as  $(N,d)\rightarrow\infty$,  and to estimate deviations of MSE from its mean, 
given by  the bound (\ref{eq:ML-MSE-dev-ub}), for finite $(N,d)$. The result (\ref{eq:ML-MSE-mean}) suggests that
the deviations of $\thetaest_{\text{ML}}$ from $\thetav_0$ 
can grow significantly with covariate correlations, 
proportional to $\Tr\left[ \Cov^{-1}\right]$, thus leading to severe inefficacy in the inference of $\thetav$.

If we assume that the noise parameter $\sigma^2$ is known, we obtain average (\ref{eq:LR-cond-F/N-aver-ML}) and variance (\ref{eq:LR-cond-F/N-Var-ML}) 
of  the conditional free energy density (\ref{eq:LR-F-cond.}) of  ML inference with \emph{finite} temperature $T=1/\beta$ and for finite $(N,d)$.  This result is independent of the  distributions of the data $\Zm$ and the noise $\epsilonv$.  
For finite $T$, the noise estimator $\hat{\sigma}^2_{\text{ML}}$, given by (\ref{eq:LR-sigma-beta}) with $\eta=0$, 
is self-averaging when $(N,d)\rightarrow\infty$. The same is true for the   free energy density (\ref{eq:LR-F/N})  
if the density of eigenvalue  \emph{spectrum}  (\ref{eq:eigenvalue-density}) of the covariance matrix $\Zm^\Tran\Zm/N$ is self-averaging.  The latter is true if   $\Zm$ is sampled from a Gaussian with $\Cov=\I_d$. In that case, and upon assuming  Gaussian noise $\epsilonv$,  the  free energy density  (\ref{eq:LR-F/N}) recovers 
the result obtained by the replica method~\cite{Coolen2020}.  
The ML estimator  $\hat{\sigma}^2_{\text{ML}}$ diverges at $\beta=\zeta$, and the  free energy density (\ref{eq:LR-F/N})  
is \emph{discontinuous} at this value of $\beta$. This is an instance where the presence of the thermal noise
with finite generic $\beta$ allows us to derive an interesting result.
Another is Eq (\ref{eq:ML-noise-estimator-finite-T}) for the finite temperature ML noise estimator.

 The additional assumption that the  true parameters $\thetav_0$ are drawn at random, 
with  mean $\nullv$ and covariance  $S^2\I_d$, allows us to derive average (\ref{eq:LR-cond-F/N-aver-MAP}) and 
variance (\ref{eq:LR-cond-F/N-Var-MAP}) of the conditional free energy (\ref{eq:LR-F-cond.})  of MAP inference  for finite $T$  and finite $(N,d)$.   We also computed the variance of the MAP estimator (\ref{eq:LR-sigma-beta}) of the noise parameter $\sigma^2$.   These  results  are again independent of the distributions of the data $\Zm$  and the noise $\epsilonv$.   We find that  the free energy (\ref{eq:LR-F/N}) is self-averaging  if the  spectrum of the empirical covariance matrix $\Zm^\Tran\Zm/N$ is self-averaging 
as $(N,d)\rightarrow\infty$.   

The above results emphasize that still much can be learned  about high-dimensional  
Bayesian linear regression from  exact calculations with standard methods. While we present this as a minimal model of inference 
in the high dimensional setting, linear regression is commonly used in many areas of research. 
For example, linear regression models are used extensively in the statistical analysis of genetic data. 
Genome-wide association studies, where the aim is to undercover effect sizes for each Single Nucleotide Polymorphisms (SNPs), 
often use extremely high-dimensional datasets. Here the number $N$ of individuals is $O(10^3)$ and the number $d$ of SNPs is $O(10^6)$ 
with correlations occurring due to the phenomenon of genetic linkage. Another biological example is the analysis of gene expression data 
where due to nature of biological pathways involved \cite{Segal2003}, correlations pose a significant challenge in uncovering 
true associations in data \cite{Balding2006}.

Many questions remain still open and we hope that this paper may contribute to future work in this direction.
Some results appear well within reach, such as the extension to sub-Gaussian noise 
for the argument that leads to (\ref{eq:ML-noise-estimator-tails-ub}), 
employing techniques  used in \cite{Rudelson2013}. 
Other results are less immediate but seem feasible, such as extending some of the ML results to MAP inference, 
starting from evaluation of the distribution of $\hat{\theta}_\text{MAP}$ 
 (\ref{eq:MAP-theta-distr-z-Norm}) for  $(N,d)\rightarrow\infty$.
Another interesting line of work would be to try to extend our present results to generalized linear models (GLMs), a very similar distribution for the estimator $\hat{\theta}_\text{MAP}$ has already been conjectured through the use of the replica method, \cite{Coolen2020}.
Other crucial investigations, such as a rigorous analytical study of the effect of 
model mismatch, appear instead to be still quite challenging 
with current techniques.  Overall we expect high dimensional linear regression to serve as a starting point to tackle more realistic scenarios, which should include among other things, correlation between data and noise, and dimensional mismatch between the teacher and the student model.
\vspace*{8mm}

\begin{acknowledgments} 
\noindent
We would like to thank Dr  Heather Battey and Sir David Cox for  reading a first draft of this manuscript and valuable suggestions. 
AM is supported by Cancer Research UK (award C45074/A26553) and the UK's Medical Research Council (award MR/R014043/1).
FA is supported by EIT Health ID 19277 (RGS@HOME) under H2020. FAL gratefully acknowledges financial support through a scholarship from Conacyt (Mexico).
\end{acknowledgments}
\appendix

\section{Ingredients}
We write  $\I_N$ for $N\times N$ identity matrix.  The  data $\mathscr{D} =\{\tv,\Zm\}$,  where   $\tv\in \mathbb{R}^{N}$ and $\Zm=\big(\z_1, \ldots,\z_N \big)$ is the $N\times d$ matrix,    is  a set    of observed   `input-output' pairs $\{(\z_1, t_1), \ldots, (\z_N,t_N)\}$ generated by the process
 \begin{eqnarray}
   \tv&=&\Zm\thetav_0 +\epsilonv\label{eq:process}. 
\end{eqnarray}
The \emph{true} parameter vector  $\thetav_0\in \mathbf{R}^d$ is unknown, and the vector $\epsilonv\in \mathbf{R}^N$ represents noise with mean $\nullv$ and covariance  $\sigma_0^2\I_N$, with also the  \emph{true} noise  parameter  $\sigma_0^2$  unknown to us.  The  (empirical) covariance matrix of the  input data is 
\begin{eqnarray}
\J[\Zm]=\Zm^\mathrm{T}\Zm   \label{def:J},
\end{eqnarray}
 where  $[\J[\Zm]]_{k\ell}=\sum_{i=1}^N z_i(k) z_i(\ell)$. To  simplify notation we will sometimes omit the dependence on $\Zm$ and write  $\J\equiv\J[\Zm]$.  The \emph{maximum a posteriori} estimator (MAP) of $\thetav_0$ in linear regression with Gaussian prior $\mathcal{N}(\nullv, \eta^{-1}\I_d)$  is 
\begin{eqnarray}
\thetaest\left[  \mathscr{D}\right]&=& \J_{\sigma^2\eta}^{-1}  \Zm^\Tran\tv \label{eq:LR-theta-MAP-App},
\end{eqnarray}
where $\J_{\eta}  =\J+\eta\I_d$.  For $\eta=0$ the above gives us the  \emph{maximum likelihood} (ML) estimator 
\begin{eqnarray}
\thetaest\left[  \mathscr{D}\right]   &=& \J^{-1}  \Zm^\Tran\tv \label{eq:LR-theta-ML}
\end{eqnarray}
We  are interested   in the \emph{high-dimensional} regime: $(N,d)\rightarrow(\infty, \infty)$ with fixed $\zeta=d/N>0$, which we will write as $(N,d)\rightarrow \infty$ to simplify notation.  

\section{Distribution of $\thetaest$ estimator  in MAP inference \label{section:theta-est-distr}}
Let us  assume that the noise parameter $\sigma^2$ is \emph{independent} of  data $ \mathscr{D}$, and  that the noise $\epsilonv$ is  sampled from the Gaussian distribution $\mathcal{N}(\nullv, \sigma_0^2\I_N)$. The distribution of the MAP estimator  (\ref{eq:LR-theta-MAP-App}) can then be computed  as follows
\begin{eqnarray}
P(\thetaest)&=&\left\langle\delta\big(\thetaest -   \J^{-1}_{\sigma^2\eta}  \Zm^\Tran\tv    \big) \right\rangle_{ \mathscr{D}}  =\left\langle\delta\big(\thetaest -   \J^{-1}_{\sigma^2\eta}  \Zm^\Tran \big( \Zm\thetav_{\!0}+\epsilonv \big)    \big) \right\rangle_{ \mathscr{D}} \nonumber\\
&=&\left\langle\left\langle  \delta\big(\thetaest -   \J^{-1}_{\sigma^2\eta}  \J\thetav_{\!0}-   \J^{-1}_{\sigma^2\eta}  \Zm^\Tran   \epsilonv    \big)   \right\rangle_{\epsilonv}\right\rangle_{\Zm}    \nonumber\\
&=&  \int\frac{\rmd \x}{(2\pi)^d} \rme^{\rmi\x^T\thetaest}\left\langle \rme^{-\rmi\x^T    \J^{-1}_{\sigma^2\eta}  \J\thetav_{\!0}       }\left\langle \rme^{-\rmi\x^T   \J^{-1}_{\sigma^2\eta}  \Zm^\Tran   \epsilonv       }\right\rangle_{\epsilonv}\right\rangle_{\Zm} \nonumber\\
&=&  \int\frac{\rmd \x}{(2\pi)^d} \left\langle \rme^{-\frac{1}{2}\sigma_0^2\x^T \J^{-1}_{\sigma^2\eta}\J^{-1}_{\sigma^2\eta}   \J \x +\rmi\x^T \big(\thetaest - \J^{-1}_{\sigma^2\eta}  \J\thetav_{\!0}     \big)      }\right\rangle_{\Zm} \nonumber\\
&=&\frac{1}{(2\pi)^d}  \Bigg\langle  \sqrt{(2\pi)^d\left\vert \big(\sigma_0^2\J^{-2}_{\sigma^2\eta}   \J\big)^{-1}  \right\vert}\int \!\rmd \x~
\mathcal{N}\big(\x\vert\nullv, \big(\sigma_0^2\J^{-2}_{\sigma^2\eta}   \J\big)^{-1}  \big)\rme^{\rmi\x^T\big(\thetaest-    \J^{-1}_{\sigma^2\eta}  \J\thetav_{\!0}\big)    }\Bigg\rangle_{\Zm}\nonumber\\
&=& \left\langle \mathcal{N}\Big(\thetaest\, \vert \,  \J^{-1}_{\sigma^2\eta}  \J\thetav_{\!0} ,\,  \sigma_0^2\J^{-2}_{\sigma^2\eta}   \J    \Big)\right\rangle_{\Zm}
\label{eq:LR-MAP-theta-distr.-1}
\end{eqnarray}
To take the limit $(N,d)\rightarrow \infty$,  we rescale $z_i(\mu)\to z_i(\mu)/\sqrt{d}$ with $z_i(\mu)={\mathcal O}(1)$. 
Now $\J=\C/\zeta$,  where $\zeta=d/N$, $[\C]_{\mu\nu}=N^{-1}\sum_i z_i(\mu)z_i(\nu)$, and 
$\J^{-1}_{\sigma^2\eta} = \zeta\C^{-1}_{\zeta\sigma^2\eta}$ giving us the distribution
\begin{eqnarray}
P(\thetaest)
&=& \left\langle \mathcal{N}\Big(\thetaest\,\vert \, \C^{-1}_{\zeta\sigma^2\eta}\C\thetav_{\!0},\, \zeta \sigma_0^2 \C^{-2}_{\zeta\sigma^2\eta}\C\Big)\right\rangle_{\Zm} \label{eq:LR-MAP-theta-distr.-3}
\end{eqnarray}
with $\C_{\zeta\sigma^2\eta}\equiv \C+\zeta\sigma^2\eta\I$.
Furthermore,  since  $\C\equiv\C[\Zm]$ and $\C^{-1}_{\zeta\sigma^2\eta}\equiv \C^{-1}_{\zeta\sigma^2\eta}[\C[\Zm]]$ we have that 
\begin{eqnarray}
P(\thetaest)&=& \left\langle \mathcal{N}\Big(\thetaest\,\vert \, \C^{-1}_{\zeta\sigma^2\eta}[\C[\Zm]]\,\C[\Zm]\thetav_{\!0},\, \zeta\sigma_0^2  \C^{-2}_{\zeta\sigma^2\eta}[\C[\Zm]]\,\C[\Zm]   \Big)\right\rangle_{\Zm} \label{eq:LR-MAP-theta-distr.-4}\\
&=&\int\rmd\C\left\{ \prod_{i=1}^N\int P(\z_i)\rmd\z_i \right\} \delta\big(\C-\C[\Zm]\big)  \mathcal{N}\Big(\thetaest\,\vert \, \C^{-1}_{\zeta\sigma^2\eta}[\C]\,\C\thetav_{\!0},\, \zeta\sigma_0^2  \C^{-2}_{\zeta\sigma^2\eta}[\C]\,\C  \Big)         \nonumber\\
&=&\int\rmd\C \int\frac{\rmd\hat{\cmatrix}}{\big(2\pi\big)^{d^2}}\left\{ \prod_{i=1}^N\int P(\z_i)\rmd\z_i \right\} 
\exp\left[\rmi\mathrm{Tr}\left\{\hat{\cmatrix}\big(\J-\J[\Zm]\big)\right\}\right]
  \mathcal{N}\Big(\thetaest\,\vert \, \C^{-1}_{\zeta\sigma^2\eta}[\C]\,\C\thetav_{\!0},\, \zeta\sigma_0^2 
  \C^{-2}_{\zeta\sigma^2\eta}[\C]\,\C  \Big)    \nonumber
\end{eqnarray}
Let us now  consider the following average, assuming  that $\z$ is sampled from the Gaussian distribution $\mathcal{N}(\nullv,\Cov)$:
\begin{eqnarray}
\left\{ \prod_{i=1}^N\int P(\z_i)\rmd\z_i \right\}               \exp\left[-\rmi\mathrm{Tr}\left\{\hat{\cmatrix}\J[\Zm]\right\}\right]
&=& \left\{ \prod_{i=1}^N\int P(\z_i)\rmd\z_i \right\}               \exp\left[ -\rmi\frac{1}{N} \mathrm{Tr}\left\{ \hat{\cmatrix} \sum_{i=1}^N{\z_i\z_i}^{\!\!\mathrm{T}}\right\}\right]\nonumber\\
&=& \left[\int P(\z)\,\rmd\z\, \rme^{ -\rmi\frac{1}{N} \mathrm{Tr}\left\{ \hat{\cmatrix}\, {\z\z}^\Tran\right\}}\right]^N =
\left\vert\I+\frac{2\rmi }{N}\Cov\hat{\cmatrix} \right\vert^{-\frac{N}{2}}. \label{eq:LR-MAP-Wishart-GF}
\end{eqnarray}
This is the characteristic function of the Wishart distribution~\cite{Eaton1983},  defined by the density
 \begin{eqnarray}
\mathcal{W}\big(\J\vert\Cov/N,d,N\big)&=&\frac{\left\vert \Cov/N\right\vert^{-\frac{N}{2}}  \left\vert \J\right\vert^{\frac{N-d-1}{2}}}{2^{\frac{N d}{2}} \pi^{\frac{d(d-1)}{4}}\prod_{\ell=1}^d\Gamma\big(\frac{N+1-\ell}{2}\big)}\rme^{-N\frac{1}{2}\Tr\big(\J{\Cov}^{-1}\big)}\label{def:Wishart}. 
\end{eqnarray}
The Wishart distribution is \emph{singular} when $d>N$.  Thus for  Gaussian $\z$  the distribution    (\ref{eq:LR-MAP-theta-distr.-1})  is  the Gaussian mixture
\begin{eqnarray}
P(\thetaest)&=& \int\!\!\rmd\C~ \mathcal{W}\!\big(\C\vert\Cov/N,d,N\big)   \mathcal{N}\Big(\thetaest\,\vert \, \C^{-1}_{\zeta\sigma^2\eta}[\C]\,\C\thetav_{\!0},\, \zeta\sigma_0^2  \C^{-2}_{\zeta\sigma^2\eta}[\C]\,\C  \Big)  \nonumber \\
&=& \int  \rmd\C~\frac{\left\vert \Cov/N\right\vert^{-\frac{N}{2}}  \left\vert \C\right\vert^{\frac{N-d-1}{2}}\rme^{-N\frac{1}{2}\Tr\big(\C{\Cov}^{-1}\big)}}{2^{\frac{N d}{2}} \pi^{\frac{d(d-1)}{4}}\prod_{\ell=1}^d\Gamma\big(\frac{N+1-\ell}{2}\big)}
\frac{\rme^{-\frac{1}{2\zeta\sigma_0^2} \big( \thetaest  - \C^{-1}_{\zeta\sigma^2\eta}\C\thetav_{\!0}   \big)^\Tran \big(\C^{-2}_{\zeta\sigma^2\eta}[\C]\,\C\big)^{-1}[\C] \big( \thetaest  - \C^{-1}_{\zeta\sigma^2\eta}\C\thetav_{\!0}   \big)   }}{  \left\vert 2\pi \zeta\sigma_0^2\C^{-2}_{\zeta\sigma^2\eta}[\C]\,\C\right\vert^{\frac{1}{2}}     } .
 \label{eq:LR-MAP-theta-distr.-z-Norm}
 \end{eqnarray}
We note that an alternative derivation of this result is provided in~\cite{Coolen2020}.

\section{Distribution of $\thetaest$ estimator in ML inference\label{section:theta-est-ML}}

Let us consider the following integral appearing in the distribution of MAP estimator (\ref{eq:LR-MAP-theta-distr.-z-Norm}):
\begin{eqnarray}
&&\hspace*{-10mm} \int  \rmd\C\,  \left\vert \C\right\vert^{\frac{N-d-1}{2}}   \rme^{-N\frac{1}{2}\Tr\big(\C{\Cov}^{-1}\big)}\frac{\rme^{-\frac{1}{2\zeta\sigma_0^2}\big( \thetaest  - \C^{-1}_{\zeta\sigma^2\eta}\C\thetav_{\!0}   \big)^\Tran \big(\C^{-2}_{\zeta\sigma^2\eta}[\C]\,\C\big)^{-1}[\C] \big( \thetaest  - \C^{-1}_{\zeta\sigma^2\eta}\C\thetav_{\!0}   \big)   }}{   \left\vert 2\pi \zeta\sigma_0^2 \C^{-2}_{\zeta\sigma^2\eta}[\C]\,\C\right\vert^{\frac{1}{2}}       } \nonumber\\
&=& \int  \rmd\C\,   \left\vert \C\right\vert^{\frac{N-d-1}{2}} \frac{\rme^{-\frac{1}{2}\Tr\left\{\C\left[ N{\Cov}^{-1}  + \C^{-1}\frac{1}{\zeta\sigma_0^2}\big( \C^{-2}_{\zeta\sigma^2\eta}[\C]\,\C\big)^{-1}[\C]\big( \thetaest  - \C^{-1}_{\zeta\sigma^2\eta}\C\thetav_{\!0}   \big) \big( \thetaest  - \C^{-1}_{\zeta\sigma^2\eta}\C\thetav_{\!0}   \big)^{\mathrm{T}  }  \right] \right\} }}{   \left\vert 2\pi \zeta\sigma_0^2 \C^{-2}_{\zeta\sigma^2\eta}[\C]\,\C \right\vert^{\frac{1}{2}}    }  \nonumber
\label{eq:LR-MAP-J-int.}
\end{eqnarray}
 For $\eta=0$, i.e. ML inference, this integral simplifies to 
\begin{eqnarray}
&&\int  \rmd\C\,   \left\vert \C\right\vert^{\frac{N-d-1}{2}}  \frac{\rme^{-\frac{1}{2}\Tr\left\{\C\left[ N{\Cov}^{-1}  + \frac{1}{\zeta\sigma_0^2}\big( \thetaest  - \thetav_{\!0}   \big) \big( \thetaest  - \thetav_{\!0}   \big)^{\mathrm{T}  }  \right] \right\} }}{    \left\vert 2\pi \zeta\sigma_0^2\C^{-1}\right\vert^{\frac{1}{2}}    }  
= \big(\frac{1}{ 2\pi \zeta \sigma_0^2} \big)^{\frac{d}{2}} \int  \rmd\C\,   \left\vert \C\right\vert^{\frac{N-d}{2}}  \rme^{-\frac{1}{2}\Tr\left\{\C\left[ N{\Cov}^{-1}  + \frac{1}{\zeta\sigma_0^2}\big( \thetaest  - \thetav_{\!0}   \big) \big( \thetaest  - \thetav_{\!0}   \big)^{\mathrm{T}  }  \right] \right\} }  \nonumber
\end{eqnarray}
and can be computed by using  the normalization identity
$\int  \rmd \J\, \mathcal{W}\big(\J\vert\Cov,d,N\big)=1$, 
from which one obtains 
\begin{eqnarray}
\int  \rmd \C\,    \left\vert \C\right\vert^{\frac{N-d-1}{2}} \rme^{-\frac{1}{2}\Tr\big(\C{\Cov}^{-1}\big)}=\left\vert \Cov\right\vert^{\frac{N}{2}} 2^{\frac{N d}{2}} \pi^{\frac{d(d-1)}{4}}\prod_{\ell=1}^d\Gamma\big(\frac{N+1-\ell}{2}\big)\label{eq:Wishart-identity}.
\end{eqnarray}
Hence also 
\begin{eqnarray}
\int  \rmd \C\,    \left\vert \C\right\vert^{\frac{N-d}{2}} \rme^{-\frac{1}{2}\Tr\big(\C{\Cov}^{-1}\big)}=\left\vert \Cov\right\vert^{\frac{N+1}{2}} 2^{\frac{(N+1) d}{2}} \pi^{\frac{d(d-1)}{4}}\prod_{\ell=1}^d\Gamma\big(\frac{N+2-\ell}{2}\big)\label{eq:Wishart-identity-N+1}
\end{eqnarray}
which 
gives us the result
\begin{eqnarray}
\int  \rmd\C\,   \left\vert \C\right\vert^{\frac{N-d}{2}}  \rme^{-\frac{1}{2}\Tr\left\{\C\left[ N{\Cov}^{-1}  + \frac{1}{\zeta\sigma_0^2}\big( \thetaest - \thetav_{\!0}   \big) \big( \thetaest  - \thetav_{\!0}   \big)^{\mathrm{T}  }  \right] \right\} }  
=\left\vert N{\Cov}^{-1} \! +\! \frac{( \thetaest \! -\! \thetav_{\!0} ) ( \thetaest \! -\! \thetav_{\!0} )^{\mathrm{T} }}{\zeta\sigma_0^2}
 \right\vert^{-\frac{N+1}{2}} \!\!\!
 2^{\frac{(N+1) d}{2}} \pi^{\frac{d(d-1)}{4}}\prod_{\ell=1}^d\Gamma\big(\frac{N\!+\!2\!-\!\ell}{2}\big).
\nonumber
\\
\end{eqnarray}
 For $\eta=0$ the distribution (\ref{eq:LR-MAP-theta-distr.-z-Norm}) thereby becomes
\begin{eqnarray}
P(\thetaest)&=& \int  \rmd\C\,\frac{\left\vert \Cov/N\right\vert^{-\frac{N}{2}}  \left\vert \C\right\vert^{\frac{N-d-1}{2}}}{2^{\frac{N d}{2}} \pi^{\frac{d(d-1)}{4}}\prod_{\ell=1}^d\Gamma\big(\frac{N+1-\ell}{2}\big)}\rme^{-N\frac{1}{2}\Tr\big(\C{\Cov}^{-1}\big)}\frac{\rme^{-\frac{1}{2}\frac{1}{\zeta\sigma_0^2}\big( \thetaest  - \thetav_{\!0}   \big)^\Tran \C \big( \thetaest  - \thetav_{\!0}   \big)   }}{  \left\vert 2\pi \zeta\sigma_0^2\C^{-1}\right\vert^{\frac{1}{2}}     } .\nonumber\\
&=&\frac{\left\vert \Cov/N\right\vert^{-\frac{N}{2}}}{   2^{\frac{N d}{2}} \pi^{\frac{d(d-1)}{4}}\prod_{\ell=1}^d\Gamma\big(\frac{N+1-\ell}{2}\big)}\int  \rmd\C\, \left\vert \C\right\vert^{\frac{N-d-1}{2}}  \rme^{-\frac{1}{2}\Tr\big(\C N{\Cov}^{-1}\big)    } \frac{\rme^{-\frac{1}{2}\frac{1}{\zeta\sigma_0^2}\big( \thetaest  - \thetav_{\!0}   \big)^\Tran \C \big( \thetaest  - \thetav_{\!0}   \big)   }}{  \left\vert 2\pi \zeta\sigma_0^2\C^{-1}\right\vert^{\frac{1}{2}}     } .\nonumber\\
&=&  \big(\frac{1}{ \pi \zeta \sigma_0^2} \big)^{\frac{d}{2}}\prod_{\ell=1}^d\frac{  \Gamma\big(\frac{N+2-\ell}{2}\big)             }{  \Gamma\big(\frac{N+1-\ell}{2}\big)}\left\vert \Cov/N\right\vert^{-\frac{N}{2}}  \left\vert N{\Cov}^{-1}  + \frac{1}{\zeta\sigma_0^2}\big( \thetaest  - \thetav_{\!0}   \big) \big( \thetaest  - \thetav_{\!0}   \big)^{\mathrm{T}  } \right\vert^{-\frac{N+1}{2}} \nonumber\\
&=&  \pi^{-\frac{d}{2}}\frac{  \Gamma\big(\frac{N+1}{2}\big)             }{  \Gamma\big(\frac{N+1-d}{2}\big)} \left\vert  \frac{\Cov}{\zeta \sigma_0^2N}  \right\vert^{\frac{1}{2}}      \big( 1  + \big( \thetaest  - \thetav_{\!0}   \big)^\Tran    \frac{\Cov}{\zeta \sigma_0^2N}\big( \thetaest  - \thetav_{\!0}   \big) \big)^{-\frac{N+1}{2}}.
 \label{eq:LR-ML-theta-distr.-calc.}%
\end{eqnarray}
The last line in above was obtained  using  the `matrix determinant lemma'. Thus, after slight rearrangement,
\begin{eqnarray}
&&P(\thetaest)= [\pi(  N\!+\!1\!-\!d )]^{-\frac{d}{2}}\frac{  \Gamma\big(\frac{N+1}{2}\big)             }{  \Gamma\big(\frac{N+1-d}{2}\big)}
 \left\vert  \frac{(1\!-\!\zeta\!+\!1/N)\Cov}{\zeta \sigma_0^2}  \right\vert^{\frac{1}{2}} \big( 1  + ( \thetaest  \!-\! \thetav_{\!0} )^\Tran   \frac{(1-\zeta+1/N)\Cov}{\zeta\sigma_0^2(N\!+\!1\!-\!d)}( \thetaest \! -\! \thetav_{\!0} ) \big)^{-\frac{N+1}{2}},    \label{eq:LR-ML-theta-distr.}
\end{eqnarray}
which is the multivariate Student's $t$-distribution, with $N+1-d$ degrees of freedom, `location'  vector $\thetav_{\!0}$ and `shape' matrix $ \zeta\sigma_0^2\Cov^{-1}/(1-\zeta+1/N)$.

\section{Statistical properties of $\hat{\sigma}^2$ estimator in ML inference\label{section:sigma-est}}

In ML inference the estimator of $\thetav$  is given by (\ref{eq:LR-theta-ML}) and the estimator of noise parameter $\sigma^2$ is given by the density 
\begin{eqnarray}
\hat{\sigma}^2 [\mathscr{D}]&=&\frac{1}{N}\left\vert\left\vert\tv-\Zm\thetaest   [\mathscr{D}]   \right\vert\right\vert^2= \frac{1}{N}\left\vert\left\vert\tv-\Zm\,  \big(\Zm^\Tran\Zm\big)^{-1}  \Zm^\Tran\tv  \right\vert\right\vert^2=\frac{1}{N}\left\vert\left\vert\big(\I_N-\Zm\,     \big(\Zm^\Tran\Zm\big)^{-1}   \Zm^\Tran\big)\tv  \right\vert\right\vert^2 \nonumber\\
&=& \frac{1}{N}\left\vert\left\vert\big(\I_N-\Zm\,  \big(\Zm^\Tran\Zm\big)^{-1}     \Zm^\Tran\big)\big(\Zm\thetav_{\!0}+\epsilonv \big)  \right\vert\right\vert^2        \nonumber\\
&=& \frac{1}{N}\left\vert\left\vert  \big(\I_N-\Zm\,  \big(\Zm^\Tran\Zm\big)^{-1}   \Zm^\Tran\big)\Zm\thetav_{\!0}+\big(\I_N-\Zm\, \J^{-1}[\Zm]\Zm^\Tran\big)\epsilonv  \right\vert\right\vert^2   \nonumber\\  
&=& \frac{1}{N}\left\vert\left\vert  \big(\I_N\!-\!\Zm\,  \big(\Zm^\Tran\Zm\big)^{-1}    \Zm^\Tran\big)\epsilonv  \right\vert\right\vert^2
=\frac{1}{N} \epsilonv^\Tran \big(\I_N\!-\!\Zm\,    \big(\Zm^\Tran\Zm\big)^{-1}   \Zm^\Tran  \big)^2\epsilonv  
=\frac{1}{N} \epsilonv^\Tran \big(\I_N\!-\!\Zm\,    \big(\Zm^\Tran\Zm\big)^{-1}   \Zm^\Tran  \big)\epsilonv
~~~~
 \label{eq:LR-ML-noise-estimator}
\end{eqnarray}
In above we used $(\I_N\!-\!\Zm\,  (\Zm^\Tran\Zm)^{-1}   \Zm^\Tran)\Zm\thetav_{\!0}= \Zm\thetav_{\!0}-\Zm\,  \big(\Zm^\Tran\Zm\big)^{-1}   \Zm^\Tran\Zm\thetav_{\!0}=\nullv$ and $(\I_N\!-\!\Zm\,    \big(\Zm^\Tran\Zm\big)^{-1}   \Zm^\Tran )^2=\I_N\!-\!\Zm\,    \big(\Zm^\Tran\Zm\big)^{-1}   \Zm^\Tran $, i.e. $\I_N\!-\!\Zm\,    \big(\Zm^\Tran\Zm\big)^{-1}   \Zm^\Tran$ is an \emph{idempotent} matrix.  For $\epsilonv$ sampled from \emph{any}  distribution with mean $\nullv$ and covariance $\sigma_0^2\I_N$,  the average and variance of $\hat{\sigma}^2 [\mathscr{D}]$ are (by Wick's theorem):
\begin{eqnarray}
\left\langle\hat{\sigma}^2 [\mathscr{D}]\right\rangle_{\epsilonv} &=&\frac{1}{N}\left\langle \epsilonv^\Tran \big(\I_N-\Zm\,    \big(\Zm^\Tran\Zm\big)^{-1}   \Zm^\Tran  \big)\epsilonv\right\rangle_{\epsilonv}
=\frac{\sigma_0^2}{N}\Tr \big(\I_N-\Zm\,    \big(\Zm^\Tran\Zm\big)^{-1}   \Zm^\Tran  \big)=\sigma_0^2\big(1-\zeta\big)\label{eq:LR-ML-noise-estimator-aver.}
\\
\left\langle\hat{\sigma}^4 [\mathscr{D}]\right\rangle_{\epsilonv}-\left\langle\hat{\sigma}^2 [\mathscr{D}]\right\rangle^2_{\epsilonv} &=& \frac{2\sigma_0^4}{N}\big(1-\zeta\big)   \label{eq:LR-ML-noise-estimator-var.}
\end{eqnarray}
Next 
we are interested in the probability of event $\hat{\sigma}^2 [\mathscr{D}]  \notin \big(\sigma_0^2(1-\zeta)-\delta, \sigma_0^2(1-\zeta)+\delta\big)$. This is given by 
\begin{eqnarray}
\mathrm{Prob} \left[   \frac{1}{N}\sum_{i=1}^N\big(t_i -\thetaest [\mathscr{D}].\z_i\big)^2\notin \big(\sigma_0^2(1-\zeta)-\delta, \sigma_0^2(1-\zeta)+\delta\big)     \right]
&=& \mathrm{Prob} \left[\sum_{i=1}^N\big(t_i -\thetaest [\mathscr{D}].\z_i\big)^2\leq N\big(\sigma_0^2(1-\zeta)-\delta\big)\right]  \nonumber\\
&&\hspace*{-10mm}+ ~\mathrm{Prob} \left[\sum_{i=1}^N\big(t_i -\thetaest [\mathscr{D}].\z_i\big)^2\geq N\big(\sigma_0^2(1-\zeta)+\delta\big)\right]  \label{eq:LR-ML-noise-estimator-tails} .
%
\end{eqnarray}
First,  we  consider the probability 
\begin{eqnarray}
\mathrm{Prob} \left[\sum_{i=1}^N\big(t_i -\thetaest [\mathscr{D}].\z_i\big)^2\geq N\big(\sigma_0^2(1-\zeta)+\delta\big)\right]
&=&\mathrm{Prob} \left[ \rme^{\frac{1}{2}\alpha\sum_{i=1}^N\big(t_i -\thetaest [\mathscr{D}].\z_i\big)^2}\geq \rme^{\frac{1}{2}\alpha N\big(\sigma_0^2(1-\zeta)+\delta\big)}\right]\nonumber\\
&& \leq\left\langle\rme^{\frac{1}{2}\alpha\sum_{i=1}^N\big(t_i -\thetaest [\mathscr{D}].\z_i\big)^2}\right\rangle_{\mathscr{D}} \rme^{-\frac{1}{2}\alpha N\big(\sigma_0^2(1-\zeta)+\delta\big)}, 
\end{eqnarray}
where i$\alpha>0$ and  we used Markov inequality to derive the upper bound.  
Let us assume that the distribution of noise is \emph{Gaussian} and consider the moment-generating  function 
\begin{eqnarray}
 \left\langle\rme^{\frac{1}{2}\alpha\sum_{i=1}^N\big(t_i -\thetaest[\mathscr{D}].\z_i\big)^2}\right\rangle_{\mathscr{D}} 
&=&  \left\langle\rme^{\frac{1}{2}\alpha\epsilonv^\Tran \big(\I_N-\Zm\,    \big(\Zm^\Tran\Zm\big)^{-1}   \Zm^\Tran  \big)\epsilonv             }\right\rangle_{\mathscr{D}}     =  \int\rmd\epsilonv\,\mathcal{N}\!\big(\epsilonv \vert\nullv,\sigma_0^2\I_N\big) \left\langle\rme^{\frac{1}{2}\alpha\epsilonv^\Tran \big(\I_N-\Zm\,    \big(\Zm^\Tran\Zm\big)^{-1}   \Zm^\Tran  \big)\epsilonv             }\right\rangle_{\Zm}   \nonumber\\
&=&  \frac{1}{\big(2\pi \sigma_0^2 \big)^{N/2} }  \left\langle \int\rmd\epsilonv\, \rme^{ -\frac{1}{2\sigma_0^2} \epsilonv^\Tran    \epsilonv + \frac{1}{2}\alpha\epsilonv^\Tran \big(\I_N-\Zm\,    \big(\Zm^\Tran\Zm\big)^{-1}   \Zm^\Tran  \big)\epsilonv             }\right\rangle_{\Zm}   \nonumber\\
&=&  \frac{1}{\big(2\pi \sigma_0^2 \big)^{N/2} }  \left\langle \int\rmd\epsilonv\, \rme^{ -\frac{1}{2} \epsilonv^\Tran\left[ \I_N/\sigma_0^2   - \alpha \big(\I_N-\Zm\,    \big(\Zm^\Tran\Zm\big)^{-1}   \Zm^\Tran  \big)\right]\epsilonv             }\right\rangle_{\Zm}   \nonumber\\
&=&  \frac{1}{\big(2\pi \sigma_0^2 \big)^{N/2} }  \left\langle\left\vert2\pi \left[ \I_N/\sigma_0^2   - \alpha \big(\I_N-\Zm\,    \big(\Zm^\Tran\Zm\big)^{-1}   \Zm^\Tran  \big)\right]^{-1}\right\vert^{\frac{1}{2}}\right\rangle_{\Zm}   \nonumber\\
&=&   \left\langle\left\vert  \I_N  - \alpha\sigma_0^2  \big(\I_N-\Zm\,    \big(\Zm^\Tran\Zm\big)^{-1}   \Zm^\Tran  \big)\right\vert^{-\frac{1}{2}}\right\rangle_{\Zm}   
=   \left\langle\left\vert  \I_N (1 - \alpha\sigma_0^2) +   \alpha\sigma_0^2\Zm\,    \big(\Zm^\Tran\Zm\big)^{-1}   \Zm^\Tran  \right\vert^{-\frac{1}{2}}\right\rangle_{\Zm}   \nonumber\\
&=&   \big(1 - \alpha\sigma_0^2\big)^{-N/2}  \left\langle  \frac{1}{\left\vert  \I_N  +  \frac{ \alpha\sigma_0^2}{1 - \alpha\sigma_0^2}\Zm\,    \big(\Zm^\Tran\Zm\big)^{-1}   \Zm^\Tran  \right\vert^{\frac{1}{2}}}\right\rangle_{\Zm} 
\end{eqnarray}
Now $ \Zm\,    \big(\Zm^\Tran\Zm\big)^{-1}   \Zm^\Tran$ is a \emph{projection} matrix, and  its eigenvalue are $\lambda_i\in\{0,1\}$, giving us 
\begin{eqnarray}
 \left\langle\rme^{\frac{1}{2}\alpha\sum_{i=1}^N\big(t_i -\thetaest[\mathscr{D}].\z_i\big)^2}\right\rangle_{\mathscr{D}}  
&=&   \big(1 - \alpha\sigma_0^2\big)^{-N/2}  \left\langle  \frac{1}{\left\vert  \I_N  +  \frac{ \alpha\sigma_0^2}{1 - \alpha\sigma_0^2}\Zm\,    \big(\Zm^\Tran\Zm\big)^{-1}   \Zm^\Tran  \right\vert^{\frac{1}{2}}}\right\rangle_{\Zm}\nonumber\\
&&=  \big(1 - \alpha\sigma_0^2\big)^{-N/2}  \left\langle  \frac{1}{ \big(   \prod_{i=1}^N \left[1+  \frac{ \alpha\sigma_0^2}{1 - \alpha\sigma_0^2}\lambda_i\big( \Zm\,    \big(\Zm^\Tran\Zm\big)^{-1}   \Zm^\Tran\big) \right]  \big)^{\frac{1}{2}}}\right\rangle_{\Zm}\nonumber\\
&&=  \big(1 - \alpha\sigma_0^2\big)^{-N/2}  \left\langle \rme^{-\frac{N}{2}  \sum_{\lambda} \frac{1}{N}\sum_{i=1}^N\delta_{\lambda; \lambda_i\big( \Zm\,    \big(\Zm^\Tran\Zm\big)^{-1}   \Zm^\Tran\big)}\log   \left[1+  \frac{ \alpha\sigma_0^2}{1 - \alpha\sigma_0^2}\lambda\right]}\right\rangle_{\Zm}\nonumber\\
&&=  \big(1 - \alpha\sigma_0^2\big)^{-N/2}  \left\langle \rme^{-\frac{N}{2}\log   \left[1+  \frac{ \alpha\sigma_0^2}{1 - \alpha\sigma_0^2}\right]  \frac{1}{N}\sum_{i=1}^N \lambda_i\big( \Zm\,    \big(\Zm^\Tran\Zm\big)^{-1}   \Zm^\Tran\big)}\right\rangle_{\Zm}\nonumber\\
 &&=  \big(1 - \alpha\sigma_0^2\big)^{-N/2}  \left\langle \rme^{-\frac{1}{2}\log   \big(1+  \frac{ \alpha\sigma_0^2}{1 - \alpha\sigma_0^2}\big)  \Tr\left\{ \Zm\,    \big(\Zm^\Tran\Zm\big)^{-1}   \Zm^\Tran\right\}}\right\rangle_{\Zm}\nonumber\\
 &&=  \big(1 - \alpha\sigma_0^2\big)^{-N/2}  \rme^{-\frac{d}{2}\log   \big(1+  \frac{ \alpha\sigma_0^2}{1 - \alpha\sigma_0^2}\big) }\nonumber\\
 &&= \rme^{-\frac{N}{2}\left[ \zeta \log  \big(1+  \frac{ \alpha\sigma_0^2}{1 - \alpha\sigma_0^2}\big)  + \log  \big(1 - \alpha\sigma_0^2\big)\right]} = \rme^{-\frac{N}{2}\big(1- \zeta\big)  \log  \big(1 - \alpha\sigma_0^2\big)}. \label{eq:LR-ML-noise-estimator-MGF}  
\end{eqnarray}
Hence 
\begin{eqnarray}
\mathrm{Prob} \left[\sum_{i=1}^N\big(t_i -\thetaest[\mathscr{D}].\z_i\big)^2\geq N\big(\sigma_0^2(1-\zeta)+\delta\big)\right]
&\leq &  \left\langle\rme^{\frac{1}{2}\alpha\sum_{i=1}^N\big(t_i -\thetaest[\mathscr{D}].\z_i\big)^2}\right\rangle_{\mathscr{D}}  \rme^{-\frac{1}{2}\alpha N\big(\sigma_0^2(1-\zeta)+\delta\big)}\nonumber\\
&\leq&   \rme^{-\frac{1}{2}N\big(1- \zeta\big)  \log  \big(1 - \alpha\sigma_0^2\big)}    \rme^{-\frac{1}{2}\alpha N\big(\sigma_0^2(1-\zeta)+\delta\big)}\nonumber\\
&=& \rme^{ -\frac{1}{2}N \left[  \big(1- \zeta\big)   \log \big(1 - \alpha\sigma_0^2\big)  +\alpha \big(\sigma_0^2(1-\zeta)+\delta\big)\right]}=\rme^{ -\frac{1}{2}N\Phi(\alpha)}.
\end{eqnarray}
The  \emph{rate} function 
\begin{eqnarray}
\Phi(\alpha)= \big(1- \zeta\big)\log \big(1 - \alpha\sigma_0^2\big)  +\alpha \big(\sigma_0^2(1-\zeta)+\delta\big)
\end{eqnarray}
 has a \emph{maximum} at   
\begin{eqnarray}
\alpha={\frac {\delta}{\sigma_0^2\, \big( (1-\zeta) \sigma_0^2+\delta\big) }},
\end{eqnarray}
such that  ${\rm max}_\alpha \Phi(\alpha)= {\big( 1-\zeta \big) \log 
 \big( {\frac {\big( 1-\zeta \big) }{ \big( 1-
\zeta \big) +\delta/\sigma_0^2}} \big) +\delta/\sigma_0^2} $, and hence 
\begin{eqnarray}
\mathrm{Prob} \left[\sum_{i=1}^N\big(t_i -\thetaest[\mathscr{D}].\z_i\big)^2\geq N\big(\sigma_0^2(1-\zeta)+\delta\big)\right]
\leq   \rme^{ -\frac{1}{2}N  \big(    {\big( 1-\zeta \big) \log 
 \big( {\frac {\big( 1-\zeta \big) }{ \big( 1-
\zeta \big) +\delta/\sigma_0^2}} \big) +\delta/\sigma_0^2} \big)} \label{eq:LR-ML-noise-estimator-upper-tail-ub}. 
\end{eqnarray}
We note that in above expression the rate function $( 1-\zeta) \log [
 ( 1-\zeta) /( 1-
\zeta+\delta/\sigma_0^2)] +\delta/\sigma_0^2 $ vanishes when $\delta/\sigma_0^2=0$, and is a monotonically increasing function of   $\delta/\sigma_0^2$. 

Second,  for $\alpha>0$ we consider the probability     
\begin{eqnarray}
\mathrm{Prob} \left[\sum_{i=1}^N\big(t_i -\thetaest[\mathscr{D}].\z_i\big)^2\leq N\big(\sigma_0^2(1-\zeta)-\delta\big)\right]&=&
\mathrm{Prob} \left[ \rme^{-\frac{1}{2}\alpha\sum_{i=1}^N\big(t_i -\thetaest[\mathscr{D}].\z_i\big)^2}\geq \rme^{-\frac{1}{2}\alpha N\big(\sigma_0^2(1-\zeta)-\delta\big)}\right]\nonumber\\
&\leq&   \left\langle\rme^{-\frac{1}{2}\alpha\sum_{i=1}^N\big(t_i -\thetaest[\mathscr{D}].\z_i\big)^2}\right\rangle_{\mathscr{D}}  \rme^{\frac{1}{2}\alpha N\big(\sigma_0^2(1-\zeta)-\delta\big)}\nonumber\\
&=& \rme^{-\frac{N}{2}\big(1- \zeta\big)  \log  \big(1 + \alpha\sigma_0^2\big)} \rme^{\frac{1}{2}\alpha N\big(\sigma_0^2(1-\zeta)-\delta\big)}= \rme^{-\frac{N}{2}\phi(\alpha)} ,
\end{eqnarray}
where 
\begin{eqnarray}
\phi(\alpha)=\big(1- \zeta\big)  \log  \big(1 + \alpha\sigma_0^2\big) -\alpha \big(\sigma_0^2(1-\zeta)-\delta\big).
\end{eqnarray}
Here we used the Markov inequality and the result (\ref{eq:LR-ML-noise-estimator-MGF}) with $\alpha\rightarrow-\alpha$.  The rate function $\phi(\alpha)$ has a maximum at  $\alpha=\delta/\sigma_0^2 ((1-\zeta) \sigma_0^2-\delta)$,  such that ${\rm max}_\alpha \phi(\alpha)=( 1-\zeta) \log[ 
 ( 1-\zeta)/
( 1-\zeta-\delta/\sigma_0^2)] -\delta/\sigma_0^2 $, and hence 
\begin{eqnarray}
\mathrm{Prob} \left[\sum_{i=1}^N\big(t_i -\thetaest[\mathscr{D}].\z_i\big)^2\leq N\big(\sigma_0^2(1-\zeta)-\delta\big)\right]\leq \rme^{-\frac{N}{2}  \big(  {  \big( 1-\zeta \big) \log 
 \big( {\frac {\big( 1-\zeta \big) }{
 \big( 1-\zeta \big)-\delta/\sigma_0^2 }} \big) -\delta/\sigma_0^2  }  \big)}\label{eq:LR-ML-noise-estimator-lower-tail-ub}. 
\end{eqnarray}
Here the rate function $( 1-\zeta) \log[ 
 ( 1-\zeta) /(1-\zeta-\delta/\sigma_0^2)] -\delta/\sigma_0^2 $  vanishes when $\delta/\sigma_0^2 =0$,  and is a monotonic increasing function of $\delta/\sigma_0^2$  when $\sigma_0^2\big(1-\zeta \big)>\delta$.

Finally, combining the inequalities (\ref{eq:LR-ML-noise-estimator-upper-tail-ub}) and (\ref{eq:LR-ML-noise-estimator-lower-tail-ub}) we obtain the inequality 
\begin{eqnarray}
\mathrm{Prob} \left[   \hat{\sigma}^2[\mathscr{D}]\notin \big(\sigma_0^2(1-\zeta)-\delta, \sigma_0^2(1-\zeta)+\delta\big)   \right]
&\leq & \rme^{-\frac{1}{2} N\left[  \big(  {  \big( 1-\zeta \big) \log 
 \big( {\frac {\big( 1-\zeta \big) }{
 \big( 1-\zeta \big)-\delta/\sigma_0^2 }} \big) -\delta/\sigma_0^2  }  \big)  \right]} \nonumber\\
   &&+ ~  \rme^{ -\frac{1}{2}N   \left[   \big(  {  \big( 1-\zeta \big) \log 
 \big( {\frac {\big( 1-\zeta \big) }{
 \big( 1-\zeta \big)+\delta/\sigma_0^2 }} \big) +\delta/\sigma_0^2  }  \big)  \right ]     } \label{eq:LR-ML-noise-estimator-tails-ub}
\end{eqnarray}
which is valid for   $\delta\in (0, \sigma_0^2 (1-\zeta))$.

\section{   Statistical properties  of MSE  in ML inference  \label{section:mse}}
In this section we consider statistical properties of  the \emph{minimum square error} (MSE)  $\frac{1}{d} \vert\vert\thetav_{\!0}- \thetaest [\mathscr{D}]    \vert\vert^2$, where $\thetav_0$ is the vector of the \emph{true}  parameters responsible for the data, and  $ \thetaest [\mathscr{D}] $ is the ML estimator (\ref{eq:LR-theta-ML}). 

\subsection{Moment generating function\label{section:mse-mgf}}
Let us consider the moment generating function 
\begin{eqnarray}
 \left\langle   \rme^{\frac{1}{2}\alpha  \vert\vert\thetav_{\!0}-\thetaest[\mathscr{D}]\vert\vert^2 }\right\rangle_{\mathscr{D}} &=&  \int\! \rmd  \thetaest ~ P(\thetaest)\,  \rme^{\frac{1}{2}\alpha  \vert\vert\thetav_{\!0}-\thetaest\vert\vert^2 }\,  \label{eq:LR-ML-MSE-MGF-calc.}\\
  &=& \big(\pi\big(  N+1-d  \big)\big)^{-\frac{d}{2}}\frac{  \Gamma\big(\frac{N+1}{2}\big)             }{  \Gamma\big(\frac{N+1-d}{2}\big)} \left\vert  \frac{(1-\zeta+1/N)\Cov}{\zeta \sigma_0^2}  \right\vert^{\frac{1}{2}} \nonumber\\
   &&~~~\times  \int\!\rmd  \thetaest~\big( 1  + \big( \thetaest  - \thetav_{\!0}   \big)^\Tran   \frac{1}{N+1-d} \frac{(1-\zeta+1/N)\Cov}{\zeta\sigma_0^2}\big( \thetaest  - \thetav_{\!0}   \big) \big)^{-\frac{N+1}{2}}   \rme^{\frac{1}{2}\alpha  \vert\vert\thetav_{\!0}-\thetaest\vert\vert^2 }\, \nonumber \\
   &=& \big(\pi\big(  N+1-d  \big)\big)^{-\frac{d}{2}}\frac{  \Gamma\big(\frac{N+1}{2}\big)             }{  \Gamma\big(\frac{N+1-d}{2}\big)} \left\vert  \frac{(1-\zeta+1/N)\Cov}{\zeta \sigma_0^2}  \right\vert^{\frac{1}{2}} \nonumber\\
   &&~~~\times  \int\!\rmd  \thetaest~\big( 1  + \thetaest^\Tran   \frac{1}{N+1-d} \frac{(1-\zeta+1/N)\Cov}{\zeta\sigma_0^2}\thetaest  \big)^{-\frac{N+1}{2}}  
    \rme^{\frac{1}{2}\alpha  \vert\vert\thetaest\vert\vert^2 }\, \nonumber \\
  &=& \int\!\rmd  \thetaest \int_{0}^\infty\!\rmd\omega~\mathcal{N}\big(\thetaest\,\vert\, \nullv,  \frac{\zeta\sigma_0^2   }{\omega(1-\zeta+1/N)}\Cov^{-1}\big) \Gamma_{\!N+1-d}\,(\omega) \rme^{\frac{1}{2}\alpha  \vert\vert\thetaest\vert\vert^2 }\, \nonumber \\
  &=& \int\!\rmd  \thetaest \int_{0}^\infty\!\rmd\omega ~\frac{\rme^{ -\frac{1}{2}  \thetaest^\Tran \big(   \frac{\omega(1-\zeta+1/N)\Cov}{\zeta\sigma_0^2  }-\alpha\I_d\big)\thetaest }}{\left\vert 2\pi\frac{\zeta\sigma_0^2   }{\omega(1-\zeta+1/N)}\Cov^{-1} \right\vert^{\frac{1}{2}}}\,   \,  \Gamma_{N+1-d}\,(\omega)\, \nonumber \\
  &=& \int_{0}^\infty\! \rmd\omega ~ \frac{  \Gamma_{N+1-d}\,(\omega)   }{\left\vert \frac{\zeta\sigma_0^2   }{\omega(1-\zeta+1/N)}\Cov^{-1} \right\vert^{\frac{1}{2}}   \left\vert \frac{\omega(1-\zeta+1/N)\Cov}{\zeta\sigma_0^2}-\alpha\I_d\right\vert^{\frac{1}{2}}}\ \nonumber \\
   &=& \int_{0}^\infty\!\rmd\omega~\frac{  \Gamma_{N+1-d}\,(\omega)   }{   \left\vert \I_d  -\alpha  \frac{\zeta\sigma_0^2   }{\omega(1-\zeta+1/N)}\Cov^{-1}  \right\vert^{\frac{1}{2}}}
   = \int_{0}^\infty  \!\rmd\omega~\rme^{-\frac{1}{2}     \log \left\vert \I_d  -\alpha  \frac{\zeta\sigma_0^2   }{\omega(1-\zeta+1/N)}\Cov^{-1}  \right\vert   } \Gamma_{N+1-d}\,(\omega) \,  \nonumber \\
       &=& \int_{0}^\infty \!\rmd\omega~\rme^{-\frac{1}{2}\sum_{\mu=1}^d     \log \big(1  -  \frac{\alpha\zeta\sigma_0^2   }{\omega(1-\zeta+1/N) \lambda_\mu(\Cov)}  \big)   } \Gamma_{N+1-d}\,(\omega) \,
\end{eqnarray}
where we encounter the \emph{gamma} distribution, for $\nu>0$,
\begin{eqnarray}
\Gamma_{\nu}\,(\omega)=\frac{\nu^{\nu/2}}{2^{\nu/2}\Gamma(\nu/2)}\omega^{\frac{\nu-2}{2}}\rme^{-\frac{1}{2}\nu\omega} 
\end{eqnarray}
We note that the above derivation was obtained using the mixture of Gaussians representation of multivariate Student t distribution~\cite{Taboga2017}.  Thus  the moment generating function is given by 
\begin{eqnarray}
\left\langle   \rme^{\frac{1}{2}\alpha  \vert\vert\thetav_{\!0}-\thetaest[\mathscr{D}]\vert\vert^2 }\right\rangle_{\mathscr{D}} &=& \int_{0}^\infty \frac{  \Gamma_{N+1-d}\,(\omega)   }{   \left\vert \I_d  -\alpha  \frac{\zeta\sigma_0^2   }{\omega(1-\zeta+1/N)}\Cov^{-1}  \right\vert^{\frac{1}{2}}}\,\rmd\omega   \label{eq:LR-ML-MSE-MGF}\\
&=& \int_{0}^\infty\!\rmd\omega~\Gamma_{N+1-d}\,(\omega)\prod_{\ell=1}^d  \big( 1  -\alpha  \frac{\zeta\sigma_0^2   }{\omega(1-\zeta+1/N) \lambda_\ell (\Cov)}  \big)^{-\frac{1}{2}}\nonumber
\end{eqnarray}
and, by the transformation $\alpha=2\rmi a$ in the above,  we also obtain the characteristic function 
\begin{eqnarray}
\left\langle   \rme^{\rmi a  \vert\vert\thetav_{\!0}-\thetaest\vert\vert^2 }\right\rangle_{\mathscr{D}}
&=& \int_{0}^\infty \Gamma_{N+1-d}\,(\omega)\,\rmd\omega\prod_{\ell=1}^d  \big( 1  -\rmi a  \frac{2\zeta\sigma_0^2   }{\omega(1-\zeta+1/N) \lambda_\ell (\Cov)}   \big)^{-\frac{1}{2}} \label{eq:LR-ML-MSE-CF}. 
\end{eqnarray}
We note that the last term  in above is the product of characteristic functions of gamma distributions. Each gamma distribution has the same `shape' parameter $1/2$ and  different scale parameter  $2\zeta\sigma_0^2  /\omega(1-\zeta+1/N) \lambda_\ell (\Cov)$.

\subsubsection{The first two moments of the MSE\label{section:mse-moments}}
Let us now consider derivatives of the moment generating function    (\ref{eq:LR-ML-MSE-MGF}) upon replacing $\alpha\to \alpha/d$. The derivative with respect to $\alpha$ then gives us 
\begin{eqnarray}
2\frac{\partial}{\partial\alpha}\left\langle   \rme^{\frac{1}{2d}\alpha  \vert\vert\thetav_{\!0}-\thetaest[\mathscr{D}]\vert\vert^2 }\right\rangle_{\mathscr{D}}  
&=&2 \int_{0}^\infty \!\rmd\omega~\Gamma_{N+1-d}\,(\omega)\,\frac{\partial}{\partial\alpha}\prod_{\ell=1}^d  \big( 1  -\alpha  \frac{\zeta\sigma_0^2   }{\omega d  (1-\zeta+1/N) \lambda_\ell (\Cov)}  \big)^{-\frac{1}{2}}\nonumber\\
 &=& \int_{0}^\infty\! \rmd\omega~ \Gamma_{N+1-d}\,(\omega) \prod_{\ell=1}^d  \big( 1  -\alpha  \frac{\zeta\sigma_0^2   }{\omega d  (1-\zeta+1/N) \lambda_\ell (\Cov)}  \big)^{-\frac{1}{2}}\nonumber\\
 &&~~~\times\frac{1}{d}\sum_{\ell=1}^d\!  \big(\! 1  -\alpha  \frac{\zeta\sigma_0^2   }{\omega d  (1-\zeta+1/N) \lambda_\ell (\Cov)} \! \big)^{-1}\!\!\! \!\frac{\zeta\sigma_0^2   }{\omega(1-\zeta+1/N) \lambda_\ell (\Cov)} \label{eq:LR-ML-MSE-MGF-1-der.}
\end{eqnarray}
For $\alpha=0$ this gives us the average 
\begin{eqnarray}
 \frac{1}{d}\left\langle   \vert\vert\thetav_{\!0}-\thetaest[\mathscr{D}]\vert\vert^2 \right\rangle_{\mathscr{D}} &=& \frac{\zeta\sigma_0^2   }{1-\zeta+1/N} \frac{1}{d}  \Tr \left\{\Cov^{-1}\right\}\int_{0}^\infty\!\rmd\omega~ \Gamma_{N+1-d}\,(\omega)\,\omega^{-1}\nonumber\\
 &=&\frac{\zeta\sigma_0^2   }{1-\zeta+1/N}  \frac{1-\zeta+1/N}{1-\zeta-1/N} \frac{1}{d}  \Tr \left[\Cov^{-1}\right] 
 =\frac{\zeta\sigma_0^2   }{1-\zeta-1/N}    \frac{1}{d}  \Tr \left[\Cov^{-1}\right]   \label{eq:LR-ML-MSE-mean}.
\end{eqnarray}
Now we consider the second derivative with respect to $\alpha$: 
\begin{eqnarray}
4\frac{\partial^2}{\partial\alpha^2}\left\langle   \rme^{\frac{1}{2d}\alpha  \vert\vert\thetav_{\!0}-\thetaest[\mathscr{D}]\vert\vert^2 }\right\rangle_{\mathscr{D}}   
&=&4 \int_{0}^\infty\!\rmd\omega~ \Gamma_{N+1-d}\,(\omega)\,\frac{\partial^2}{\partial\alpha^2}\prod_{\ell=1}^d  \big( 1  -\alpha  \frac{\zeta\sigma_0^2   }{\omega d(1-\zeta+1/N) \lambda_\ell (\Cov)}  \big)^{-\frac{1}{2}}\nonumber\\
 &=&2 \int_{0}^\infty\!\rmd\omega~\Gamma_{N+1-d}\,(\omega)\,\frac{\partial}{\partial\alpha}\prod_{\ell=1}^d  \big( 1  -\alpha  \frac{\zeta\sigma_0^2   }{\omega d(1-\zeta+1/N) \lambda_\ell (\Cov)}  \big)^{-\frac{1}{2}}\nonumber\\
 &&~~~\times\frac{1}{d}\sum_{\ell=1}^d  \big( 1  -\alpha  \frac{\zeta\sigma_0^2   }{\omega d(1-\zeta+1/N) \lambda_\ell (\Cov)}  \big)^{-1} \frac{\zeta\sigma_0^2   }{\omega(1-\zeta+1/N) \lambda_\ell(\Cov)} \nonumber\\
&=& \int_{0}^\infty\!\rmd\omega~ \Gamma_{N+1-d}\,(\omega)\,\Bigg\{2\frac{\partial}{\partial\alpha}\prod_{\ell=1}^d  \big( 1  -\alpha  \frac{\zeta\sigma_0^2   }{\omega d(1-\zeta+1/N) \lambda_\ell (\Cov)}  \big)^{-\frac{1}{2}}\nonumber\\
 &&~~~\times\frac{1}{d}\sum_{\ell=1}^d  \big( 1  -\alpha  \frac{\zeta\sigma_0^2   }{\omega d(1-\zeta+1/N) \lambda_\ell (\Cov)}  \big)^{-1} \frac{\zeta\sigma_0^2   }{\omega(1-\zeta+1/N) \lambda_\ell(\Cov)} \nonumber\\
  &&~~~+2    \prod_{\ell=1}^d  \big( 1  -\alpha  \frac{\zeta\sigma_0^2   }{\omega d(1-\zeta+1/N) \lambda_\ell (\Cov)}  \big)^{-\frac{1}{2}}\nonumber\\
 &&~~~\times\frac{\partial}{\partial\alpha}\frac{1}{d}\sum_{\ell=1}^d  \big( 1  -\alpha  \frac{\zeta\sigma_0^2   }{\omega d(1-\zeta+1/N) \lambda_\ell (\Cov)}  \big)^{-1} \frac{\zeta\sigma_0^2   }{\omega(1-\zeta+1/N) \lambda_\ell(\Cov)}\Bigg\} \nonumber\\
&=& \int_{0}^\infty\!\rmd\omega~ \Gamma_{N+1-d}\,(\omega)\, \prod_{\ell=1}^d  \big( 1  -\alpha  \frac{\zeta\sigma_0^2   }{\omega d(1-\zeta+1/N) \lambda_\ell (\Cov)}  \big)^{-\frac{1}{2}}\nonumber\\
 &&~~~\times\Bigg\{  \left[\frac{1}{d}\sum_{\ell=1}^d  \big( 1  -\alpha  \frac{\zeta\sigma_0^2   }{\omega d(1-\zeta+1/N) \lambda_\ell (\Cov)}  \big)^{-1} \frac{\zeta\sigma_0^2   }{\omega(1-\zeta+1/N) \lambda_\ell(\Cov)} \right]^2\nonumber\\
  &&~~~+\frac{2}{d^2}\sum_{\ell=1}^d  \big( 1  -\alpha  \frac{\zeta\sigma_0^2   }{\omega d(1-\zeta+1/N) \lambda_\ell (\Cov)}  \big)^{-2} \left[\frac{\zeta\sigma_0^2   }{\omega(1-\zeta+1/N) \lambda_\ell(\Cov)}  \right]^2 \Bigg\} \nonumber\\\label{eq:LR-ML-MSE-MGF-2-der.}
\end{eqnarray}
Evaluation at $\alpha=0$ gives us the second moment
\begin{eqnarray}
\left\langle \frac{1}{d^2}  \vert\vert\thetav_{\!0}-\thetaest [\mathscr{D}] \vert\vert^4\right\rangle_{\mathscr{D}}   
 &=& \big(\frac{\zeta\sigma_0^2   }{(1-\zeta+1/N) } \big)^2\int_{0}^\infty\!\rmd\omega~ \Gamma_{N+1-d}\,(\omega)\,\omega^{-2}  \left[  \big(\frac{1}{d} \sum_{\ell=1}^d \frac{1  }{ \lambda_\ell(\Cov)} \big)^2+\frac{2}{d^2}\sum_{\ell=1}^d  \big(\frac{1   }{ \lambda_\ell(\Cov)}  \big)^2 \right] \nonumber\\
 &=& \big(\frac{\zeta\sigma_0^2   }{(1-\zeta+1/N) } \big)^2\int_{0}^\infty\!\rmd\omega~ \Gamma_{N+1-d}\,(\omega)\,\omega^{-2}\left[  \big(\frac{1}{d} \Tr\left[\Cov^{-1}\right] \big)^2+\frac{2}{d^2} \Tr\left[\Cov^{-2}\right]  \right] \nonumber\\
 &=& \big(\frac{\zeta\sigma_0^2   }{(1-\zeta+1/N) } \big)^2     \frac{   (1-\zeta+1/N)^2   }{(1-\zeta-1/N)(1-\zeta-3/N) } \left[  \big(\frac{1}{d} \Tr\left[\Cov^{-1}\right] \big)^2+\frac{2}{d^2} \Tr\left[\Cov^{-2}\right]  \right] \nonumber\\
&=&  \frac{   \zeta^2\sigma_0^4   }{(1-\zeta-1/N)(1-\zeta-3/N) }  \left[  \big(\frac{1}{d} \Tr\left[\Cov^{-1}\right] \big)^2+\frac{2}{d^2} \Tr\left[\Cov^{-2}\right]  \right] \label{eq:LR-ML-MSE-MGF-2-moment}
\end{eqnarray}
Now upon combining the mean  (\ref{eq:LR-ML-MSE-mean})  and the second moment (\ref{eq:LR-ML-MSE-MGF-2-moment}) we obtain the variance 
of the random variable $ \frac{1}{d}  \vert\vert\thetav_{\!0}-\thetaest[\mathscr{D}]\vert\vert^2$ for $(N,d)\rightarrow\infty$:%
\begin{eqnarray}
\left\langle \frac{1}{d^2}  \vert\vert\thetav_{\!0}-\thetaest[\mathscr{D}]\vert\vert^4\right\rangle_{\mathscr{D}} -  \frac{1}{d^2}\left\langle   \vert\vert\thetav_{\!0}-\thetaest[\mathscr{D}]\vert\vert^2 \right\rangle^2_{\mathscr{D}} 
&=&  \frac{   \zeta^2\sigma_0^4   }{(1-\zeta-1/N)(1-\zeta-3/N) }  \left[  \big(\frac{1}{d} \Tr\left[\Cov^{-1}\right] \big)^2+\frac{2}{d^2} \Tr\left[\Cov^{-2}\right]  \right] \nonumber\\
&&~~~~~~~~~-\frac{\zeta^2\sigma_0^4   }{\big(1-\zeta-1/N\big)^2}    \big(\frac{1}{d}  \Tr \left[\Cov^{-1}\right] \big)^2  \nonumber\\
&&\hspace*{-60mm} =\left[  \frac{   \zeta^2\sigma_0^4   }{(1-\zeta-1/N)(1-\zeta-3/N) }   -\frac{\zeta^2\sigma_0^4   }{\big(1-\zeta-1/N\big)^2} \right]   \frac{1}{d^2}  \Tr^2 \left[\Cov^{-1}\right]  + \frac{   \zeta^2\sigma_0^4   }{(1-\zeta-1/N)(1-\zeta-3/N) }   \frac{2}{d^2} \Tr\left[\Cov^{-2}\right]  \nonumber\\
&&~~~~~= 2\big(\frac{   \zeta\sigma_0^2   }{1-\zeta}\big)^2   \frac{1}{d^2} \Tr\left[\Cov^{-2}\right].  \label{eq:LR-ML-MSE-Var.}
\end{eqnarray}

\subsubsection{Properties of the MGF  for large $(N,d)$\label{appendix:mse-mgf-asympt}}

Let us consider the moment generating function (\ref{eq:LR-ML-MSE-MGF}) of the MSE for the covariance matrix $\Cov=\lambda\I_d$. The mean MSE $\langle   \frac{1}{d}  \vert\vert\thetav_{\!0}-\thetaest[\mathscr{D}]\vert\vert^2 \rangle_{\mathscr{D}}$ is given by  
\begin{eqnarray}
\mu(\lambda)&=& \zeta\sigma_0^2 /(1-\zeta)\lambda \label{eq:LR-ML-MSE-Cov-I-mean} 
\end{eqnarray}
using equation  (\ref{eq:LR-ML-MSE-mean}) for large $(N,d)$.   The MGF of the MSE is given by  
\begin{eqnarray}
\left\langle  \rme^{\frac{1}{2}\alpha  \vert\vert\thetav_{\!0}-\thetaest[\mathscr{D}]\vert\vert^2 }\right\rangle_{\mathscr{D}} &=& \int_{0}^\infty \Gamma_{N+1-d}\,(\omega)\,\rmd\omega\frac{1}{\big( 1  -\alpha  \frac{\zeta\sigma_0^2   }{\omega(1-\zeta+1/N) \lambda}  \big)^{\frac{d}{2}}}= \int_{0}^\infty\!\rmd\omega~ \Gamma_{N+1-d}\,(\omega)\,  \big( \frac{\omega(1-\zeta+1/N) \lambda}{\omega(1-\zeta+1/N) \lambda-\alpha\zeta\sigma_0^2}    \big)^{  \frac{d}{2}    }  
 \nonumber
 \\&
=& \int_{0}^\infty\!\rmd\omega~ \Gamma_{N+1-d}\,(\omega)\,     \Big( \frac{\omega}{\omega-  \alpha\frac{\zeta\sigma_0^2}{(1-\zeta+1/N) \lambda}   }    \Big)^{  \frac{d}{2}    } 
 = \int_{0}^\infty \!\rmd\omega~\Gamma_{N+1-d}\,(\omega)\big( \frac{\omega}{\omega-  \alpha  \mu(\lambda)   }    \big)^{  \frac{d}{2}    } \label{eq:LR-ML-MSE-MGF-Cov-I}  ,
\end{eqnarray}
where in the last line we assumed $(N,d)\to\infty$ . Let us now consider the integral 
\begin{eqnarray}
 \frac{1}{N}\log\int_{0}^\infty\!\rmd\omega~ \Gamma_{N+1-d}\,(\omega)\,\big( \frac{\omega}{\omega-\alpha    }    \big)^{  \frac{d}{2}    } 
  &=& \frac{1}{N}\log\frac{(N+1-d)^{(N+1-d)/2}}{2^{(N+1-d)/2}\Gamma((N+1-d)/2)}\nonumber\\
  &&~~~~~~~~+\frac{1}{N}\log \int_{0}^\infty \!\rmd\omega~\omega^{\frac{N-d-1}{2}}\rme^{-\frac{1}{2}(N-d+1)\omega}  \big( \frac{\omega}{\omega-\alpha    }    \big)^{  \frac{d}{2}    } \nonumber\\
   &&\hspace*{-40mm} =\frac{1}{N}\log\big(\frac{1-\zeta}{2}N\big)^{\frac{1-\zeta}{2}N}/\,\Gamma\big(\frac{1-\zeta}{2}N\big)
   +\frac{1}{N}\log \int_{0}^\infty\!\rmd\omega~ \rme^{ N\left[ \frac{1}{2} (1-\zeta) \log\omega    -\frac{1}{2}(1-\zeta)\,\omega+  \frac{\zeta}{2}    \log  \big( \frac{\omega}{\omega-\alpha    }    \big)\right]} 
   \nonumber\\
  &&\hspace*{-40mm} =         \frac{1-\zeta}{2}   +\frac{1}{2N}  \log  \big( \frac{1-\zeta}{4\pi} N \big) +O \big( {N}^{-2} \big)      
  + \frac{1}{2}\phi_{-}(\omega_0^{-}) +\frac{1}{2N}  \log  \big( \frac{4\pi}{N (-\phi_{-}^{\prime\prime}(\omega_0^{-}))}  \big)+O \big( {N}^{-2} \big)\nonumber\\
   &&\hspace*{-40mm}=  \frac{1-\zeta}{2} + \frac{1}{2}\phi_{-}(\omega_0^{-})+\frac{1}{2N}  \log  \big(              \frac{\zeta-1}{\phi_{-}^{\prime\prime}(\omega_0^{-})}        \big)+O \big( {N}^{-2} \big)\label{eq:Laplace-Integral}, 
   \end{eqnarray}
   where $\omega_0^{-}=\argmax_{\omega\in(0, \infty)} \phi_{-}(\omega)$, with  the function  
   \begin{eqnarray}
   \phi_{-}(\omega)=  (1-\zeta) \log\omega    -(1-\zeta)\,\omega+  \zeta\log  \big( \frac{\omega}{\omega-\alpha    }    \big)
   \label{def:phi-}
   \end{eqnarray}
   We note $\phi_{-}(\omega)$ has a maximum when the solution of    
   \begin{eqnarray}
   \phi_{-}^\prime(\omega)
&=&{\frac { \big( 1-\zeta \big) {\omega}^{2}- \big( \alpha+1 \big) 
 \big( 1-\zeta \big) \omega+\alpha}{ \big( \alpha-\omega
 \big) \omega}}=0\label{eq:phi-max}
   \end{eqnarray}
%
 %
%
satisfies the condition  $\phi_{-}^{\prime\prime}(\omega)>0$ given by the inequality ${\omega}^{2}\zeta- \big( \omega-\alpha \big) ^{2}<0$. The latter is satisfied when  $\omega\in\big(0, ( 1-\sqrt {\zeta}) \alpha/(1-\zeta)\big)\cup\big( ( 1+\sqrt {\zeta}) \alpha/(1-\zeta), \infty\big)$ for  $\zeta\in[\,0,1)$ and when  $\omega\in\big(0, \alpha/2\big)$ for $\zeta=1$.  However, the difference $\omega-\alpha$  for $\zeta\in[\,0,1)$ is negative on the interval $\big(0,( 1-\sqrt {\zeta}) \alpha/(1-\zeta)\big)$, so  $\phi_{-}(\omega)$ is undefined. The same is also true for $\big(0, \alpha/2\big)$ at  $\zeta=1$ thus leaving us only with the interval $\big( ( 1+\sqrt {\zeta}) \alpha/(1-\zeta), \infty\big)$ with  $\zeta\in(0,1)$. 
   
   The equation (\ref{eq:phi-max}) has real solutions   $\frac{1+\alpha}{2}  \pm {  \sqrt {\big( \frac{1+\alpha}{2}\big) ^{2}  -\frac{\alpha}{ \big( 1-\zeta \big)}      }}  $  when the inequality $\big( \frac{\alpha+1}{2} \big)^{2}/ \alpha\geq1/ \big( 1-\zeta \big)$  is satisfied. This holds when $\alpha\in\big(0, (1+\zeta-2\,\sqrt {\zeta})/(1-\zeta)\big)\cup \big((1+\zeta+2\,\sqrt {\zeta})/(1-\zeta), \infty\big)$, but for  $\zeta\in(0,1)$ only one solution $\omega_0=\frac{1+\alpha}{2}  + {  \sqrt {\big( \frac{1+\alpha}{2}\big) ^{2}  -\frac{\alpha}{ \big( 1-\zeta \big)}      }}  $ belongs  to the interval   $\big( ( 1+\sqrt {\zeta} ) \alpha)/(1-\zeta), \infty\big)$,   when $\alpha\in\big(0, (1+\zeta-2\,\sqrt {\zeta})/(1-\zeta)\big)$, i.e. it corresponds to the maximum of $\phi_{-}(\omega)$.  
   
We note that upon substituting  $ \alpha\to \alpha\mu(\lambda)$, the factor appearing in the integral (\ref{eq:LR-ML-MSE-MGF-Cov-I}), we obtain 
\begin{eqnarray}
\omega_0^{-}&=&\frac{1+\alpha  \mu(\lambda)}{2}               +              \sqrt { \big(    \frac{1+\alpha  \mu(\lambda)}{2}       \big) ^{2}   -    {\frac {\alpha\,\mu(\lambda)}{ \big( 1-\zeta \big)}}} \label{def:omega-}
\end{eqnarray}
for  $\alpha\  \in\big(0, {\lambda(1+\zeta-2\,\sqrt {\zeta})/\zeta\sigma_0^2 } \big)$ and $\zeta\in(0,1)$.  Now, using (\ref{eq:Laplace-Integral}), the moment generating function (\ref{eq:LR-ML-MSE-MGF-Cov-I}) is for large $N$ found to become 
\begin{eqnarray}
\left\langle   \rme^{\frac{1}{2}\alpha  \vert\vert\thetav_{\!0}-\thetaest[\mathscr{D}]\vert\vert^2 }\right\rangle_{\mathscr{D}} &=&  \sqrt{             \frac{\zeta-1}{\phi_{-}^{\prime\prime}(\omega_0^{-})}     }\, \rme^{\frac{1}{2}N \big(  1-\zeta + \phi_{-}(\omega_0^{-}) \big) +O \big( 1/N \big)}   \label{eq:LR-ML-MSE-MGF-Cov-I-large-N}.
\end{eqnarray}

\subsection{Deviations from the mean \label{appendix:deviations}}

We are interested in the probability of the event $ \frac{1}{d}\vert\vert\thetav_{\!0}-\thetaest[\mathscr{D}]\vert\vert^2  \notin \big(  \mu(\lambda)    -\delta,  \mu(\lambda) +\delta\big)$. This is given by 
\begin{eqnarray}
\mathrm{Prob} \left[  \frac{1}{d}\vert\vert\thetav_{\!0}-\thetaest[\mathscr{D}]\vert\vert^2  \notin \big(  \mu(\lambda)    -\delta,  \mu(\lambda) +\delta\big)     \right]&=&\mathrm{Prob} \left[\vert\vert\thetav_{\!0}-\thetaest[\mathscr{D}]\vert\vert^2 \leq d\big( \mu(\lambda)  -   \delta\big)\right]  
\nonumber
\\&&
+~ \mathrm{Prob} \left[   \vert\vert\thetav_{\!0}-\thetaest[\mathscr{D}]\vert\vert^2 \geq d\big(  \mu(\lambda)   +\delta\big)       \right]    \label{eq:LR-ML-MSE-Cov-I-dev.}.
\end{eqnarray}
First,  for $\alpha>0$  we  consider the probability 
\begin{eqnarray}
\mathrm{Prob} \left[  \vert\vert\thetav_{\!0}-\thetaest[\mathscr{D}]\vert\vert^2 \geq d\big(  \mu(\lambda)   +\delta\big)               \right]
&=& \mathrm{Prob} \left[ \rme^{\frac{1}{2}\alpha     \vert\vert\thetav_{\!0}-\thetaest[\mathscr{D}]\vert\vert^2      }\geq \rme^{\frac{1}{2}\alpha   d\big(   \mu(\lambda)  +\delta\big)    }\right]
\leq  \left\langle\rme^{\frac{1}{2}\alpha   \vert\vert\thetav_{\!0}-\thetaest[\mathscr{D}]\vert\vert^2       }  \right\rangle_{\mathscr{D}}  \rme^{-\frac{1}{2}\alpha d\big(  \mu(\lambda)   +\delta\big)       }\nonumber\\
&=&  \sqrt{             \frac{\zeta-1}{\phi_{-}^{\prime\prime}(\omega_0^{-})}     }\, \rme^{\frac{1}{2}N \big(  1-\zeta + \phi_{-}(\omega_0^{-}) \big)+O \big( 1/N \big)}   \rme^{-\frac{1}{2}\alpha \zeta N\big(   \mu(\lambda)   +\delta\big)       }\nonumber\\
&=&   \sqrt{             \frac{\zeta-1}{\phi_{-}^{\prime\prime}(\omega_0^{-})}     }\, \rme^{-\frac{1}{2}N \left[ \zeta-1 -  \phi_{-}(\omega_0^{-})+\alpha \zeta \big(   \mu(\lambda)   +\delta\big) \right] +O \big( 1/N \big)}   
\end{eqnarray}
From this it follows that for $N\rightarrow\infty$ we have 
\begin{eqnarray}
-\frac{2}{N}\log\mathrm{Prob} \left[  \vert\vert\thetav_{\!0}-\thetaest[\mathscr{D}]\vert\vert^2 \geq d\big(  \mu(\lambda)   +\delta\big)               \right]\geq    \zeta-1 -  \phi_{-}(\omega_0^{-})+\alpha\, \zeta \big(  \mu(\lambda)  +\delta\big)+O\big(1/N\big) \label{eq:LR-ML-MSE-Prob-lower-tail}
\end{eqnarray}
 We seek an upper bound with respect to $\alpha$, but it is not clear how to implement this analytically for any $\alpha$.  However, for small $\alpha$ the function (divided by $\zeta$) appearing  in the right-hand side of (\ref{eq:LR-ML-MSE-Prob-lower-tail}) has the following Taylor expansion: 
 \begin{eqnarray}
&&\big( \mu(\lambda)-1+\delta   \big) \alpha+ 
{\frac { \zeta\big(  \mu(\lambda)-1 \big) ^{2}-1  }{
2(1-\zeta)}}{\alpha}^{2}+{\frac {\big( \mu-1 \big) ^{3}{\zeta}^{2}+ \big( 3\,{\mu}
^{3}-3\,{\mu}^{2}-3\,\mu+2 \big) \zeta-1}{3 \big( 1-\zeta \big) ^{2
}}}{\alpha}^{3} 
+O
 \big( {\alpha}^{4} \big) , \label{eq:LR-ML-MSE-Prob-lower-tail-small-alpha}
\end{eqnarray}
so if $\mu(\lambda)+\delta>1$, the  first term in this expansion is positive and hence   if  $\alpha>0$  is sufficiently small then the RHS of  (\ref{eq:LR-ML-MSE-Prob-lower-tail}) is positive.  We note that for $\mu(\lambda)\geq1$, where $\mu(\lambda)=\zeta\sigma_0^2 /(1-\zeta)\lambda$,  the value of $\delta>0$ can be made arbitrary small, but for $\mu<1$ the positivity of first term in (\ref{eq:LR-ML-MSE-Prob-lower-tail-small-alpha}) is dependent on $\delta$.

Second,  for $\alpha>0$  we  consider the probability 
\begin{eqnarray}
\mathrm{Prob} \left[  \vert\vert\thetav_{\!0}-\thetaest[\mathscr{D}]\vert\vert^2 \leq d\big(  \mu(\lambda)   -\delta\big)               \right]
&=& \mathrm{Prob} \left[ \rme^{-\frac{1}{2}\alpha     \vert\vert\thetav_{\!0}-\thetaest[\mathscr{D}]\vert\vert^2      }\geq \rme^{-\frac{1}{2}\alpha   d\big(   \mu(\lambda)  -\delta\big)    }\right]
\leq  \left\langle\rme^{-\frac{1}{2}\alpha   \vert\vert\thetav_{\!0}-\thetaest[\mathscr{D}]\vert\vert^2       }\right\rangle_{\mathscr{D}}  \rme^{\frac{1}{2}\alpha d\big(  \mu(\lambda)   -\delta\big)       }\nonumber\\
&=&  \sqrt{             \frac{\zeta-1}{\phi_{+}^{\prime\prime}(\omega_0^{+})}     }\, \rme^{\frac{1}{2}N \big(  1-\zeta + \phi_{+}(\omega_0^{+}) \big)+O \big( 1/N \big)}   \rme^{\frac{1}{2}\alpha \zeta N\big(   \mu(\lambda)   -\delta\big)       }\nonumber\\
&=&  \sqrt{             \frac{\zeta-1}{\phi_{+}^{\prime\prime}(\omega_0^{+})}     }\, \rme^{-\frac{1}{2}N \left[ \zeta-1 -  \phi_{+}(\omega_0^{+})-\alpha \zeta \big(   \mu(\lambda)   -\delta\big) \right] +O \big( 1/N \big)}\label{eq:LR-ML-MSE-Prob-upper-tail-1}, 
\end{eqnarray}
where the function $\phi_{+}$, defined as 
 \begin{eqnarray}
 \phi_{+}(\omega)= (1-\zeta) \log\omega    -(1-\zeta)\,\omega+  \zeta\log  \big( \frac{\omega}{\omega+\alpha    }    \big)
 \label{def:phi+},
 \end{eqnarray}
 has a maximum at 
 \begin{eqnarray}
 \omega_0^{+}=\frac{1}{2}  \big(1-\alpha+\sqrt { \big( \alpha-1 \big) ^{2}+4 \alpha/(1-\zeta) }\big)
 \label{def:omega+}. 
 \end{eqnarray}
 Now for $\alpha\to \alpha\mu(\lambda)$  we have in the exponential of (\ref{eq:LR-ML-MSE-Prob-upper-tail-1}):
 \begin{eqnarray}
\zeta-1 -  \phi_{+}(\omega_0^{+})-\alpha \zeta \big(   \mu(\lambda)  \! -\!\delta\big)
&=& \delta\,\zeta\,\alpha -  {\frac {   {\mu}^{2}(\lambda)\zeta}{2\big(1\!-\!\zeta\big)}}{\alpha}^{2}                +            {\frac {{\mu}^{3}(\lambda) \zeta\,\big( \zeta\!+\!1 \big)     }{3\big( 1\!-\!\zeta \big) ^{2}}} {\alpha}^{3}
-      {\frac {{\mu}^{4}(\lambda)\zeta\, \big( {\zeta}^{2}\!+\!3\,\zeta\!+\!1 \big)}{4 \big( 1\!-\!\zeta \big) ^{3}}} {\alpha}^{4} +O(\alpha^5).~~~
 \end{eqnarray}
This  is positive for sufficiently small $\alpha$, and hence the lower bound in the inequality 
\begin{eqnarray}
-\frac{2}{N}\log\mathrm{Prob} \left[  \vert\vert\thetav_{\!0}-\thetaest[\mathscr{D}]\vert\vert^2 \leq d\big(  \mu(\lambda)   -\delta\big)               \right]
\geq   \zeta-1 -  \phi_{+}(\omega_0^{+})-\alpha \zeta \big(   \mu(\lambda)   -\delta\big)  +O \big( 1/N \big)
\label{eq:LR-ML-MSE-Prob-upper-tail-2} 
\end{eqnarray}
 is positive for any $\delta\in(0,\mu(\lambda))$ and sufficiently small $\alpha$, when $N\rightarrow\infty$. 
 
 Now combining (\ref{eq:LR-ML-MSE-Prob-lower-tail}) with (\ref{eq:LR-ML-MSE-Prob-upper-tail-2}) allows us bound the  probability  (\ref{eq:LR-ML-MSE-Cov-I-dev.}) as follows
\begin{eqnarray}
\mathrm{Prob} \left[  \frac{1}{d}\vert\vert\thetav_{\!0}-\thetaest[\mathscr{D}]\vert\vert^2  \notin \big(  \mu(\lambda)    -\delta,  \mu(\lambda) +\delta\big)     \right]&\leq& \mathrm{C}_{-}\rme^{-\frac{N}{2}\left[ \zeta-1 -  \phi_{-}(\omega_0^{-})+\alpha\, \zeta \big(  \mu(\lambda)  +\delta\big) \right]}  
\nonumber
\\[-2mm]&&
+~ \mathrm{C}_{+}\rme^{-\frac{N}{2}\left[  \zeta-1 -  \phi_{+}(\omega_0^{+})-\alpha \zeta \big(   \mu(\lambda)   -\delta\big)\right]}    \label{eq:LR-ML-MSE-Cov-I-dev.-ub}.
\end{eqnarray}
for some \emph{constants} $\mathrm{C}_{\pm}$  and some sufficiently small $\alpha>0$. We note that  for the first term in the above upper bound to  vanish,  as $N\rightarrow\infty$,  for arbitrary small  $\delta$ it is sufficient that   $\mu(\lambda)\geq1$, where $\mu(\lambda)=\zeta\sigma_0^2 /(1-\zeta)\lambda$, but  for $\mu(\lambda)<1$  the value of $\delta$ must be such that $\delta>1-\mu(\lambda)$.  The second term in the upper bound is vanishing for any  $\delta\in(0,\mu(\lambda))$.

Finally we consider deviations of $ \frac{1}{d} \vert\vert\thetav_{\!0}-\thetaest[\mathscr{D}]\vert\vert^2  $ from its mean $\mu\big(\Cov\big)= [\zeta\sigma_0^2 /(1\!-\!\zeta\!-\!N^{-1})]    \frac{1}{d}  \Tr\left[ \Cov^{-1}\right]$, derived in  (\ref{eq:LR-ML-MSE-mean}). To this end we consider the probability of event  $ \frac{1}{d}\vert\vert\thetav_{\!0}-\thetaest[\mathscr{D}]\vert\vert^2 \notin \big( \mu\big(\Cov\big) -\delta,\mu\big(\Cov\big) +\delta    \big)$ given by the sum
\begin{eqnarray}
&&\hspace*{-25mm} 
\mathrm{Prob} \left[  \vert\vert\thetav_{\!0}-\thetaest[\mathscr{D}]\vert\vert^2         \leq d\big( \mu\big(\Cov\big)  -   \delta\big)\right]  
+ \mathrm{Prob} \left[   \vert\vert\thetav_{\!0}-\thetaest[\mathscr{D}]  \vert\vert^2 \geq d\big(  \mu\big(\Cov\big)  +\delta\big)       \right] \nonumber\\
&\leq&      \left\langle\rme^{-\frac{1}{2}\alpha   \vert\vert\thetav_{\!0}-\thetaest [\mathscr{D}] \vert\vert^2       } \right\rangle_{\mathscr{D}}  \rme^{\frac{1}{2}\alpha d\big(  \mu\big(\Cov\big)   -\delta\big)       }
+   \left\langle\rme^{\frac{1}{2}\alpha   \vert\vert\thetav_{\!0}-\thetaest[\mathscr{D}]\vert\vert^2       }\right\rangle_{\mathscr{D}}  \rme^{-\frac{1}{2}\alpha d\big(  \mu\big(\Cov\big)   +\delta\big)       }
\end{eqnarray}
for $\alpha>0$. Let us order the eigenvalues of $\Cov$ in a such a way that $\lambda_1 (\Cov)\leq \lambda_2 (\Cov)\leq\cdots\leq\lambda_d (\Cov)$ then, using (\ref{eq:LR-ML-MSE-MGF}),  for the moment generating functions in above we obtain 
\begin{eqnarray}
\left\langle   \rme^{-\frac{1}{2}\alpha  \vert\vert\thetav_{\!0}-\thetaest[\mathscr{D}]\vert\vert^2 }   \right\rangle_{\mathscr{D}}
&=& \int_{0}^\infty \Gamma_{N+1-d}\,(\omega)\,\rmd\omega\prod_{\ell=1}^d  \big( 1  +\alpha  \frac{\zeta\sigma_0^2   }{\omega(1-\zeta+1/N) \lambda_\ell (\Cov)}  \big)^{-\frac{1}{2}}\nonumber\\
 &\leq& \int_{0}^\infty \Gamma_{N+1-d}\,(\omega)\,\rmd\omega \big( 1  +\alpha  \frac{\zeta\sigma_0^2   }{\omega(1-\zeta+1/N) \lambda_1 (\Cov)}  \big)^{-\frac{d}{2}}\label{eq:LR-ML-MSE-MGF-ub-1}
\end{eqnarray}
and 
\begin{eqnarray}
\left\langle   \rme^{\frac{1}{2}\alpha  \vert\vert\thetav_{\!0}-\thetaest[\mathscr{D}]\vert\vert^2 }   \right\rangle_{\mathscr{D}} 
&=& \int_{0}^\infty \Gamma_{N+1-d}\,(\omega)\,\rmd\omega\prod_{\ell=1}^d  \big( 1  -\alpha  \frac{\zeta\sigma_0^2   }{\omega(1-\zeta+1/N) \lambda_\ell (\Cov)}  \big)^{-\frac{1}{2}}\nonumber\\
 &\leq& \int_{0}^\infty \Gamma_{N+1-d}\,(\omega)\,\rmd\omega\big( 1  -\alpha  \frac{\zeta\sigma_0^2   }{\omega(1-\zeta+1/N) \lambda_d (\Cov)}  \big)^{-\frac{d}{2}}\label{eq:LR-ML-MSE-MGF-ub-2}.
\end{eqnarray}
Furthermore,  by the inequalities  $1/ \lambda_d (\Cov)\leq\frac{1}{d}  \Tr\left[ \Cov^{-1}\right]\leq 1/ \lambda_1 (\Cov)$, the mean obeys $\mu\big(\lambda_d (\Cov)\big) \leq \mu\big(\Cov\big) \leq  \mu\big(\lambda_1 (\Cov)\big)$, where $\mu\big(\lambda\big)=\zeta\sigma_0^2 /(1\!-\!\zeta\!-\!N^{-1}) \lambda$. The latter combined with the upper bounds in (\ref{eq:LR-ML-MSE-MGF-ub-1}) and (\ref{eq:LR-ML-MSE-MGF-ub-2}) gives us 
\begin{eqnarray}
\mathrm{Prob} \left[  \frac{1}{d}\left\vert\left\vert\thetav_{\!0}-\thetaest[\mathscr{D}]\right\vert\right\vert^2\!  \notin\! \big( \mu\big(\Cov\big) \!   -\!\delta,  \mu\big(\Cov\big) \!+\!\delta\big)     \right]
&\leq&    
 \int_{0}^\infty\!\!\rmd\omega~ \Gamma_{N+1-d}\,(\omega)\big( 1  \!+\!\frac{\alpha\zeta\sigma_0^2   }{\omega(1-\zeta+1/N) \lambda_1 (\Cov)}  \big)^{-\frac{d}{2}} \rme^{\frac{1}{2}\alpha d\big(  \mu\big(\lambda_1 (\Cov)\big)   -\delta\big)       }\nonumber\\
&&\hspace*{-20mm} +   
 \int_{0}^\infty\!\!\rmd\omega~ \Gamma_{N+1-d}\,(\omega)\big( 1 \! -\!\frac{\alpha\zeta\sigma_0^2   }{\omega(1-\zeta+1/N) \lambda_d (\Cov)}  \big)^{-\frac{d}{2}}
 \rme^{-\frac{1}{2}\alpha d\big(  \mu\big(  \lambda_d (\Cov)    \big)   +\delta\big)       }
\end{eqnarray}
Finally, using similar  steps to those that led us earlier to (\ref{eq:LR-ML-MSE-Cov-I-dev.-ub}) give
\begin{eqnarray}
&&\mathrm{Prob} \left[  \frac{1}{d}\left\vert\left\vert\thetav_{\!0}-\thetaest[\mathscr{D}]\right\vert\right\vert^2  \notin \big( \mu\big(\Cov\big)    -\delta,  \mu\big(\Cov\big) +\delta\big)     \right]
\leq\mathrm{C}_{-}\rme^{-N \Phi_{-}\left[\alpha, \mu(\lambda_d ), \delta\right]  }  + \mathrm{C}_{+}\rme^{-N\Phi_{+}\!\left[\alpha, \mu(\lambda_1 ), \delta\right]} ,
\label{eq:extra}
\end{eqnarray}
for some constants $\mathrm{C}_{\pm}$  and some sufficiently small $\alpha>0$.  In the above we have defined the (rate) functions 
\begin{eqnarray}
\Phi_{-}\!\left[\alpha, \mu(\lambda_d ), \delta\right]&=&\frac{1}{2}\Big[\zeta-1 -  \phi_{-}(\omega_0^{-}  )+\alpha\, \zeta \big(  \mu(\lambda_d (\Cov))   +\delta\big)\Big],
\\
\Phi_{+}\!\left[\alpha, \mu(\lambda_1 ), \delta\right]&=& \frac{1}{2}\Big[ \zeta-1 -  \phi_{+}(\omega_0^{+} )-\alpha \zeta \big(   \mu\big(\lambda_1 (\Cov)\big)   -\delta\big)\Big].
\end{eqnarray}
Here $\phi_{-}(\omega_0^{-}  )$ is defined by (\ref{def:phi-}) and (\ref{def:omega-}) with $\mu(\lambda)$ replaced by $\mu(\lambda_d)$,  and
 $\phi_{+}(\omega_0^{+})$ is defined by (\ref{def:phi+}) and (\ref{def:omega+}) with $\alpha$ replaced by $\alpha\mu(\lambda_1)$. We note that  for the first term in the upper bound (\ref{eq:extra}) to  vanish,  as $(N,d)\rightarrow\infty$,  for arbitrary small  $\delta$, it is sufficient that   $\mu(\lambda_d)\geq1$, where $\mu(\lambda)=\zeta\sigma_0^2 /(1-\zeta)\lambda$, but  for $\mu(\lambda_d)<1$  for this to happen the $\delta$ must be such that $\delta>1-\mu(\lambda_d)$.  The second term in the upper bound is vanishing for any  $\delta\in(0,\mu(\lambda))$. 
 
\section{Statistical properties of free energy\label{section:free-energy}}
In this section we consider statistical properties of the (conditional) free energy  
\begin{eqnarray}
 F_{\beta,\sigma^2}\left[\mathscr{D}\right]&=&\frac{d}{2\beta} +\frac{1}{2\sigma^2}\tv^\Tran\big(\I_N-\Zm \J^{-1}_{\sigma^2\eta} \Zm^{T}\big)\tv   -\frac{1}{2\beta} \log\left\vert 2\pi\rme\sigma^2 \beta^{-1}  \J_{\sigma^2\eta}^{-1}\right \vert  \label{eq:LR-F-cond.-App}
\end{eqnarray}
assuming that $\sigma^2$ is independent from data $\mathscr{D}$.

\subsection{The average of  free energy}
Let us  consider the average free energy
\begin{eqnarray}
 \left\langle F_{\beta,\sigma^2}\left[\mathscr{D}\right] \right\rangle_{\mathscr{D}}&=&\frac{d}{2\beta} +\frac{1}{2\sigma^2} \left\langle\tv^\Tran\big(\I_N-\Zm \J^{-1}_{\sigma^2\eta} \Zm^{T}\big)\tv  \right\rangle_{\mathscr{D}}   -\frac{1}{\beta}\frac{1}{2} \left\langle \log\vert 2\pi\rme\sigma^2 \beta^{-1}  \J_{\sigma^2\eta}^{-1} \vert \right\rangle_{\mathscr{D}}\nonumber\\  
&=&\frac{d}{2\beta} +\frac{1}{2\sigma^2} \left\langle \left\langle\tv^\Tran\big(\I_N-\Zm \J^{-1}_{\sigma^2\eta} \Zm^{T}\big)\tv  \right\rangle_{\epsilonv}\right\rangle_{\Zm}   -\frac{1}{\beta}\frac{1}{2} \left\langle \log\vert 2\pi\rme\sigma^2 \beta^{-1}  \J_{\sigma^2\eta}^{-1} \vert \right\rangle_{\Zm} 
%
\label{eq:LR-F-cond.-aver.-calc.-1}
\end{eqnarray}
Now,  assuming   that   the noise vector $\epsilonv$ has mean $\nullv$  and covariance $\sigma_0^2\I_N$,  we can work out
\begin{eqnarray}
\left\langle\tv^\Tran\big(\I_N-\Zm \J^{-1}_{\sigma^2\eta} \Zm^{T}\big)\tv  \right\rangle_{\epsilonv}
&=&\left\langle  \Tr\left[ \big(\I_N-\Zm \J^{-1}_{\sigma^2\eta} \Zm^{T}\big)\tv\tv^\Tran\right ]   \right\rangle_{\epsilonv}
=  \Tr\left[ \big(\I_N-\Zm \J^{-1}_{\sigma^2\eta} \Zm^{T}\big) \left\langle\big(\Zm\thetav_{\!0}+\epsilonv\big)   \big(\Zm\thetav_{\!0}+\epsilonv\big)^\Tran \right\rangle_{\epsilonv} \right ] \nonumber\\
&&=  \Tr\left[ \big(\I_N-\Zm \J^{-1}_{\sigma^2\eta} \Zm^{T}\big) \big(\Zm\thetav_{\!0}  \thetav_{\!0}^\Tran    \Zm^\Tran+ 2\Zm\thetav_{\!0}\left\langle\epsilonv^\Tran \right\rangle_{\epsilonv}  +\left\langle\epsilonv\epsilonv^\Tran \right\rangle_{\epsilonv} \big) \right ]  \nonumber\\
&&=  \Tr\left[ \big(\I_N-\Zm \J^{-1}_{\sigma^2\eta} \Zm^{T}\big) \big(\Zm\thetav_{\!0}  \thetav_{\!0}^\Tran    \Zm^\Tran  +\sigma_0^2\I_N \big) \right ]  \nonumber\\
&&=  \thetav_{\!0}^\Tran    \Zm^\Tran  \big(\I_N-\Zm \J^{-1}_{\sigma^2\eta} \Zm^{T}\big) \Zm\thetav_{\!0}   + \sigma_0^2\Tr\left[ \I_N-\Zm \J^{-1}_{\sigma^2\eta} \Zm^{T}   \right ]  \nonumber\\
&&=  \thetav_{\!0}^\Tran    \big(\J-\J \J^{-1}_{\sigma^2\eta} \J\big) \thetav_{\!0}   + \sigma_0^2 \big(N-\Tr\left[\J \J^{-1}_{\sigma^2\eta}  \right ] \big),
\end{eqnarray}
and hence  the average free energy is given by 
\begin{eqnarray}
 \left\langle F_{\beta,\sigma^2}\left[\mathscr{D}\right] \right\rangle_{\mathscr{D}} &=&\frac{d}{2\beta}+\frac{1}{2\sigma^2}      \thetav_{\!0}^\Tran     \left\langle\big(\J-\J \J^{-1}_{\sigma^2\eta} \J\big) \right\rangle_{\Zm}\! \thetav_{\!0}   + \frac{\sigma_0^2}{2\sigma^2}\! \big(N-\left\langle\Tr\left[\J \J^{-1}_{\sigma^2\eta}  \right ] \right\rangle_{\Zm} \big) -\frac{1}{2\beta}\! \left\langle \log\left\vert 2\pi\rme\sigma^2 \beta^{-1}  \J_{\sigma^2\eta}^{-1}\right \vert \right\rangle_{\Zm}.
 ~~~~\label{eq:LR-F-cond.-aver.-calc.-2}
\end{eqnarray}

\subsection{The variance of free energy\label{appendix:free-energy-variance}}
We turn to the variance of the free energy $F_{\beta,\sigma^2}\left[\mathscr{D}\right]$.  To this end we exploit the free energy equality $F=U-TS$, which gives us  $\Var (F)=\Var(U-TS)=\Var(U)-2T\Covar(U,S)+T^2\Var(S)$. The latter applied  to (\ref{eq:LR-F-cond.-App}) leads to
\begin{eqnarray}
 \Var \big(F_{\beta,\sigma^2}\left[\mathscr{D}\right]\big)&=&\Var \big(E[\mathscr{D}]\big)+T^2 \Var\big(\mathrm{S}[\mathscr{D}]\big) - 2T\,\Covar\big( E[ \mathscr{D}], \mathrm{S}[\mathscr{D}]  \big)\label{eq:LR-Var-F}
\end{eqnarray}
Let us consider the energy variance 
\begin{eqnarray}
\Var \big(E[\mathscr{D}]\big)&=&\Var \big(\frac{d}{2\beta} +\frac{1}{2\sigma^2} \tv^\Tran\big(\I_N-\Zm \J^{-1}_{\sigma^2\eta} \Zm^{T}\big)\tv  \big)\label{eq:LR-Var-E-calc.-1}\\
&=&\frac{1}{4\sigma^4}\Var \big( \big(\Zm\thetav_{\!0}+\epsilonv\big)^\Tran\big(\I_N-\Zm \J^{-1}_{\sigma^2\eta} \Zm^{T}\big)\big(\Zm\thetav_{\!0}+\epsilonv\big)  \big)\nonumber
\end{eqnarray}
If we define $\vv=\Zm\thetav_{\!0}$ and $\A=(\I_N-\Zm \J^{-1}_{\sigma^2\eta} \Zm^{T})$ then the above is of the form 
\begin{eqnarray}
\Var \big( \big(\vv+\epsilonv\big)^\Tran\A\big(\vv+\epsilonv\big)  \big)
&=&\Var \big( \vv^\Tran\A\vv+2\,\epsilonv^\Tran\A\vv+\epsilonv^\Tran\A\epsilonv\big)\nonumber\\
&=&\Var \big( \vv^\Tran\A\vv\big)+4\,\Var \big(\epsilonv^\Tran\A\vv\big)+\Var \big(\epsilonv^\Tran\A\epsilonv\big)\nonumber\\
&&+2\left[2\,\Covar \big( \vv^\Tran\A\vv,\epsilonv^\Tran\A\vv    \big)+\Covar \big( \vv^\Tran\A\vv,\epsilonv^\Tran\A\epsilonv    \big)+2\,\Covar \big(\epsilonv^\Tran\A\vv  ,\epsilonv^\Tran\A\epsilonv\big)\right]\label{eq:Var-id}.
\end{eqnarray}
In the following we will use the following identities:
\begin{eqnarray}
\Zm^\Tran\A\Zm=\J-\J \J^{-1}_{\sigma^2\eta}\J &~~~~~~&
\Tr [\A]=N-\Tr[\J \J^{-1}_{\sigma^2\eta}],\\
\Zm^\Tran\A^2\Zm
= \J\big(\I_d -\J^{-1}_{\sigma^2\eta} \J\big)^2
&& \Tr\left[ \A^2\right]=N-2\Tr[\J \J^{-1}_{\sigma^2\eta}]  +\Tr [(\J \J^{-1}_{\sigma^2\eta})^2]  
\end{eqnarray}
We can 
 now compute  each term in  (\ref{eq:Var-id}): 
\begin{eqnarray}
\Var \big( \vv^\Tran\A\vv\big)&=&\Var \big(  \thetav_{\!0}^\Tran \Zm^\Tran  \big(\I_N-\Zm \J^{-1}_{\sigma^2\eta} \Zm^{T}\big)  \Zm\thetav_{\!0}    \big)
=\Var \big(\thetav_{\!0}^\Tran   \big(\J-\J \J^{-1}_{\sigma^2\eta} \J\big) \thetav_{\!0}    \big)
\nonumber
\\&
=&  \left\langle\big(\thetav_{\!0}^\Tran   \big(\J-\J \J^{-1}_{\sigma^2\eta} \J\big) \thetav_{\!0}    \big)^2 \right\rangle_{\Zm}
   -  \left\langle\thetav_{\!0}^\Tran   \big(\J-\J \J^{-1}_{\sigma^2\eta} \J\big) \thetav_{\!0}    \right\rangle^2_{\Zm}
\\[2mm]
\Var \big(\epsilonv^\Tran\A\vv\big)&=&
\left\langle\left\langle\big(\epsilonv^\Tran\A\vv\big)^2\right\rangle_{\epsilonv}\right\rangle_{\Zm} -  \left\langle\left\langle \epsilonv^\Tran\A\vv     \right\rangle_{\epsilonv}\right\rangle^2_{\Zm} 
=\left\langle\left\langle\big(\epsilonv^\Tran\A\vv\big)^2\right\rangle_{\epsilonv}\right\rangle_{\Zm}
=\left\langle    \vv^\Tran\A^\Tran \left\langle \epsilonv \epsilonv^\Tran \right\rangle_{\epsilonv}\A\vv     \right\rangle_{\Zm}
\nonumber
\\&=& 
\sigma_0^2  \left\langle    \vv^\Tran\A^2\vv     \right\rangle_{\Zm}= \sigma_0^2  \left\langle\thetav_{\!0}^\Tran \Zm^\Tran  \big(\I_N-\Zm \J^{-1}_{\sigma^2\eta} \Zm^{T}\big)^2  \Zm\thetav_{\!0}\right\rangle_{\Zm}\nonumber \\
&=&\sigma_0^2  \thetav_{\!0}^\Tran  \left\langle \J\big( \I_d -  \J^{-1}_{\sigma^2\eta}\J\big)^2  \right\rangle_{\Zm}    \thetav_{\!0}
\\[2mm]
\Var \big(\epsilonv^\Tran\A\epsilonv\big)&=&\left\langle\left\langle\big(\epsilonv^\Tran\A\epsilonv\big)^2 \right\rangle_{\epsilonv}\right\rangle_{\Zm} -\left\langle\left\langle\epsilonv^\Tran\A\epsilonv \right\rangle_{\epsilonv}\right\rangle^2_{\Zm}
= \sigma_0^4\big(\left\langle\Tr^2[\A]\right\rangle_{\Zm} + 2 \left\langle\Tr[\A^2]\right\rangle_{\Zm}- \left\langle\Tr[\A]\right\rangle^2_{\Zm}\big)
\nonumber \\
&=& \sigma_0^4\Big[ \left\langle \big(N\!-\!\Tr[\J \J^{-1}_{\sigma^2\eta}]\big)^2  \right\rangle_{\Zm} + 2 \left\langle N\!-\!2\Tr[\J \J^{-1}_{\sigma^2\eta}]  +\Tr [\big(\J \J^{-1}_{\sigma^2\eta}\big)^2]   \right\rangle_{\Zm}- \left\langle  N\!-\!\Tr[\J \J^{-1}_{\sigma^2\eta}]    \right\rangle^2_{\Zm}\Big]
\\[2mm]
\Covar \big( \vv^\Tran\A\vv,\epsilonv^\Tran\A\vv    \big)&=&\left\langle\left\langle  \vv^\Tran\A\vv\epsilonv^\Tran\A\vv  \right\rangle_{\epsilonv}\right\rangle_{\Zm}-\left\langle\left\langle   \vv^\Tran\A\vv \right\rangle_{\epsilonv}\right\rangle_{\Zm}       \left\langle \left\langle   \epsilonv^\Tran\A\vv    \right\rangle_{\epsilonv}\right\rangle_{\Zm}=0
\\[2mm]
\Covar \big( \vv^\Tran\A\vv,\epsilonv^\Tran\A\epsilonv    \big)&=&\left\langle\left\langle\vv^\Tran\A\vv\epsilonv^\Tran\A\epsilonv    \right\rangle_{\epsilonv}\right\rangle_{\Zm} - \left\langle\left\langle\vv^\Tran\A\vv\right\rangle_{\epsilonv}\right\rangle_{\Zm}  \left\langle\left\langle\epsilonv^\Tran\A\epsilonv    \right\rangle_{\epsilonv}\right\rangle_{\Zm}
\nonumber \\
&=& \left\langle\vv^\Tran\A\vv  \Tr\A\left\langle\epsilonv\epsilonv^\Tran    \right\rangle_{\epsilonv}\right\rangle_{\Zm} - \left\langle\vv^\Tran\A\vv \right\rangle_{\Zm}  \left\langle\left\langle\epsilonv^\Tran\A\epsilonv    \right\rangle_{\epsilonv}\right\rangle_{\Zm}
\nonumber \\
&=& \sigma_0^2\Big[\left\langle    \thetav_{\!0}^\Tran   \big(\J\!-\!\J \J^{-1}_{\sigma^2\eta} \J\big) \thetav_{\!0}      \big(N\!-\!\Tr[ \J\J^{-1}_{\sigma^2\eta}]\big)\right\rangle_{\Zm}- \left\langle  \thetav_{\!0}^\Tran   \big(\J\!-\!\J \J^{-1}_{\sigma^2\eta} \J\big) \thetav_{\!0}     \right\rangle_{\Zm}  \left\langle \big(N\!-\!\Tr[ \J\J^{-1}_{\sigma^2\eta}]\big)    \right\rangle_{\Zm}\Big]~~~
\\[2mm]
\Covar \big(\epsilonv^\Tran\A\vv  ,\epsilonv^\Tran\A\epsilonv\big)&=& \left\langle\left\langle \epsilonv^\Tran\A\vv\epsilonv^\Tran\A\epsilonv\right\rangle_{\epsilonv}\right\rangle_{\Zm}  -  \left\langle\left\langle \epsilonv^\Tran\A\vv \right\rangle_{\epsilonv}\right\rangle_{\Zm}      \left\langle\left\langle\epsilonv^\Tran\A\epsilonv    \right\rangle_{\epsilonv}\right\rangle_{\Zm}=0 
\end{eqnarray}
In deriving the above identities  we also used the following result 
\begin{eqnarray}
\left\langle \big(\epsilonv^\Tran\A\epsilonv\big)^2\right\rangle_{\epsilonv}
&=&\sum_{i_1,\ldots,i_4} \left\langle  \epsilon_{i_1}\epsilon_{i_2}\epsilon_{i_3}\epsilon_{i_4 }\right\rangle_{\epsilonv}\A_{i_1i_2}\A_{i_3i_4}
=\sigma_0^4\sum_{i_1,\ldots,i_4} \Big[
 \delta_{i_1,i_2}\delta_{i_3, i_4 }%
+ \delta_{i_1, i_3}\delta_{i_2, i_4 }
+ \delta_{i_1, i_4}\delta_{i_2, i_3 }
\Big]
\A_{i_1i_2}\A_{i_3i_4}\nonumber\\
&=&\sigma_0^4\sum_{i_1,i_2} \left\{\A_{i_1i_1} \A_{i_2 i_2}  + 2\A^2_{i_1i_2}\right\}
=\sigma_0^4 \big(\Tr^2[\A] + 2\Tr[\A^2]\big).\label{eq:Wicks-theorem}
\end{eqnarray}
Using all of the above results in (\ref{eq:LR-Var-E-calc.-1}) we obtain 
\begin{eqnarray}
4\sigma^4\Var \big(E[\mathscr{D}]\big)
&=&\left\langle\big(\thetav_{\!0}^\Tran   \big(\J\!-\!\J \J^{-1}_{\sigma^2\eta} \J\big) \thetav_{\!0}    \big)^2 \right\rangle_{\Zm}
 \!  -  \left\langle\thetav_{\!0}^\Tran   \big(\J\!-\!\J \J^{-1}_{\sigma^2\eta} \J\big) \thetav_{\!0}    \right\rangle^2_{\Zm} 
 +4\sigma_0^2  \thetav_{\!0}^\Tran  \left\langle \J\big( \I_d -  \J^{-1}_{\sigma^2\eta}\J\big)^2  \right\rangle_{\Zm}    \thetav_{\!0}\nonumber\\
&&+\sigma_0^4\Bigg[ \left\langle \big(N\!-\!\Tr\left[\J \J^{-1}_{\sigma^2\eta}\right]\big)^2  \right\rangle_{\Zm}
 + 2 \left\langle  \big(N\!-\!2\Tr\left[\J \J^{-1}_{\sigma^2\eta} \right]  +\Tr \left[\big(\J \J^{-1}_{\sigma^2\eta}\big)^2\right] \big)   \right\rangle_{\Zm}
 - \left\langle  \big(N\!-\!\Tr\left[\J \J^{-1}_{\sigma^2\eta}\right]\big)    \right\rangle^2_{\Zm}\Bigg]\nonumber\\
 &&+ 2\sigma_0^2\Bigg(\left\langle    \thetav_{\!0}^\Tran   \big(\J\!-\!\J \J^{-1}_{\sigma^2\eta} \J\big) \thetav_{\!0}      \big(N\!-\!\Tr[ \J\J^{-1}_{\sigma^2\eta}]\big)\right\rangle_{\Zm}- \left\langle  \thetav_{\!0}^\Tran   \big(\J\!-\!\J \J^{-1}_{\sigma^2\eta} \J\big) \thetav_{\!0}     \right\rangle_{\Zm}  \left\langle \big(N\!-\!\Tr[ \J\J^{-1}_{\sigma^2\eta}]\big)    \right\rangle_{\Zm}\Bigg) \label{eq:LR-Var-E-calc.-2}.
\end{eqnarray}
The entropy variance is given by 
\begin{eqnarray}
\Var\big(\mathrm{S}[\mathscr{D}]\big)&=&\Var\big(\frac{1}{2} \log\vert 2\pi\rme\sigma^2 \beta^{-1}  \J_{\sigma^2\eta}^{-1} \vert          \big)
=\frac{1}{4}\left\langle \log^2\left\vert  \J_{\sigma^2\eta}^{-1} \right\vert \right\rangle_{\Zm}
-  \frac{1}{4} \left\langle \log\left\vert  \J_{\sigma^2\eta}^{-1} \right\vert \right\rangle^2_{\Zm}\label{eq:LR-Var-S}
\end{eqnarray}
and the covariance is
\begin{eqnarray}
\Covar\big( E[ \mathscr{D}], \mathrm{S}[\mathscr{D}]  \big)
&=&\Covar\big(\frac{d}{2\beta} +\frac{1}{2\sigma^2} \tv^\Tran\big(\I_N-\Zm \J^{-1}_{\sigma^2\eta} \Zm^{T}\big)\tv  , \frac{d}{2} \log(2\pi\rme\sigma^2 \beta^{-1})  +\frac{1}{2}\log\left\vert\J_{\sigma^2\eta}^{-1}\right \vert          \big)\nonumber\\
&=&\frac{1}{4\sigma^2} \Covar\Bigg(\tv^\Tran\big(\I_N-\Zm \J^{-1}_{\sigma^2\eta} \Zm^{T}\big)\tv  ,\log\left\vert\J_{\sigma^2\eta}^{-1}\right \vert          \Bigg)\nonumber\\
&=&\frac{1}{4\sigma^2} \left\langle\left\langle\tv^\Tran\big(\I_N-\Zm \J^{-1}_{\sigma^2\eta} \Zm^{T}\big)\tv\right\rangle_{\epsilonv} \log\left\vert\J_{\sigma^2\eta}^{-1}\right \vert      \right\rangle_{\Zm}-\frac{1}{4\sigma^2}\left\langle\left\langle\tv^\Tran\big(\I_N-\Zm \J^{-1}_{\sigma^2\eta} \Zm^{T}\big)\tv \right\rangle_{\epsilonv}   \right\rangle_{\Zm}\left\langle\log\left\vert\J_{\sigma^2\eta}^{-1}\right \vert  \right\rangle_{\Zm}   \nonumber\\
&=&\frac{1}{4\sigma^2} \left\langle  \left[  \thetav_{\!0}^\Tran    \big(\J-\J \J^{-1}_{\sigma^2\eta} \J\big) \thetav_{\!0}  + \sigma_0^2 \big(N-\Tr\left[\J \J^{-1}_{\sigma^2\eta}  \right ] \big) \right]  \log\left\vert\J_{\sigma^2\eta}^{-1}\right \vert      \right\rangle_{\Zm}\nonumber\\
&&~~~~~-\frac{1}{4\sigma^2}\Big\langle\Big[  \thetav_{\!0}^\Tran    \big(\J-\J \J^{-1}_{\sigma^2\eta} \J\big) \thetav_{\!0} + \sigma_0^2 \big(N-\Tr\left[\J \J^{-1}_{\sigma^2\eta}  \right ] \big) \Big]  \Big\rangle_{\Zm}\left\langle\log\left\vert\J_{\sigma^2\eta}^{-1}\right \vert  \right\rangle_{\Zm}   \label{eq:LR-Covar-E-S}
\end{eqnarray}
Finally,  using the results (\ref{eq:LR-Var-E-calc.-2}),  (\ref{eq:LR-Var-S}) and  (\ref{eq:LR-Covar-E-S}), the variance of free energy follows from equation   (\ref{eq:LR-Var-F}).

\subsection{Free energy of ML inference}
 For $\eta=0$ we simply have  $\J_{\sigma^2\eta}=\J$ and  the average free energy (\ref{eq:LR-F-cond.-aver.-calc.-2}) is given by 
\begin{eqnarray}
 \left\langle F_{\beta,\sigma^2}\left[\mathscr{D}\right] \right\rangle_{\mathscr{D}} &=&\frac{d}{2\beta}+\frac{1}{2\sigma^2}      \thetav_{\!0}^\Tran     \left\langle\big(\J-\J \J^{-1}\J\big) \right\rangle_{\Zm} \thetav_{\!0}    +  \frac{\sigma_0^2}{2\sigma^2} \big(N-\left\langle\Tr\left[\J \J^{-1}  \right ] \right\rangle_{\Zm} \big)   -\frac{1}{2\beta} \left\langle \log\vert 2\pi\rme\sigma^2 \beta^{-1}  \J^{-1} \vert \right\rangle_{\Zm}\nonumber\\
&=&\frac{d}{2\beta} + \ \frac{\sigma_0^2}{2\sigma^2} \big(N-d \big)  -\frac{d}{2\beta}\log(2\pi\rme\sigma^2)+\frac{d}{2\beta}\log(\beta)   +\frac{1}{2\beta}\left\langle\log   \vert\J \vert \right\rangle_{\Zm}\label{eq:LR-F-cond.-aver.-ML-1}
\end{eqnarray}
Once more we put   $\Zm\to \Zm/\sqrt{d}$ where now $z_i(\mu)={\mathcal O}(1)$ for all $(i,\mu)$, then we obtain the average free energy density
\begin{eqnarray}
\frac{1}{N} \left\langle F_{\beta,\sigma^2}\left[\mathscr{D}\right] \right\rangle_{\mathscr{D}}
%
 %
 &=&\frac{1}{2} \frac{\sigma_0^2}{\sigma^2} \big(1-\zeta \big)+\frac{\zeta}{2\beta}\log\big(\frac{\beta}{2\pi\sigma^2\zeta}\big) +\frac{\zeta}{2\beta}\int\!\rmd \lambda\,\left\langle\rho_d (\lambda\vert\Zm )\right\rangle_{\Zm} \log \lambda \label{eq:LR-F-cond.-aver.-ML-density}
\end{eqnarray}
with the density of eigenvalues of the $d\times d$ empirical covariance matrix $\C=\frac{1}{N}\Zm^\Tran\Zm$, 
\begin{eqnarray}
\rho_d (\lambda \vert\Zm)&=&\frac{1}{d}\sum_{\ell=1}^d\delta\big(\lambda-\lambda_\ell(\C)\big)\label{def:eigenvalue-density}. 
\end{eqnarray}
We can similarly compute  the variance of (\ref{eq:LR-F-cond.-App}) for $\eta=0$ and $\Zm\to \Zm/\sqrt{d}$.  Firstly, we consider the energy variance   (\ref{eq:LR-Var-E-calc.-2}) which is given by 
$\Var \big(E[\mathscr{D}]\big)=\frac{\sigma_0^4}{2\sigma^4} \big(N-d\big)$, from which we deduce    
\begin{eqnarray}
\Var \big(\frac{E[\mathscr{D}]}{N}\big)&=& \frac{\sigma_0^4}{2\sigma^4} \frac{\big(1-\zeta\big)}{N}\label{eq:LR-Var-E-ML}.
\end{eqnarray}
Secondly,  the entropy variance (\ref{eq:LR-Var-S}) is given by 
\begin{eqnarray}
\Var\big(\mathrm{S}[\mathscr{D}]\big)&=&\frac{1}{4}\left\langle \log^2\left\vert  \J\right\vert \right\rangle_{\Zm}-  \frac{1}{4} \left\langle \log\left\vert  \J \right\vert \right\rangle^2_{\Zm}
=\frac{d^2}{4}\Big[\Big\langle           \big(\frac{1}{d}\sum_{\ell=1}^d\log\lambda_\ell\big)^2        \Big\rangle_{\Zm} -  \Big\langle \frac{1}{d}\sum_{\ell=1}^d\log\lambda_\ell \Big\rangle^2_{\Zm}\Big]\nonumber\\
&=&\frac{d^2}{4}\Big[\Big\langle           \Big(  \int\!\rmd \lambda\,\rho_d (\lambda\vert\Zm ) \log \lambda             \Big)^2        \Big\rangle_{\Zm}    -  \Big\langle\int\!\rmd \lambda\,\rho_d (\lambda\vert\Zm ) \log \lambda\Big\rangle^2_{\Zm}\Big]
\end{eqnarray}
and hence
\begin{eqnarray}
\Var\big(\frac{\mathrm{S}[\mathscr{D}]}{N}\big)&=&\frac{\zeta^2}{4}\Big[\Big\langle           \big(  \int\!\rmd \lambda\,\rho_d (\lambda\vert\Zm ) \log \lambda             \big)^2        \Big\rangle_{\Zm}  -  \Big\langle\int\!\rmd \lambda\,\rho_d (\lambda\vert\Zm ) \log \lambda\Big\rangle^2_{\Zm}\Big]\label{eq:LR-Var-S-ML}
\end{eqnarray}
Finally, we consider the covariance (\ref{eq:LR-Covar-E-S}) which gives us  
\begin{eqnarray}
\Covar\big( E[ \mathscr{D}], \mathrm{S}[\mathscr{D}]  \big)
&=&\frac{\sigma_0^2}{4\sigma^2}  \big(N-d \big)  \left\langle \log\left\vert\J_{\sigma^2\eta}^{-1}\right \vert      \right\rangle_{\Zm}    -\frac{\sigma_0^2}{4\sigma^2} \big(N-d \big) \left\langle\log\left\vert\J_{\sigma^2\eta}^{-1}\right \vert  \right\rangle_{\Zm} =0  \label{eq:LR-Covar-E-S-ML}
\end{eqnarray}
Using above results in identity (\ref{eq:LR-Var-F}) we obtain the desired variance of the free energy density: 
\begin{eqnarray}
 \Var \big(\frac{F_{\beta,\sigma^2}\left[\mathscr{D}\right]}{N}\big)
&=&\Var \big(\frac{E[\mathscr{D}]}{N}\big)+T^2 \Var\big(\frac{\mathrm{S}[\mathscr{D}] }{N}    \big)\nonumber\\
 &=& \frac{\sigma_0^4}{2\sigma^4} \frac{1\!-\!\zeta}{N} +   \frac{\zeta^2}{4\beta^2}\int\!\rmd\lambda\rmd\tilde{\lambda}\left[\left\langle \rho_d (\lambda\vert\Zm )\,\rho_d (\tilde{\lambda}\vert\Zm )\right\rangle_{\Zm} \!-\! \left\langle \rho_d (\lambda\vert\Zm )\right\rangle_{\Zm} \big\langle \rho_d (\tilde{\lambda}\vert\Zm )\big\rangle_{\Zm}  \right] \log (\lambda) \log (\tilde{\lambda})
 \label{eq:LR-Var-F/N-ML}
\end{eqnarray}

\subsection{Free energy of MAP inference}
Let us assume that the true parameters $\thetav_{\!0}$ are drawn randomly, with mean $\nullv$ and covariance matrix $S^2\I_d$. We may then compute the average of (\ref{eq:LR-F-cond.-aver.-calc.-2}): 
\begin{eqnarray}
 \left\langle\left\langle F_{\beta,\sigma^2}\left[\mathscr{D}\right] \right\rangle_{\mathscr{D}}\right\rangle_{\thetav_{\!0}}  
 &=&\frac{d}{2\beta}+\frac{1}{2\sigma^2}     \left\langle \left\langle \thetav_{\!0}^\Tran     \big(\J\!-\!\J \J^{-1}_{\sigma^2\eta} \J\big)  \thetav_{\!0} \right\rangle_{\thetav_{\!0}} \right\rangle_{\Zm}    +  \frac{\sigma_0^2}{2\sigma^2} \big(N-\left\langle\Tr\left[\J \J^{-1}_{\sigma^2\eta}  \right ] \right\rangle_{\Zm} \big)  -\frac{1}{2\beta} \left\langle \log\vert 2\pi\rme\sigma^2 \beta^{-1}  \J_{\sigma^2\eta}^{-1} \vert \right\rangle_{\Zm}   \nonumber\\
 &&\hspace*{-15mm} =\frac{d}{2\beta}+\frac{1}{2\sigma^2}     \left\langle     \Tr\left[ \big(\J\!-\!\J \J^{-1}_{\sigma^2\eta} \J\big)  \left\langle \thetav_{\!0}\thetav_{\!0}^\Tran \right\rangle_{\thetav_{\!0}} \right]\right\rangle_{\Zm}    + \frac{\sigma_0^2}{2\sigma^2} \big(N\!-\!\left\langle\Tr\left[\J \J^{-1}_{\sigma^2\eta}  \right ] \right\rangle_{\Zm} \big)  -\frac{1}{2\beta} \left\langle \log\vert 2\pi\rme\sigma^2 \beta^{-1}  \J_{\sigma^2\eta}^{-1} \vert \right\rangle_{\Zm}\nonumber\\
&&\hspace*{-15mm} =\frac{d}{2\beta}+\frac{S^2}{2\sigma^2}     \left\langle     \Tr\left[ \big(\J\!-\!\J \J^{-1}_{\sigma^2\eta} \J\big)   \right]\right\rangle_{\Zm}   + \frac{\sigma_0^2}{2\sigma^2} \big(N\!-\!\left\langle\Tr\left[\J \J^{-1}_{\sigma^2\eta}  \right ] \right\rangle_{\Zm} \big) -\frac{1}{2\beta} \left\langle \log\vert 2\pi\rme\sigma^2 \beta^{-1}  \J_{\sigma^2\eta}^{-1} \vert \right\rangle_{\Zm}\nonumber\\
&&\hspace*{-15mm} =\frac{d}{2\beta}+\frac{S^2}{2\zeta\sigma^2}     \left\langle     \Tr\left[ \big(\C\!-\!\C \C^{-1}_{\zeta\sigma^2\eta} \C\big)   \right]\right\rangle_{\Zm}
   +  \frac{\sigma_0^2}{2\sigma^2} \big(N\!-\!\left\langle\Tr\left[\C \C^{-1}_{\zeta\sigma^2\eta}  \right ] \right\rangle_{\Zm} \big) -\frac{1}{2\beta} \left\langle \log\vert 2\pi\rme\sigma^2 \beta^{-1}  \zeta\C_{\zeta\sigma^2\eta}^{-1} \vert \right\rangle_{\Zm}.
\label{eq:LR-F-cond.-aver.-MAP-calc.-1}
\end{eqnarray}
 In the last line  we have set $\Zm\to \Zm/\sqrt{d}$, with $z_i(\mu)={\mathcal O}(1)$ for all $(i,\mu)$,   and used  $\J=\C/\zeta$ and $\C^{-1}_{\sigma^2\eta} = \zeta\C^{-1}_{\zeta\sigma^2\eta}$,  where $\C=\Zm^\Tran\Zm/N$.
Let us consider the matrix product $\C \C^{-1}_{\eta} =\C(\C\!+\!\eta\I_d)^{-1}=((\C\!+\!\eta\I_d)\C^{-1})^{-1}=(\I_d\!+\!\eta\C^{-1})^{-1}$, giving
\begin{eqnarray}
 \Tr\left[\C \C^{-1}_{\eta}  \right ]&=&\Tr\left[\big(\I_d+\eta\C^{-1}\big)^{-1}\right]=\sum_{\ell=1}^d1/\lambda_\ell\big( \I_d+\eta\C^{-1}   \big)=\sum_{\ell=1}^d1/\big( 1+\eta\lambda_\ell\big(\C^{-1}   \big)\big)\nonumber\\
  &=&\sum_{\ell=1}^d1/\big( 1+\eta\lambda^{-1}_\ell\big(\C   \big)\big)
=\sum_{\ell=1}^d\lambda_\ell\big(\C   \big)/\big( \lambda_\ell\big(\C   \big)+\eta\big). 
 \end{eqnarray}
Similarly we can write
\begin{eqnarray}
\C^2 \C^{-1}_{\eta}& =&\C^2 \big(\C+\eta\I_d\big)^{-1}=\C\big(\I_d+\eta\C^{-1}\big)^{-1}=\big(\big(\I_d+\eta\C^{-1}\big)\C^{-1}\big)^{-1}=\big(\C^{-1}+\eta\C^{-2}\big)^{-1}
 \end{eqnarray}
 The matrices $\C^{2} $ and $   \big(\C+\eta\I_d\big)^{-1} $ obviously {commute},  so 
\begin{eqnarray}
 \Tr\left[ \C \C^{-1}_{\eta} \C   \right]&=&\Tr\left[ \C^2 \C^{-1}_{\eta}  \right]=\sum_{\ell=1}^d\lambda_\ell\big( \C^{2}  \big)/\lambda_\ell\big(\C_{\eta}     \big)=\sum_{\ell=1}^d\lambda^2_\ell\big( \C\big)/\lambda_\ell\big(\C_{\eta}     \big).  
 \end{eqnarray}
 Now
  $\C_{\eta} =\C+\eta\I_d=\C(\I_d+\eta\C^{-1})=(\I_d+\eta\C^{-1})\C$, 
 and hence 
\begin{eqnarray}
 \lambda_\ell\big(\C_{\eta}     \big)&=&\lambda_\ell\big(\C\big) \lambda_\ell\big(  \I_d+\eta\C^{-1}  \big)=\lambda_\ell\big(\C\big)   \big(1+\eta\lambda_\ell\big(\C^{-1}  \big)\big)=\lambda_\ell\big(\C\big)+\eta. 
 \end{eqnarray}
 Thus $ \Tr\left[ \C \C^{-1}_{\eta} \C   \right]=\sum_{\ell=1}^d\lambda^2_\ell\big( \C\big)/\big(\lambda_\ell\big(\C\big)+\eta\big)$.
Finally, we consider the inverse 
\begin{eqnarray}
\C^{-1}_{\eta} &=&\big(\C+\eta\I_d\big)^{-1}=\big(\I_d+\eta\C^{-1}\big)^{-1}\C^{-1}=\C^{-1}\big(\I_d+\eta\C^{-1}\big)^{-1},
 \end{eqnarray}
  The matrices $\C^{-1}$ and $\big(\I_d+\eta\C^{-1}\big)^{-1}$ obviously {commute}, and hence the $\ell$-th eigenvalue of $\C_\eta^{-1}$ is given by
\begin{eqnarray}
  \lambda_\ell\big( \J^{-1}_{\eta}   \big)&=&1/(\lambda_\ell\big(\J  \big)+\eta).
 \end{eqnarray}
Using the above results on the relevant matrices  in (\ref{eq:LR-F-cond.-aver.-MAP-calc.-1}) allows us to compute the average free energy:
\begin{eqnarray}
 \left\langle\left\langle\frac{1}{N} F_{\beta,\sigma^2}\left[\mathscr{D}\right] \right\rangle_{\mathscr{D}}\right\rangle_{\thetav_{\!0}}  
&=&\frac{\zeta}{2\beta}+\frac{S^2}{2\zeta\sigma^2}  \frac{1}{N}   \left\langle     \Tr\left[ \big(\J-\J \J^{-1}_{\zeta\sigma^2\eta} \J\big)   \right]\right\rangle_{\Zm}+ \frac{\sigma_0^2}{2\sigma^2} \big(1- \frac{1}{N}\left\langle\Tr\left[\J \J^{-1}_{\zeta\sigma^2\eta}  \right ] \right\rangle_{\Zm} \big)  \nonumber\\
&&~~~~~~~~~~~~~~~~~~~~~-\frac{1}{2\beta}  \frac{1}{N}\left\langle \log\vert 2\pi\rme\sigma^2 \beta^{-1}  \zeta\J_{\zeta\sigma^2\eta}^{-1} \vert \right\rangle_{\Zm}\nonumber\\
&=&\frac{\zeta}{2\beta}+\frac{S^2}{2\zeta\sigma^2}  \frac{d}{Nd}      \sum_{\ell=1}^d\left\langle \lambda_\ell\big(\J\big)-\frac{\lambda^2_\ell\big( \J\big)}{\lambda_\ell\big(\J\big)+\zeta\sigma^2\eta}
 \right\rangle_{\Zm} +  \frac{\sigma_0^2}{2\sigma^2} \big(1- \frac{d}{Nd}\left\langle  \sum_{\ell=1}^d \frac{\lambda_\ell\big( \J\big)}{\lambda_\ell\big(\J\big)+\zeta\sigma^2\eta}
  \right\rangle_{\Zm} \big)  \nonumber\\
&&~~~~~~~~~~-\frac{\zeta}{2\beta}   \log\big( 2\pi\rme\sigma^2 \beta^{-1}  \zeta \big)+\frac{\zeta}{2\beta} \frac{1}{d}\left\langle\sum_{\ell=1}^d \log \big(\lambda_\ell\big(\J  \big)+ \zeta\sigma^2\eta   \big)  \right\rangle_{\Zm}\nonumber\\
 &=&\frac{\zeta}{2\beta}+\frac{S^2   \zeta\eta  }{2}   \int\!\rmd \lambda~ \left\langle     \rho_d (\lambda\vert\Zm ) \right\rangle_{\Zm}  \frac{  \lambda        }{\lambda+\zeta\sigma^2\eta} 
+  \frac{\sigma_0^2}{2\sigma^2} \Big(1- \zeta   \int\! \rmd \lambda~  \left\langle     \rho_d (\lambda\vert\Zm ) \right\rangle_{\Zm}\frac{\lambda}{\lambda+\zeta\sigma^2\eta}\,
  \Big)  \nonumber\\
&&~~~~~~~~~~~~~~~-\frac{\zeta}{2\beta}   \log\big( 2\pi\rme\sigma^2 \beta^{-1}  \zeta \big)+\frac{\zeta}{2\beta}  \int\! \rmd \lambda~ \left\langle     \rho_d (\lambda\vert\Zm ) \right\rangle_{\Zm}\log \big(\lambda+ \zeta\sigma^2\eta   \big)     \label{eq:LR-F-cond.-aver.-MAP-calc.-2}
\end{eqnarray}
and hence  the average free energy density  is given by 
\begin{eqnarray}
 \left\langle\left\langle\frac{1}{N} F_{\beta,\sigma^2}\left[\mathscr{D}\right] \right\rangle_{\mathscr{D}}\right\rangle_{\thetav_{\!0}}   &=&\frac{\zeta}{2\beta}+\frac{S^2   \zeta\eta  }{2}   \int\!\rmd \lambda~ \left\langle     \rho_d (\lambda\vert\Zm ) \right\rangle_{\Zm}  \frac{  \lambda        }{\lambda+\zeta\sigma^2\eta} 
+  \frac{\sigma_0^2}{2\sigma^2} \Big(1- \zeta   \int\! \rmd \lambda~  \left\langle     \rho_d (\lambda\vert\Zm ) \right\rangle_{\Zm}\frac{\lambda}{\lambda+\zeta\sigma^2\eta}
  \Big)  \nonumber\\
&&-\frac{\zeta}{2\beta}   \log\big( 2\pi\rme\sigma^2 \beta^{-1}  \zeta \big)
+\frac{\zeta}{2\beta}  \int\! \rmd \lambda~ \left\langle     \rho_d (\lambda\vert\Zm ) \right\rangle_{\Zm}\log \big(\lambda+ \zeta\sigma^2\eta   \big)     \label{eq:LR-F/N-cond.-aver.-MAP}
\end{eqnarray}
In the derivation of the above results we used the following simple eigenvalue identities 
\begin{eqnarray}
 \lambda_\ell\big( \C^{-1}_{\eta}  \big)=\frac{1}{\lambda_\ell\big(\C  \big)+\eta},~~~~~
 \lambda_\ell\big(\C \C^{-1}_{\eta}  \big)= \frac{\lambda_\ell\big(\C   \big)}{    \lambda_\ell\big(\C   \big)+\eta},~~~~~
   \lambda_\ell\big(\C\!-\!\C \C^{-1}_{\eta} \C   \big) = \frac{\eta\lambda_\ell\big(\C   \big)}{    \lambda_\ell\big(\C   \big)+\eta} \label{eq:LR-eigen-identities-1},
\end{eqnarray}
 from one also obtains
\begin{eqnarray}
  \Tr\left[ \J-\J \J^{-1}_{\eta} \J   \right] =d\eta \int\!\rmd \lambda~\rho_d (\lambda\vert\Zm )  \frac{  \lambda        }{\lambda+\eta} ,~~~~~~
 \Tr\left[\J \J^{-1}_{\eta}  \right ] =d \int\! \rmd \lambda~    \rho_d (\lambda\vert\Zm )  \frac{  \lambda        }{\lambda+\eta} 
   \label{eq:LR-Tr-identities-1}
\end{eqnarray}

Finally, we  compute the  variance of  the free energy  (\ref{eq:LR-F-cond.-App}). First, we consider the energy variance   (\ref{eq:LR-Var-E-calc.-2}) which is 
\begin{eqnarray}
4\sigma^4\Var \big(E[\mathscr{D}]\big) &=&\left\langle\left\langle\big(\thetav_{\!0}^\Tran   \big(\J\!-\!\J \J^{-1}_{\sigma^2\eta} \J\big) \thetav_{\!0}    \big)^2 \right\rangle_{\Zm}\right\rangle_{\thetav_{\!0}}        -  \left\langle\left\langle\thetav_{\!0}^\Tran   \big(\J\!-\!\J \J^{-1}_{\sigma^2\eta} \J\big) \thetav_{\!0} \right\rangle_{\thetav_{\!0}}   \right\rangle^2_{\Zm} +4\left\langle\sigma_0^2  \thetav_{\!0}^\Tran  \left\langle \J\big( \I_d -  \J^{-1}_{\sigma^2\eta}\J\big)^2  \right\rangle_{\Zm}    \thetav_{\!0}\right\rangle_{\thetav_{\!0}}\nonumber\\
&&+\sigma_0^4\Bigg[ \left\langle \big(N\!-\!\Tr\left[\J \J^{-1}_{\sigma^2\eta}\right]\big)^2  \right\rangle_{\Zm}\!+ 2 \left\langle  \big(N\!-\!2\Tr\left[\J \J^{-1}_{\sigma^2\eta} \right]  +\Tr \left[\big(\J \J^{-1}_{\sigma^2\eta}\big)^2\right] \big)   \right\rangle_{\Zm}\!- \left\langle  \big(N\!-\!\Tr\left[\J \J^{-1}_{\sigma^2\eta}\right]\big)    \right\rangle^2_{\Zm}\Bigg]\nonumber\\
 &&+ 2\sigma_0^2\Bigg[\Big\langle  \Big\langle  \thetav_{\!0}^\Tran   \big(\J\!-\!\J \J^{-1}_{\sigma^2\eta} \J\big) \thetav_{\!0}      \big(N\!-\!\Tr[ \J\J^{-1}_{\sigma^2\eta}]\big)\Big\rangle_{\!\thetav_{\!0}}\Big\rangle_{\!\Zm}\!- \Big\langle \Big\langle \thetav_{\!0}^\Tran   \big(\J\!-\!\J \J^{-1}_{\sigma^2\eta} \J\big) \thetav_{\!0}  \Big\rangle_{\!\thetav_{\!0}}   \Big\rangle_{\Zm}  \Big\langle \big(N-\Tr[ \J\J^{-1}_{\sigma^2\eta}]\big)    \Big\rangle_{\Zm}\Bigg]
 \nonumber\\
 \nonumber\\
 &&\hspace*{-15mm} = S^4\left\langle  
   \Tr^2\!\left[  \big(\J\!-\!\J \J^{-1}_{\sigma^2\eta} \J\big)  \right]\!+ 2\Tr\left[ \big(\J\!-\!\J \J^{-1}_{\sigma^2\eta} \J\big)^2\right]
 \right\rangle_{\Zm}\!-  S^4\left\langle  \Tr \left[\J\!-\!\J \J^{-1}_{\sigma^2\eta} \J\right]   \right\rangle^2_{\Zm}\!+4\sigma_0^2S^2    \left\langle  \Tr\left[\J\big( \I_d \!-\!  \J^{-1}_{\sigma^2\eta}\J\big)^2 \right] \right\rangle_{\Zm}    \nonumber\\
&&+\sigma_0^4\Bigg[ \left\langle \big(N\!-\!\Tr\left[\J \J^{-1}_{\sigma^2\eta}\right]\big)^2  \right\rangle_{\Zm}+ 2 \left\langle  \big(N\!-\!2\Tr\left[\J \J^{-1}_{\sigma^2\eta} \right]  +\Tr \left[\big(\J \J^{-1}_{\sigma^2\eta}\big)^2\right] \big)   \right\rangle_{\Zm}
- \left\langle  \big(N\!-\!\Tr\left[\J \J^{-1}_{\sigma^2\eta}\right]\big)    \right\rangle^2_{\Zm}\Bigg]\nonumber\\
 &&+ 2\sigma_0^2S^2\Bigg(\left\langle    \Tr \left[\J\!-\!\J \J^{-1}_{\sigma^2\eta} \J\right]       \big(N\!-\!\Tr[ \J\J^{-1}_{\sigma^2\eta}]\big)\right\rangle_{\Zm}
  - \left\langle  \Tr \left[\J\!-\!\J \J^{-1}_{\sigma^2\eta} \J\right]    \right\rangle_{\Zm}  \left\langle \big(N\!-\!\Tr[ \J\J^{-1}_{\sigma^2\eta}]\big)    \right\rangle_{\Zm}\Bigg) \label{eq:LR-Var-E-calc.-3}
\end{eqnarray}
where we have  used (\ref{eq:Wicks-theorem}).  Hence  for $\Zm\to \Zm/\sqrt{d}$ with $z_i(\mu)={\mathcal O}(1)$, we have 
\begin{eqnarray}
4\sigma^4\Var \big(\frac{E[\mathscr{D}] }{N}   \big)
 &=& \frac{S^4}{\zeta^2N^2}\left\langle  
   \Tr^2\left[  \big(\C\!-\!\C \C^{-1}_{\zeta\sigma^2\eta} \C\big)  \right]+ 2\Tr\left[ \big(\C\!-\!\C \C^{-1}_{\zeta\sigma^2\eta} \C\big)^2\right]
 %
 \right\rangle_{\Zm}-  \frac{S^4}{\zeta^2N^2}\left\langle 
  \Tr \left[\C\!-\!\C \C^{-1}_{\zeta\sigma^2\eta} \C\right] 
    \right\rangle^2_{\Zm} \nonumber\\
&&+\frac{4\sigma_0^2S^2}{\zeta N^2}    \left\langle  \Tr\left[\C\big( \I_d\! -\!  \C^{-1}_{\zeta\sigma^2\eta}\C\big)^2 \right] \right\rangle_{\Zm}    \nonumber\\
&&\hspace*{-15mm} +\frac{\sigma_0^4}{N^2}\Bigg[ \left\langle \big(N\!-\!\Tr\left[\C \C^{-1}_{\zeta\sigma^2\eta}\right]\big)^2  \right\rangle_{\Zm}\!+ 2 \left\langle  \big(N\!-\!2\Tr\left[\C \C^{-1}_{\zeta\sigma^2\eta} \right]  +\Tr \left[\big(\C \C^{-1}_{\zeta\sigma^2\eta}\big)^2\right] \big)   \right\rangle_{\Zm}\!
- \left\langle  \big(N\!-\!\Tr\left[\C \C^{-1}_{\zeta\sigma^2\eta}\right]\big)    \right\rangle^2_{\Zm}\Bigg]\nonumber\\
 &&\hspace*{-15mm} + 2\sigma_0^2S^2\frac{1}{\zeta N^2}\Bigg[\left\langle    \Tr \left[\C\!-\!\C \C^{-1}_{\zeta\sigma^2\eta} \C\right]       \big(N\!-\!\Tr[ \C\C^{-1}_{\zeta\sigma^2\eta}]\big)\right\rangle_{\Zm}- \left\langle  \Tr \left[\C\!-\!\C \C^{-1}_{\zeta\sigma^2\eta} \C\right]    \right\rangle_{\Zm}  \left\langle \big(N\!-\!\Tr[ \C\C^{-1}_{\zeta\sigma^2\eta}]\big)    \right\rangle_{\Zm}\Bigg] \nonumber\\
 &=& 
 \frac{S^4}{\zeta^2N^2}\left\langle  
  \left[d\zeta\sigma^2\eta\int\! \rmd\lambda ~  \frac{ \rho_d (\lambda\vert\Zm )  \lambda        }{\lambda+\zeta\sigma^2\eta} \right]^2
   + 2d (\zeta\sigma^2\eta)^2\int\!  \rmd\lambda ~\rho_d (\lambda\vert\Zm )  \Big(\frac{  \lambda        }{\lambda\!+\!\zeta\sigma^2\eta}\Big)^2 
 \right\rangle_{\Zm}\nonumber\\
 &&-  \frac{S^4}{\zeta^2N^2}\left\langle  
  d\zeta\sigma^2\eta\int\!\rmd\lambda   ~  \frac{\rho_d (\lambda\vert\Zm )  \lambda        }{\lambda\!+\!\zeta\sigma^2\eta} 
   \right\rangle^2_{\Zm} +\frac{4\sigma_0^2S^2}{\zeta N^2}    \left\langle 
d(\zeta\sigma^2\eta)^2\int\!  \rmd\lambda~\frac{   \rho_d (\lambda\vert\Zm )  \lambda        }{\big(\lambda\!+\!\zeta\sigma^2\eta\big)^2} 
  \right\rangle_{\Zm}    \nonumber\\
&&+\sigma_0^4\Bigg[ \left\langle \Big(1-\frac{d}{N}
\int\! \rmd\lambda~    \frac{\rho_d (\lambda\vert\Zm )   \lambda        }{\lambda\!+\!\zeta\sigma^2\eta} 
\Big)^2  \right\rangle_{\Zm}
- \left\langle  1-\frac{d}{N}
 \int\!\rmd\lambda~     \frac{\rho_d (\lambda\vert\Zm )   \lambda        }{\lambda\!+\!\zeta\sigma^2\eta} 
 \right\rangle^2_{\Zm}\nonumber\\
&&\hspace*{10mm}+ \frac{2}{N^2} \left\langle  N-2
d\int\!   \rmd\lambda~ \frac{ \rho_d (\lambda\vert\Zm )   \lambda        }{\lambda\!+\!\zeta\sigma^2\eta}
  +d\int\! \rmd\lambda~    \rho_d (\lambda\vert\Zm )  \big(\frac{  \lambda        }{\lambda\!+\!\zeta\sigma^2\eta}\big)^2 
  \right\rangle_{\Zm}\Bigg]\nonumber\\
 &&+ 2\sigma_0^2S^2\frac{1}{\zeta}\Bigg[\left\langle   \Big( \frac{d}{N}
  \zeta\sigma^2\eta\int\!\rmd\lambda~     \frac{\rho_d (\lambda\vert\Zm )  \lambda        }{\lambda\!+\!\zeta\sigma^2\eta} 
\Big)
       \Big(1-\frac{d}{N}
 \int\! \rmd\lambda~ \frac{   \rho_d (\lambda\vert\Zm )   \lambda        }{\lambda\!+\!\zeta\sigma^2\eta} 
 \Big)\right\rangle_{\Zm}\nonumber\\
 &&~~~~~~~~~~ - \left\langle  \frac{d}{N}
  \zeta\sigma^2\eta\int\! \rmd\lambda~ \frac{   \rho_d (\lambda\vert\Zm )   \lambda        }{\lambda\!+\!\zeta\sigma^2\eta} 
    \right\rangle_{\Zm}  \left\langle 1-\frac{d}{N}
 \int\! \rmd\lambda~ \frac{   \rho_d (\lambda\vert\Zm )   \lambda        }{\lambda\!+\!\zeta\sigma^2\eta} 
  \right\rangle_{\Zm}\Bigg] \nonumber\\
 &=& \zeta^2\big(S^4\sigma^4\eta^2 \! +\!\sigma_0^4 \!-\! 2\sigma_0^2S^2 \sigma^2\eta\big)\Bigg[ \left\langle  
 \Big(\int\! \rmd\lambda~ \frac{   \rho_d (\lambda\vert\Zm )  \lambda        }{\lambda\!+\!\zeta\sigma^2\eta} \Big)^2\right\rangle_{\Zm}
  - \left\langle  
  \int\! \rmd\lambda~\frac{    \rho_d (\lambda\vert\Zm )  \lambda        }{\lambda\!+\!\zeta\sigma^2\eta} 
   \right\rangle^2_{\Zm} \Bigg]+O(1/N)
 \label{eq:LR-Var-E/N-calc.-1}
\end{eqnarray}
Hence
\begin{eqnarray}
\Var \big(\frac{ E[\mathscr{D}]}{N}    \big)
&=&\frac{ \zeta^2\big(S^4\sigma^4\eta^2  +\sigma_0^4 - 2\sigma_0^2S^2 \sigma^2\eta\big) }{4\sigma^4}     \Bigg[ \left\langle  
 \Big(\int\! \rmd\lambda~  \frac{  \rho_d (\lambda\vert\Zm )   \lambda        }{\lambda\!+\!\zeta\sigma^2\eta} \rmd \lambda\Big)^2\right\rangle_{\!\Zm} \!    - \left\langle  
  \int\!   \rmd\lambda~  \frac{\rho_d (\lambda\vert\Zm )  \lambda        }{\lambda\!+\!\zeta\sigma^2\eta} \rmd \lambda
   \right\rangle^2_{\!\Zm} \Bigg]+{\mathcal O}(\frac{1}{N}).~~~~
 \label{eq:LR-Var-E/N-MAP}. 
\end{eqnarray}
Furthermore,  the covariance (\ref{eq:LR-Covar-E-S}) becomes
\begin{eqnarray}
4\sigma^2\Covar\big( E[ \mathscr{D}], \mathrm{S}[\mathscr{D}]  \big)
&=& \Covar\Big(\tv^\Tran\big(\I_N-\Zm \J^{-1}_{\sigma^2\eta} \Zm^{T}\big)\tv  ,\log\left\vert\J_{\sigma^2\eta}^{-1}\right \vert          \Big)\nonumber\\
%
%
&=&\Big\langle  \Big[ \left\langle \thetav_{\!0}^\Tran    \big(\J-\J \J^{-1}_{\sigma^2\eta} \J\big) \thetav_{\!0}\right\rangle_{\thetav_{\!0}} + \sigma_0^2 \big(N-\Tr\left[\J \J^{-1}_{\sigma^2\eta}  \right ] \big) \Big]  \log\left\vert\J_{\sigma^2\eta}^{-1}\right \vert      \Big\rangle_{\Zm}\nonumber\\
&&-\Big\langle\Big[ \left\langle \thetav_{\!0}^\Tran    \big(\J-\J \J^{-1}_{\sigma^2\eta} \J\big) \thetav_{\!0}\right\rangle_{\thetav_{\!0}} + \sigma_0^2 \big(N-\Tr\left[\J \J^{-1}_{\sigma^2\eta}  \right ] \big) \Big]  \Big\rangle_{\Zm}\left\langle\log\left\vert\J_{\sigma^2\eta}^{-1}\right \vert  \right\rangle_{\Zm}   \nonumber\\
 &=&\Big\langle  \Big[  S^2\Tr\left[\J-\J \J^{-1}_{\sigma^2\eta} \J\right]   + \sigma_0^2 \big(N-\Tr\left[\J \J^{-1}_{\sigma^2\eta}  \right ] \big) \Big]  \log\left\vert\J_{\sigma^2\eta}^{-1}\right \vert      \Big\rangle_{\Zm}\nonumber\\
&&-\Big\langle\Big[   S^2\Tr \left[\J-\J \J^{-1}_{\sigma^2\eta} \J\right]  + \sigma_0^2 \big(N-\Tr\left[\J \J^{-1}_{\sigma^2\eta}  \right ] \big) \Big]  \Big\rangle_{\Zm}\left\langle\log\left\vert\J_{\sigma^2\eta}^{-1}\right \vert  \right\rangle_{\Zm} 
  \label{eq:LR-Covar-E-S-MAP-calc.-1}
\end{eqnarray}
From this, upon setting  $\Zm\to \Zm/\sqrt{d}$ with $z_i(\mu)={\mathcal O}(1)$ for all $(i,\mu)$ then follows the  result
\begin{eqnarray}
4\sigma^2\Covar\big( E[ \mathscr{D}]/N, \mathrm{S}[\mathscr{D}]/N  \big)
 &=&\Big\langle  \Big[  S^2\frac{1}{\zeta N}\Tr\left[\J-\J \J^{-1}_{\zeta\sigma^2\eta} \J\right]   + \sigma_0^2 \big(1-\frac{1}{N}\Tr\left[\J \J^{-1}_{\zeta\sigma^2\eta}  \right ] \big) \Big]  \frac{1}{N}\log\left\vert\zeta\J_{\zeta\sigma^2\eta}^{-1}\right \vert      \Big\rangle_{\Zm}\nonumber\\
&&-\Big\langle\Big[   S^2\frac{1}{\zeta N}\Tr \left[\J-\J \J^{-1}_{\zeta\sigma^2\eta} \J\right]  + \sigma_0^2 \big(1-\frac{1}{N}\Tr\left[\J \J^{-1}_{\zeta\sigma^2\eta}  \right ] \big) \Big]  \Big\rangle_{\Zm}\frac{1}{N}\left\langle\log\left\vert\zeta\J_{\zeta\sigma^2\eta}^{-1}\right \vert  \right\rangle_{\Zm}  \nonumber\\
 &=& \zeta  \int\!\rmd\lambda\rmd\tilde{\lambda} ~\left[   \langle \rho_d (\lambda\vert\Zm ) \rangle_{\Zm} \langle\rho_d (\tilde{\lambda}\vert\Zm ) \rangle_{\Zm}- 
 \langle\rho_d (\lambda\vert\Zm )    \rho_d (\tilde{\lambda}\vert\Zm ) \rangle_{\Zm}\right]\nonumber\\
 &&~~~~~~\times\Big[ 
   \frac{  S^2\lambda   \zeta\sigma^2\eta      }{\lambda\!+\!\zeta\sigma^2\eta} 
   + \sigma_0^2 \big(1\!-\!
    \frac{  \lambda   \zeta     }{\lambda\!+\!\zeta\sigma^2\eta} 
 \big) \Big]\log\big(\tilde{\lambda}\!+\! \zeta\sigma^2\eta\big)
   \label{eq:LR-Covar-E-S-MAP}
\end{eqnarray}
Finally,  the entropy variance  (\ref{eq:LR-Var-S}) for $\Zm\to\Zm/\sqrt{d}$  is given by 
\begin{eqnarray}
\Var\big(\frac{\mathrm{S}[\mathscr{D}]}{N}\big) &=&\frac{\zeta^2}{4} \!\int    \!\rmd\lambda\rmd\tilde{\lambda} ~ \left[\langle \rho_d (\lambda\vert\Zm )\rho_d (\tilde{\lambda}\vert\Zm )\rangle_{\Zm} \!-\! \langle \rho_d (\lambda\vert\Zm )\rangle_{\Zm} \langle \rho_d (\tilde{\lambda}\vert\Zm )\rangle_{\Zm}  \right]
 \log(\lambda\!+\! \zeta\sigma^2\eta)  \log(\tilde{\lambda}\!+\! \zeta\sigma^2\eta). \label{eq:LR-Var-S-MAP}
\end{eqnarray}
Using all of the above results in (\ref{eq:LR-Var-F}) we finally obtain the variance of the free energy density:
\begin{eqnarray}
 \Var \Big(  \frac{ F_{\beta,\sigma^2}\left[\mathscr{D}\right] }      {N}     \Big)    
&=&\int\!\rmd\lambda\rmd\tilde{\lambda}~\left[\langle \rho_d (\lambda\vert\Zm )\rho_d (\tilde{\lambda}\vert\Zm )\rangle_{\Zm} -\langle \rho_d (\lambda\vert\Zm )\rangle_{\Zm} \langle \rho_d (\tilde{\lambda}\vert\Zm )
\rangle_{\Zm}  \right]\nonumber\\
&&\times\Bigg[\frac{ \zeta^2}{4\sigma^4}\big(S^4\sigma^4\eta^2  \!+\!\sigma_0^4 \!-\! 2\sigma_0^2S^2 \sigma^2\eta\big)\frac{  \lambda        }{\lambda\!+\!\zeta\sigma^2\eta} \frac{  \tilde{\lambda}        }{\tilde{\lambda}\!+\!\zeta\sigma^2\eta}+\frac{1}{4}T^2\zeta^2 \log\big(\lambda\!+\! \zeta\sigma^2\eta\big) \log\big(\tilde{\lambda}\!+\! \zeta\sigma^2\eta\big) \nonumber\\
&&~~~~~~~~~~~~~~~- \frac{T\zeta}{2\sigma^2} \Big( \frac{  S^2 \lambda   \zeta\sigma^2\eta      }{\lambda\!+\!\zeta\sigma^2\eta} 
   + \sigma_0^2 \Big(1\!-\! \frac{  \lambda   \zeta     }{\lambda\!+\!\zeta\sigma^2\eta} \Big) \Big)\log\big(\tilde{\lambda}+ \zeta\sigma^2\eta\big)\Bigg]
 +{\mathcal O}(1/N) \label{eq:LR-Var-F/N-MAP}.
\end{eqnarray}

\section{Self-averaging  of the MAP estimator of $\hat{\sigma}^2$ \label{section:sigma-self-aver}}

Let us consider the equation 
\begin{eqnarray}
\sigma^2&=& \frac{\beta}{(\beta-\zeta)}\frac{1}{N}\left\vert\left\vert \tv -\Zm \thetaest\left[  \mathscr{D}\right]   \right\vert\right\vert^2 -\frac{\sigma^4\eta}{(\beta-\zeta)}\frac{1}{N}\Tr\left[  \J^{-1}_{\sigma^2\eta}  \right]      + \frac{2\sigma^4\beta}{(\beta-\zeta)N}\frac{\partial}{\partial\sigma^2}\log P(\sigma^2)  \label{eq:sigma-est-beta-1}.
\end{eqnarray}
Upon inserting the 
 MAP estimator (\ref{eq:LR-theta-MAP-App}) and the short-hand  $\sigma^2=v$  we  obtain 
\begin{eqnarray}
v
&=& \frac{\beta}{(\beta-\zeta)}\frac{1}{N}     \tv^\Tran  \left( \I_N -\Zm    \J_{v \eta}^{-1}  \Zm^\Tran\right)^2\tv                 -\frac{v^2\eta}{(\beta-\zeta)}\frac{1}{N}\Tr\left[  \J^{-1}_{v\eta}  \right]      + \frac{2v^2\beta}{(\beta-\zeta)N}\frac{\partial}{\partial v}\log P(v)\label{eq:sigma-est-beta-2}.
\end{eqnarray}
To solve this equation for $v$ we define the following recursion, with the short-hand   $\partial_{v}   \equiv  \frac{\partial}{\partial v}$:
\begin{eqnarray}
%
%
v_{t+1}&=& \frac{\beta}{(\beta-\zeta)}\frac{1}{N}     \tv^\Tran  \left( \I_N -\Zm    \J_{v_t \eta}^{-1}  \Zm^\Tran\right)^2\tv                 -\frac{v_t^2\eta}{(\beta-\zeta)}\frac{1}{N}\Tr\left[  \J^{-1}_{v_t\eta}  \right]      + \frac{2v_t^2\beta}{(\beta-\zeta)N}  \partial_{v}  \log P(v)\vert_{v=v_t} \label{eq:sigma-est-beta-recursion},
\end{eqnarray}
Since $\tv=\Zm  \thetav_0 + \epsilonv$ this recursion, of which the desired estimator is the fixed-point,  has the general form 
\begin{eqnarray}
%
%
v_{t+1}&=&  \Psi \left[v_{t}\vert \Zm,  \thetav_{0},  \epsilonv\, \right]  \label{eq:random-recursion}. 
\end{eqnarray}
Thus  for any choice of $\{ \Zm, \thetav_{0}, \epsilonv\}$, i.e. which play the role of `disorder',  the function $\Psi$   is a random  non-linear operator  acting on $v_{t}$.  If the   initial   value $v_{0}$ is  \emph{independent} of the disorder, then the next value $v_{1}$  is independent from a \emph{particular} realisation of disorder, i.e. $v_{1}$  is \emph{self-averaging},  as soon as the operator $\Psi$ is self-averaging, i.e. if 
\begin{eqnarray}
\lim_{(N,d)\to\infty}\left\langle\Psi^2 \left[v_{0}\vert \Zm,  \thetav_{0},  \epsilonv\, \right] \right\rangle_{\Zm, \thetav_{0}, \epsilonv}-\left\langle\Psi^2 \left[v_{0}\vert \Zm,  \thetav_{0},  \epsilonv\, \right] \right\rangle^2_{\Zm, \thetav_{0}, \epsilonv}=0 \label{eq:nl-operator-self-aver.}. 
\end{eqnarray}
By induction, all $v_{t}$ with $t\geq 1$ will then be self-averaging, and (\ref{eq:random-recursion}) can for $(N,d)\to\infty$ be  replaced by the following deterministic map, whose fixed-point will be the asymptotic estimator  $\hat{\sigma}^2$ that is then guaranteed to be self-averaging:
\begin{eqnarray}
v_{t+1}&=& \left\langle \Psi \left[v_{t}\vert \Zm,  \thetav_{0},  \epsilonv\, \right] \right\rangle_{\Zm, \thetav_{0}, \epsilonv} \label{eq:aver-recursion}
\end{eqnarray}

To prove the self-averaging property of $\Psi$ we assume that the true parameters $\thetav_{\!0}$ and noise $\epsilonv$ have mean $\nullv$ and the covariance  matrices $S^2\I_d$ and $\sigma_0^2\I_N$, respectively.  Let us first consider the average 
\begin{eqnarray}
 \left\langle \Psi \left[v_{0}\vert \Zm,  \thetav_{0},  \epsilonv\, \right] \right\rangle_{\Zm, \thetav_{0}, \epsilonv} 
&=& \frac{\beta}{N(\beta\!-\!\zeta)}    \left\langle  \tv^\Tran  \left( \I_N\! -\!\Zm    \J_{v_0\eta}^{-1}  \Zm^\Tran\right)^2\tv \right\rangle_{\Zm, \thetav_{0}, \epsilonv}                 -\frac{v_0^2\eta}{N(\beta\!-\!\zeta)}\left\langle\Tr\left[  \J^{-1}_{v_0\eta}  \right]  \right\rangle_{\Zm } + \frac{2v_0^2\beta}{(\beta\!-\!\zeta)N}  \partial_{v}  \log P(v)\vert_{v=v_0}  \nonumber\\
&=&     \frac{\beta\sigma_0^2}{N(\beta\!-\!\zeta)} \left\langle  \Tr\left[ \left( \I_N \!-\!\Zm    \J_{v_0\eta}^{-1}  \Zm^\Tran\right)^2\right] \right\rangle_{\Zm} + \frac{\beta}{N(\beta-\zeta)} \left\langle\thetav^\Tran_{0} \Zm^\Tran \left( \I_N \!-\!\Zm    \J_{v_0\eta}^{-1}  \Zm^\Tran\right)^2\Zm\thetav_{0}  \right\rangle_{\Zm, \thetav_{0}}\nonumber\\
 &&~~~~~~~~~~~~~~~~~~~~~~~~~~~~~~~~~~~~~~~~~~~-\frac{v_0^2\eta}{N(\beta\!-\!\zeta)} \left\langle\Tr\left[  \J^{-1}_{v_0\eta}  \right]  \right\rangle_{\Zm }  + \frac{2v_0^2\beta}{(\beta\!-\!\zeta)N}  \partial_{v}  \log P(v)\vert_{v=v_0}  \nonumber\\
&=&     \frac{\beta}{N(\beta\!-\!\zeta)} \sigma_0^2\left\langle  \Tr\left[ \left( \I_N \!-\!\Zm    \J_{v_0\eta}^{-1}  \Zm^\Tran\right)^2\right] \right\rangle_{\Zm} + \frac{\beta}{N(\beta\!-\!\zeta)} S^2\left\langle   \Tr\left[ \J \left( \I_d \!-\!   \J_{v_0\eta}^{-1}  \J\right)^2\right]  \right\rangle_{\Zm } \nonumber\\
 &&~~~~~~~~~~~~~~~~~~~~~~~~~~~~~~~~~~~~~~~~~~~-\frac{v_0^2\eta}{N(\beta\!-\!\zeta)} \left\langle\Tr\left[  \J^{-1}_{v_0\eta}  \right]  \right\rangle_{\Zm }  + \frac{2v_0^2\beta}{(\beta\!-\!\zeta)N}  \partial_{v}  \log P(v)\vert_{v=v_0}  \nonumber\\
 &=&     \frac{\beta  \sigma_0^2}{\beta-\zeta} \left( 1\!-\!\zeta \!+\! \frac{1}{N} \left\langle  \Tr\big[ \big( \I_d\! -\!\J   \J_{v_0\eta}^{-1}   \big)^2\big] \right\rangle_{\Zm} \right)
  + \frac{\beta  S^2}{N(\beta\!-\!\zeta)}\left\langle   \Tr\big[ \J \big( \I_d \!-\!   \J_{v_0\eta}^{-1}  \J\big)^2\big]  \right\rangle_{\Zm } \nonumber\\
 &&~~~~~~~~~~~~~~~~~~~~~~~~~~~~~~~~~~~~~~~~~~~-\frac{v_0^2\eta}{N(\beta\!-\!\zeta)}\left\langle\Tr\left[  \J^{-1}_{v_0\eta}  \right]  \right\rangle_{\Zm }  + \frac{2v_0^2\beta}{(\beta\!-\!\zeta)N}  \partial_{v}  \log P(v)\vert_{v=v_0}  
 \label{eq:nl-operator-aver-calc}
\end{eqnarray}
We then  make the usual substitution $\Zm\to\Zm/\sqrt{d}$. with $z_i(\mu)={\mathcal O}(1)$ for all $(i,\mu)$,  and we define the average  $\rho_d (\lambda )=\left\langle \rho_d (\lambda\vert\Zm ) \right\rangle_{\Zm}$ of the eigenvalue density $ \rho_d (\lambda \vert\Zm)=\frac{1}{d}\sum_{\ell=1}^d\delta\left(\lambda-\lambda_\ell\left({\Zm^{\Tran}} {\Zm}/N\right)\right)$ of the empirical covariance matrix. Then the above average  becomes
\begin{eqnarray}
 \left\langle \Psi \left[v_{0}\vert \Zm,  \thetav_{0},  \epsilonv\, \right] \right\rangle_{\Zm, \thetav_{0}, \epsilonv} 
 &=&     \frac{\beta  \sigma_0^2}{\beta\!-\!\zeta} \left( 1\!-\!\zeta\! +\!                             
  \zeta \int\!\rmd\lambda~   \rho_d (\lambda)           \Big(  \frac{\zeta v_0\eta }{\lambda +\zeta v_0\eta} \Big)^2 
  \right)
 + \frac{\beta  S^2}{\beta\!-\!\zeta}   \int\!\rmd\lambda~   \rho_d (\lambda)           \Big(  \frac{\zeta v_0\eta }{\lambda +\zeta v_0\eta} \Big)^2\!\lambda
 \nonumber\\
 &&~~~~~~~~~~~~~~~~~~~~~~~~~~~~~~~~~~~~~-\frac{v_0^2\eta\zeta^2}{\beta\!-\!\zeta}     \int\!\rmd\lambda~  \frac{  \rho_d (\lambda)  }{\lambda \!+\!\zeta v_0\eta}      + \frac{2v_0^2\beta}{(\beta\!-\!\zeta)N}  \partial_{v}  \log P(v)\vert_{v=v_0}  
 \label{eq:nl-operator-aver}
\end{eqnarray}
Finally we turn to the variance 
\begin{eqnarray}
\Var( \Psi \left[v_{0}\vert \Zm,  \thetav_{0},  \epsilonv\, \right])&=& 
\left(\frac{\beta}{\beta\!-\!\zeta}\right)^2 \!  \frac{1}{N^2}    \Var\left( \tv^\Tran  \left( \I_N \!-\!\Zm    \J_{v_0 \eta}^{-1}  \Zm^\Tran\right)^2\tv \right)      + \frac{v_0^4\eta^2}{(\beta\!-\!\zeta)^2}\frac{1}{N^2}\Var\left(\Tr\left[  \J^{-1}_{v_0\eta}  \right] \right)    \nonumber\\
     && ~~~~~~~~~~~~~~~~~~~~~~~~~~~~~~~~~~~ -\frac{2\beta v_0^2\eta}{(\beta\!-\!\zeta)^2}\frac{1}{N^2}       \Covar\left(   \tv^\Tran  \left( \I_N \!-\!\Zm    \J_{v_0 \eta}^{-1}  \Zm^\Tran\right)^2\tv,                    \Tr\left[  \J^{-1}_{v_0\eta}  \right]    \right) 
  \label{eq:nl-operator-var-calc-1}.
\end{eqnarray}
Computing in this expression the relevant averages over the random variables  $\Zm$,  $\thetav_{0}$ and  $\epsilonv$, with the familiar substitution $\Zm\to\Zm/\sqrt{d}$,  gives us the following result
\begin{eqnarray}
\Var( \Psi \left[v_{0}\vert \Zm,  \thetav_{0},  \epsilonv\, \right])&=& \!\int\!  \rmd\lambda\rmd\tilde{\lambda}~C_d(\lambda, \tilde{\lambda})\Bigg\{ \left(\frac{\beta}{\beta\!-\!\zeta}\right)^2   \left(  \frac{\zeta v_0\eta }{\lambda \!+\!\zeta v_0\eta} \right)^2 \!    \left(  \frac{\zeta v_0\eta }{\tilde{\lambda} \!+\!\zeta v_0\eta} \right)^2    \left[ \sigma^4_0\zeta^2  \! +\!S^4\lambda\tilde{\lambda} \!+\! 2\sigma_0^2S^2\zeta\lambda \right  ]  \nonumber\\
&&+ \left(\frac{v^2\zeta^2\eta}{\beta\!-\!\zeta}\right)^2      \frac{  1    }{  \left(\lambda \!+\!\zeta v_0\eta\right) \big( \tilde{\lambda} \!+\!\zeta v_0\eta\big)} -\frac{2\beta v_0^2 \zeta^2 \eta  }{(\beta\!-\!\zeta)^2}    \left(  \frac{\zeta v_0\eta }{\lambda \!+\!\zeta v_0\eta} \right)^2    \frac{\sigma_0^2\zeta  \!+\! S^2 \lambda}{\tilde{\lambda}\! +\!\zeta v_0\eta}  \Bigg\}
 \nonumber\\
&&+\frac{2}{N} \left(\frac{\beta}{\beta\!-\!\zeta}\right)^2   \int\!\rmd\lambda~  \rho_d (\lambda)  \Bigg\{  \sigma_0^4  \Big[  1\!-\!\zeta \!+\!\zeta   \Big(  \frac{\zeta v_0\eta }{\lambda \!+\!\zeta v_0\eta} \Big)^4 \Big]  +          \frac{S^4}{\zeta }   \Big(  \frac{\zeta v_0\eta }{\lambda \!+\!\zeta v_0\eta} \Big)^4\lambda^2 \Bigg\} 
  \label{eq:nl-operator-var},
\end{eqnarray}
 with the correlation function $C_d(\lambda, \tilde{\lambda})=\big\langle\rho_d (\lambda\vert\Zm )\rho_d (\tilde{\lambda}\vert\Zm )\!\big\rangle_{\Zm}-\big\langle\rho_d (\lambda\vert\Zm )\!\big\rangle_{\Zm}\! \big\langle \rho_d (\tilde{\lambda}\vert\Zm )\!\big\rangle_{\Zm}$. Clearly, if the spectrum $\rho_d (\lambda\vert\Zm )$ is self-averaging when $(N,d)\to\infty$, then the correlation function will vanish in this limit, and hence  $ \Psi \left[v_{0}\vert \Zm,  \thetav_{0},  \epsilonv\, \right]$ will be self-averaging.

\bibliographystyle{apsrev4-1}



%

\end{document}